\input amstex

\def\b1{\text{\bf 1}}

\def\Z{R_{\tr}}

\def\CA{{\Cal A}}
\def\CB{{\Cal B}}
\def\CC{{\Cal C}}
\def\CD{{\Cal D}}

\def\CE{{\Cal E}}
\def\CH{{\Cal H}}
\def\CI{{\Cal I}}

\def\CG{{\Cal G}}
\def\CL{{\Cal L}}
\def\CM{{\Cal M}}

\def\CO{{\Cal O}}
\def\CP{{\Cal P}}

\def\CS{{\Cal S}}
\def\CT{{\Cal T}}
\def\CV{{\Cal V}}

\def\CZ{{\Cal Z}}

\def\gr{\text{gr}}

\def\End{\text{End}}
\def\Hom{\text{Hom}}
\def\Sym{\text{Sym}}

\def\#{\,\check{}}

\def\codim{{\text{codim}}}

\def\id{\text{id}}
\def\tr{\text{tr}}

\def\Coker{\text{Coker}}

\def\Res{\text{Res}}
\def\Spec{\text{Spec}}

\def\limleft{\mathop{\vtop{\ialign{##\crcr
  \hfil\rm lim\hfil\crcr
  \noalign{\nointerlineskip}\leftarrowfill\crcr
  \noalign{\nointerlineskip}\crcr}}}}
\def\limright{\mathop{\vtop{\ialign{##\crcr
  \hfil\rm lim\hfil\crcr
  \noalign{\nointerlineskip}\rightarrowfill\crcr
  \noalign{\nointerlineskip}\crcr}}}}
\def\Ho{\text{Hot}}

\def\LA{\mathop{\vtop{\ialign{##\crcr
  \hfil\rm $\CA$\hfil\crcr
  \noalign{\nointerlineskip}\rightarrowfill\crcr
  \noalign{\nointerlineskip}\crcr}}}}

\def\LF{\mathop{\vtop{\ialign{##\crcr
  \hfil\rm $F$\hfil\crcr
  \noalign{\nointerlineskip}\rightarrowfill\crcr
  \noalign{\nointerlineskip}\crcr}}}}

\def\LI{\mathop{\vtop{\ialign{##\crcr
  \hfil\rm $\CI$\hfil\crcr
  \noalign{\nointerlineskip}\rightarrowfill\crcr
  \noalign{\nointerlineskip}\crcr}}}}

\def\LP{\mathop{\vtop{\ialign{##\crcr
  \hfil\rm $\CP$\hfil\crcr
  \noalign{\nointerlineskip}\rightarrowfill\crcr
  \noalign{\nointerlineskip}\crcr}}}}


\def\hra{\hookrightarrow}
\def\iso{\buildrel\sim\over\rightarrow}

\def\lra{\longrightarrow}

\def\limright{\mathop{\vtop{\ialign{##\crcr
  \hfil\rm lim\hfil\crcr
  \noalign{\nointerlineskip}\rightarrowfill\crcr
  \noalign{\nointerlineskip}\crcr}}}}

\parskip=6pt

\documentstyle{amsppt}
\document
\magnification=1100
\NoBlackBoxes

\bigskip

\bigskip
\centerline{\bf  A DG GUIDE TO VOEVODSKY'S MOTIVES}
\bigskip
\centerline{A.~Beilinson and V.~Vologodsky}
\medskip
\centerline{The University of Chicago}

\bigskip

\hskip 5.2cm {\it To Joseph Bernstein with love and gratitude.}

\bigskip

 Let $\CV ar$ be the category of complex algebraic varieties, $T\! op$ that of nice  topological spaces,  $D_{ab}$ be the derived category of  finite complexes of finitely generated abelian groups. One has tensor functors $\CV ar\to T\! op \to D_{ab}$, the first assigns to a variety its space equipped with the classical  topology, the second one is the singular chain complex functor
 (the tensor structure for the first two categories is given by the direct product). The basic objective of the motive theory is to fill in a commutative square 
$$\spreadmatrixlines{2\jot}
\matrix
\CV ar&
\to &
\CD_ \CM \\  
\downarrow&&\downarrow \\
T \! op &
    \to & D_{ab}
\endmatrix
\tag 0.0$$
where $D_\CM$ - the category of motives - is a rigid tensor triangulated category defined, together with the upper horizontal arrow, in a purely geometric way (so that the base field $\Bbb C$ can be replaced by any field), and the right  vertical arrow is a tensor triangulated functor (which absorbs all the transcendence of the singular chains).

The known constructions  (due to Hanamura, Levine, and Voevodsky) proceed by first embedding $\CV ar$ into a larger DG category, and then define $D_\CM$ as its appropriate quotient. Voevodsky's  generators and relations are especially  neat. 
A rough idea: 
consider  the localization of the DG category freely generated by the category of topological manifolds modulo the relations  that kill the complexes of types  $\Bbb Z [\Delta \times X] \to \Bbb Z [X]$ and $\Bbb Z [U\cap V ]\to \Bbb Z[U]\oplus \Bbb Z [V]\to \Bbb Z [X]$  (here $\Delta$ is an interval, $\{ U,V\}$ is an open covering of $X$); the singular chains functor yields then an equivalence between this toy category of ``topological motives" and $D_{ab}$. To define $\CD_\CM$, one formally imitates this construction in the algebro-geometric setting  with an important modification:  mere combinations of true algebraic maps  should be replaced from scratch by a larger group of finite correspondences, i.e., multi-valued maps (which is irrelevant in the topological setting).

What follows is a concise  exposition of Voevodsky's theory that  covers  principal points of   \cite{Vo1--3} and  \cite{MVW} with the notable exception of   
 comparison results (relating the motivic cohomology  with Bloch's higher Chow groups and  Milnor's K-groups, and the \'etale localized motives with finite coefficients with Galois modules, see \cite{MVW} 
19.1, 5.1,  \cite{Vo2} 3.3.3). 
The more advanced subjects of $\Bbb A^1$-homotopy theory, the proof of the Milnor-Bloch-Kato conjecture,  and  the array of  motivic dreams, are not touched. 

We take the time to spell out the basic constructions on the DG category level (\cite{Vo2} and \cite{MVW} consider mere triangulated category structure). For the present material,  this has the limited advantage 
of  making formulas like (4.4.1), (4.4.2) possible, but seems to be necessary
  for some future developments (such as understanding of the motivic descent). 
The required generalities are recalled in (the lengthy) \S 1;  we use concrete Keller's construction of homotopy DG quotients (see \cite{K} and \cite{Dr} \S 4) which suffices for immediate purposes. In truth, unfettered DG functoriality requires better understanding of the world of DG categories (e.g.~to be able 
to determine a DG category \`a la Yoneda by the homotopy DG 2-functor it represents). Hopefully,  the
necessary material can be found in Lurie's treatise \cite{L1}, \cite{L2}.

We are grateful to  Volodya Drinfeld and Andrey Suslin for teaching us, respectively, DG categories and motives. Our deep thanks due to  Volodya Drinfeld, Dennis Gaitsgory, Madhav Nori, Chuck Weibel, and the referee for their truly generous help in turning a raw draft into an article. The research was partially supported by NSF grant DMS-0401164.

\bigskip

\centerline{\bf   \S 1.  Homological algebra recollections}  

\medskip

The basic reference for  triangulated categories is \cite{Ve};  for the DG story it is  \cite{Dr}. The format of 1.3 helps to formulate Voevodsky's main technical result (see  4.4). 

{\bf 1.1.} Below $R$ is a fixed commutative ring, say, $R=\Bbb Z$. We  play with $R$-categories.
 Let $A$ be one. For a full subcategory $I\subset A$ its {\it right orthogonal complement} is $I^\bot := \{ M\in A : \Hom (N,M )=0\,\, \forall N\in I\}$; the left orthogonal complement ${}^\bot I$ is defined in the dual way.   We say that $A$ is {\it cocomplete} if it is closed under  direct sums of arbitrary cardinality. For such an $A$,   an object $M\in A$ is said to be {\it compact} if for every set of objects $\{ N_\alpha \}$ one has $\oplus\,\Hom (M,N_\alpha )\iso \Hom (M,\oplus\, N_\alpha )$.  

We denote by  $A^{\text{op}}$ the dual category; $M\mapsto M^{\text{op}}$ is  the contravariant ``identity" functor $A\to A^{\text{op}}$. $A^\kappa$ is {\it the idempotent completion} (a.k.a.~{\it Karoubianization}) of $A$.
 For $A$  essentially small, $A$-mod is the category of $A$-modules, i.e., $R$-linear functors $A\to R$-mod (the category of $R$-modules); this is a cocomplete  abelian $R$-category. 

 If $D$ is a triangulated category, then such is $D^\kappa$ (see \cite{BS}).

{\it Remarks.} (i) According to  \cite{Th}, arbitrary  triangulated categories can be described in terms of Karoubian ones as follows. Let $D$ be a small Karoubian triangulated category. For a subgroup $H$ of the Grothendieck group $K(D)$ let $D_H \subset D$ be the subcategory of objects whose classes lie in $H$. Then  the map $H\mapsto D_H$ is a bijection between the set of subgroups of $K(D)$ and that of strictly full triangulated subcategories $C \subset D$ such that ${C}^{\kappa}=D$.

(ii) In many situations it is easy to check if a given triangulated category is Karoubian. E.g., in \cite{BS} it is shown that the bounded derived category of a Karoubian exact category is Karoubian; for other examples, see 1.4.1 and 1.5.6.

{\bf 1.2. Admissible subcategories.} Let $D$ be a triangulated ($R$-)category,  $I\subset D$ be a strictly full triangulated subcategory.  $I$ is said to be {\it thick} if any direct summand of an object of $I$ lies in $I$.

For $M\in D$ a {\it right $I$-localization triangle} for $M$ is an exact triangle $M_I \buildrel{\mu}\over\to M\buildrel{\nu}\over\to M_{I^\bot}$ with $M_I \in I$ and $M_{I^\bot}\in I^\bot$. If such a triangle exists, then it is unique (up to a uniquely defined isomorphism). If it exists for every $M\in D$, then  $I$  is  said to be {\it right admissible}. This  amounts to either of the following properties:

(i) The embedding functor $i_* : I \to D$ admits a right adjoint $i^! : D \to I$.

(ii) $I$ is a thick subcategory of $D$, the Verdier quotient $D /I$ is well defined, and the projection functor $j^* : D \to D/I $ admits a right adjoint $j_* : D/I \to D$.

If this happens, then  $M_I = i_* i^! M$ and $M_{I^\bot}= j_* j^* M$. Also   $I^\bot $ is a full triangulated subcategory of $D$, $I \iso {}^\bot (I^\bot )$, the projection $I^\bot \to D /I$ is an equivalence of categories, and $j_*$ is the composition of the inverse equivalence and the embedding $I^\bot \hra D$. 
The endofunctor $\CC_I = j_* j^*$ of $D$, $\CC_I (M)= M_{I^\bot}$, equipped with the morphism of functors $\nu =\nu_I :$ Id$_D \to \CC_I$, is called  {\it  the (right) $I$-localization}. Notice that the pair $(\CC_I ,\nu_I )$ admits no non-trivial automorphisms

{\it Example.} Let $X$ be a topological space, $i : Y\hra X$ its closed subspace, $j: U \hra X$ the complement to $Y$. Let $D$ be the derived category of sheaves of abelian groups on  $X$, and $I$ that on $Y$. Then $i_* : I\to D$ is a fully faithful functor whose essential image is right admissible, and $j^*$ identifies $D/I$ with the derived category of sheaves on $U$. The same is true in the $\CO$- and $\CD$-module setting of algebraic geometry.

There is a dual notion of left-admissible subcategory; $I \leftrightarrow I^\bot$ is a 1-1 correspondence  between the ``sets" of  right- and and left-admissible subcategories of $D$. The duality $D\mapsto D^{\text{op}}$ interchanges left- and right-admissible subcategories. The functors $(j_* , j^* , i_* ,i^! )$ for $I^{\bot \text{op}}\subset D^{\text{op}}$ are the same as functors $(i_* ,i^! , j_* ,j^* )$ for $I\subset D$.

Suppose  we have a t-structure on $D$ (see \cite{BBD} 1.3), and let
$I\subset D$ be a right-admissible subcategory. The next lemma is due to D.~Gaitsgory:

\proclaim{\quad Lemma} (i)  $\CC_I $ is left exact if and only if  $D/I$ admits a t-structure such that $j^*$ is exact. The latter t-structure is unique.

(ii) $\CC_I $ is right exact if and only if $D/I$ admits a t-structure such that $j_*$ is t-exact. The latter t-structure is unique.
  \hfill$\square$
\endproclaim

{\it Exercise.} 
Condition (i) for $I\subset D$ implies (ii) for $I^{\bot\text{op}}\subset D^{\text{op}}$ (and the opposite t-structure).

{\bf 1.3. Compatibility.} Two right-admissible subcategories $I_1 ,I_2$ of  $ D$ are said to be {\it compatible} if every  $M\in D$ can be fitted into a commutative diagram
$$
\spreadmatrixlines{3\jot}
\matrix
* &\to &   * & \to & * \\
\uparrow &&\uparrow &&\uparrow \\
*  &\to & M & \to & *\\
\uparrow &&\uparrow&&\uparrow \\
* &\to &  * &\to & *
   \endmatrix
\tag 1.3.1
$$
with exact rows and columns such that the bottom row lies in $\CI_1$, the top row in $I_1^\bot$, the left column in $I_2$, and the right column in $I_2^\bot$. In other words, the columns are (right) $I_1$-localization triangles, and the rows are $I_2$-localization ones.

\proclaim{\quad Lemma}  Compatibility of $I_1$, $I_2$ amounts to either of the next properties:

(i) For every  $M\in D$ there exist morphisms $M_0 \to M_1 \to M$ such that $M_0 \in I_1 \cap I_2$, $\CC one (M_0 \to M_1) \in ( I_1 \cap I_2^\bot ) \times (I_1^\bot \cap I_2 ) \subset D$, $\CC one (M_1 \to M)\in I_1^\bot \cap I_2^\bot$.

(ii) The $I_2$-localization  $\CC_2 = j_{2*} j_2^*$ preserves both $I_1$ and $I_1^\bot$.

(ii)$' $ The endofunctor $i_{2*} i_2^!$ preserves both $I_1$ and $I_1^\bot$.

(iii) $\CC_2$ preserves $I_1$, and $\CC_1$ preserves $I_2$.

(iii)$' $  $i_{2*} i_2^!$ preserves $I_1$, and $i_{1*} i_1^!$ preserves $I_2$.

(iv) The endofunctors  $\CC_1 \CC_2$ and $\CC_2 \CC_1$ are isomorphic.
 \hfill$\square$
\endproclaim

Let $I_{12}$ be the triangulated subcategory of $D$  strongly generated by $I_1$ and $I_2$ (see 1.4.2), so  $I_{12}^\bot = I_1^\bot \cap I_2^\bot$. If $I_1$, $I_2$ are compatible, then $I_{12}$ is right admissible,   the $I_{12}$-localization functor $\CC_{12}$ is canonically isomorphic to  $\CC_1 \CC_2$ and $\CC_2 \CC_1$,\footnote{ Precisely, $(\CC_{12},\nu_{12})= (\CC_1 \CC_2 , \nu_1 \nu_2 )=(\CC_2 \CC_1 ,\nu_2 \nu_1 )$.} and the projections  $I_1^\bot /I_1^\bot \cap I_2 \leftarrow I_1^\bot \cap I_2^\bot
\to I_2^\bot /I_2^\bot \cap I_1$ are  equivalences.

{\it Examples.} (i) In the  example of 1.2, the categories $I=I_Y$ for various closed $Y$'s are pairwise compatible.  (ii) If $D$ is the derived category of $\CD$-modules on the line $\Bbb A^1$, then the subcategory of $\CD$-modules supported at a point and its  Fourier transform  are {\it not} compatible.

{\bf 1.4. Cocomplete vs.~small $\&$ Karoubian.} The results of 1.4.1, 1.4.2 are due to A.~Neeman \cite{N}, those of 1.4.3 to  V.~Drinfeld.

1.4.1. Let $D$ be a triangulated category.   For a diagram $M_0 \buildrel{i_0}\over\to M_1 \buildrel{i_1}\over\to M_2   \buildrel{i_2}\over\to  \ldots$ in $D$ we set $hocolim(M_a ,i_a ):= \CC one (\oplus M_a \buildrel\nu\over\to \oplus M_a )$ where $\nu (m_a ):= m_a  - i_a (m_a )$. Notice that, as  cones do, the object  $hocolim(M_a ,i_a )$ is defined up to a {\it non-canonical} isomorphism. There are canonical morphisms $\mu_a : M_a \to hocolim (M_a ,i_a )$ such that $\mu_{a+1}i_a =\mu_a$; for every collection of morphisms $\phi_a  :M_a \to N$ such that $\phi_{a+1}i_a = \phi_a$ there exists a morphism $\phi : 
hocolim (M_a ,i_a )\to N$ such that $\phi \mu_a =\phi_a $.

\proclaim{\quad Lemma} 
If $D$ is closed under countable direct sums, then $D$ is Karoubian.
\endproclaim

{\it Proof.} 
The image of an idempotent $\pi \in \End M$  is $M_\pi := hocolim (M_a ,i_a )$, where $M_a \equiv M$, $i_a \equiv \pi$, equipped with the structure maps $ M_\pi \to M$ induced from $\phi_a =\pi$ and $ M\to M_\pi$ equal to $\mu_0$.  \hfill$\square$

1.4.2. For a class $S$ of objects of $D$ {\it the triangulated subcategory strongly  generated by $S$} is
 the smallest strictly full triangulated subcategory  $\langle S\rangle $ that contains $S$. If $D$ is cocomplete, then  {\it the triangulated subcategory   generated by $S$} is the smallest strictly full  triangulated  subcategory of $ D$ that contains $S$ and is  closed under arbitrary direct sums.

 If $I\subset D$ is a thick subcategory, then we set $D/^\kappa I := (D/I)^\kappa = (D^\kappa /I^\kappa )^\kappa$.

If $D$ is a cocomplete triangulated category, then $D^{\text{perf}}\subset D$ denotes the full subcategory of compact, a.k.a.~{\it perfect}, objects. This is a  Karoubian subcategory.

\proclaim{\quad Proposition} 
Suppose $D$ is cocomplete and let $S\subset D$ be a set of perfect objects. Set $I^o := \langle S\rangle$, and let  $I$ be the triangulated subcategory   generated by $S$.  

(i)  $I$  is right admissible, and $I^{per\!f}=(I^{o})^{\kappa}$. 

(ii) The category $D/I$ is cocomplete, the adjoint functors $(j^* , j_* ): D \leftrightarrows D/I$, $(i_* ,i^! ): I \leftrightarrows D$ commute with arbitrary direct sums, and $j^* (D^{\text{perf}})\subset (D/I)^{\text{perf}}$.  

(iii)  If $D$ is   generated by a set of perfect objects, then $D^{\text{perf}}/^\kappa I^{\text{perf}}  \iso (D/I)^{\text{perf}}$.

(iv) Every object of $I$ can be represented as $hocolim (M_a ,i_a )$ (see 1.4.1) where $M_0$ and each $\CC one (i_a )$ are direct sums of translations of objects from $S$.
\endproclaim

{\it Proof.}  We can assume that $S$ is closed with respect to translations: $S[\pm 1]=S$. Then $S^\bot =I^\bot$.

(a) For $M\in D$ let us construct  a right $I$-localization triangle with $M_I$ as in (iv). 

 Suppose we have a diagram  $(M_a ,i_a )$, as in 1.4.1, in  $I$ and morphisms $\phi_a : M_a \to M$ such that $\phi_{a+1} i_{a}=\phi_a$ with
 the next property:  If $P\in S$,  then every morphism $\xi : P\to M$ factors through $\phi_0$, and for every  $\psi: P\to M_a $ such that $\phi_a \psi =0$ one has $i_a \psi =0$. Set $M_I :=
 hocolim (M_a , i_a )$; let $\phi : M_I \to M$
 be a morphism such that  $\phi \mu_a =\phi_a $.   Then $M_{I^\bot}:=\CC one (\phi )\in S^\bot = I^\bot$, so $M_I \buildrel{\phi}\over\to M\to M_{I^\bot}$ is a right $I$-localization triangle.
 
It remains to construct such a  $(M_a , i_a ,\phi_a )$ with $(M_a ,i_a )$ as in (iv).  We do this by induction.
 Choose $\phi_0 : M_0 \to M$ so that $M_0$ is a direct sum of  objects in $S$ and every $\xi$ as above factor through $\phi_0$.
Given $\phi_a: M_a \to M$, we define $M_{a+1}$ as $\CC one (\eta_{a} : C_{a} \to M_a )$ where $C_{a}$ is a direct sum of  objects in $S$, $\phi_a \eta_{a} =0$, and  every $\psi$ as above factors through $\eta_{a}$. Then $i_a : M_a \to M_{a+1}$ is the evident morphism, and $\phi_{a+1}$ is any morphism $M_{a+1}\to M$ such that $\phi_{a+1}i_a =\phi_a$. 

(b)   Let us show that  for $M\in I$ every morphism $\theta : N\to M$ with $N\in D^{\text{perf}}\cup I^{\text{perf}} $ factors through an object of $I^o$.

Let $J\subset I$ be the subcategory of $M$'s that satisfy this property. It contains $I^o$ and is closed under translations and arbitrary direct sums. It remains to check that $J$ is closed under cones. Let $\Delta_M = (M\to M' \to M'' )$ be an exact triangle with $M' ,M'' \in J$; we want to show that $M\in J$. The composition $N\buildrel\theta\over\to M\to M'$ factors through some $N'\in I^o$, so there is an exact triangle $\Delta_N = (N\to N' \to N'')$ and a morphism of triangles $(\theta, \alpha ,\beta ): \Delta_N \to \Delta_M$. Since $N'' \in D^{\text{perf}}\cup I^{\text{perf}} $, one can write $\beta$ as composition $N'' \buildrel\gamma\over\to L \buildrel\delta\over\to M''$ where $L\in I^o$. There is an exact triangle $\Delta_L = (K\to N' \to L )$ and morphisms of triangles $\Delta_N \buildrel{(\epsilon  ,\id_{N'} , \gamma )}\over\lra \Delta_L \buildrel{(\xi , \alpha ,\delta )}\over\lra \Delta_M$. The composition $N'' [-1] \to N \buildrel{\theta - \xi\epsilon}\over\lra M$ (the first arrow comes from $\Delta_N$) vanishes, so $\theta -\xi\epsilon$ factors through $N'\in I^o$. Thus $\theta$ factors through  
 $N' \oplus K\in I^o$, q.e.d.

(c) Applying (b) to $\id_N$, $N\in I^{\text{perf}}$, we see that $I^{\text{perf}}=(I^{o})^{\kappa}$;
we have proved (i). Since the direct sum in $D$ of an arbitrary family of objects of $I^\bot$ lies in $ S^\bot = I^\bot$, we get (ii). Let us prove (iii). It follows from (b) that the functor $D^{\text{perf}}/I^{\text{perf}} \to (D/I)^{\text{perf}}\subset D/I$ is  fully faithful. By (i) (with $I$ replaced by $D$), $D^{\text{perf}}$ is essentially small. Its image  generates $D/I$, so by (i) (with $I$ replaced by $D/I$,  $S$ by $j^* (D^{\text{perf}})$),  $(D^{\text{perf}}/I^{\text{perf}} )^\kappa = \langle j^* (D^{\text{perf}})\rangle^\kappa = (D/I)^{\text{perf}}$, q.e.d.
 \hfill$\square$

1.4.3. Suppose $D$ is a cocomplete triangulated category and $I_1^o ,I_2^o \subset D^{\text{perf}}$ are two {\it small} full  triangulated subcategories. Let $I_1 , I_2 \subset D$ be the triangulated subcategories  generated by $I_1^o ,I_2^o $. By 1.4.2, $I_1 , I_2$ are right admissible.

\proclaim{\quad Proposition}  $I_1 , I_2$ are compatible if and only if  the following condition holds:  Every morphism $\phi :M_1 \to M_2$, where $M_i \in I_i^o$, can be represented as a composition $M_1 \to N_2 \to N_1 \to M_2$ with $N_i \in I_i^o$, and the same is true if we  interchange $I_i^o$.
\endproclaim

{\it Proof.} (i) Suppose $I_1$, $I_2$ are compatible; let us check the condition. Write  $\phi$ as the composition $M_1 \buildrel{\psi}\over\to P\to M_2$ where $P:= i_{1*} i^!_1 M_2 $. Since $P\in I_2$ and $M_1 \in D^{\text{perf}}$, by part (b) of the proof in 1.4.2, we can write $\psi$ as a composition $M_1 \to N_2 \buildrel{\chi}\over\to P$ with $N_2 \in I^o_2$. Similarly, since $P\in I_1$ and $N_2 \in D^{\text{perf}}$, $\chi$ can be written as a composition $N_2 \to N_1 \to P$ with $N_1 \in I_1^o$, and we are done.

(ii) Suppose the condition holds; let us show that $I_1$, $I_2$ are compatible. By 1.3(iii)$'$, it suffices to check that for any $F \in I_2$ the object $P:=i_{1*} i^!_1 F$ lies in $I_2$. We will prove  that 
any morphism $\theta : M\to P$ with $M\in D^{\text{perf}}$ factors through an object of $I_2$. This implies the claim. Indeed, 
replacing $D$ by the subcategory generated by $I_1$, $I_2$, we can assume that $D$ is  generated by $D^{\text{perf}}$. By 1.4.2(ii), $P\in I_2$ iff for any $M' \in D^{\text{perf}}$ one has $\Hom_{D/I_2}(M', P )=0$. Any morphism $M' \to P $ in $D/I_2$ can be lifted to a morphism $\theta : M \to P$ in $D$, where $M := \CC one (\xi )[-1]$ for some $\xi : M'\to N$, $N\in I_2$; by part (b) of the proof in 1.4.2, we can assume that $N\in I_2^o$, so $M \in D^{\text{perf}}$, and $\theta$ vanishes in $D/I_2$ by the assertion, q.e.d. 

 By part (b) of the proof in 1.4.2,  our $\theta : M \to   P$ factors through some $M_1 \in I_1^o$, and  the composition $M_1 \buildrel{\theta}\over\to  P \to F$ factors through some $M_2 \in I_2^o$. Let us factorize the map $M_1 \to M_2$ as in the condition. By the universal property, the composition $N_1 \to M_2 \to F$ factors through $ P $, and  $M_1 \to P $ equals the composition  $M_1 \to  N_2 \to N_1 \to P $. We have factorized $\theta$ as  $M\to  N_2 \to P $, q.e.d.
\hfill$\square$

{\bf 1.5. DG categories.}   With respect to triangulated categories, DG categories  play the same role as topological spaces for objects of the homotopy category.

 1.5.1. Let $\CA$ be a DG ($R$-)category. Thus morphisms between objects of $\CA$ are complexes of $R$-modules, and the composition maps are morphisms of complexes.  The {\it homotopy category} $\Ho \, \CA$ and the {\it graded homotopy category} $\Ho^\cdot \CA$ have same objects as $\CA$, and  morphisms $\Ho (M,M'):= \Hom_{\Ho\, \CA}(M,M')=H^0 \Hom (M,M')$, 
 $\Ho^\cdot (M,M'):= \Hom_{\Ho^\cdot \CA}(M,M')=\oplus H^n \Hom (M,M')$. 
 A {\it DG morphism} $\phi$ in $\CA$ is a closed morphism  of degree 0. Such $\phi$ yields a morphism in $\Ho\, \CA$; if the latter is an isomorphism, then $\phi$ is said to be 
{\it homotopy equivalence}. 

Our  $\CA$ is  {\it DG Karoubian} if it contains the image of any idempotent DG endomorphism. For an arbitrary $\CA$ we denote by $\CA^\kappa$ its DG idempotent completion (defined in the evident way), which is a DG Karoubian DG category.

If $\CA$,  $\CB$ are DG categories, then a DG ($R$-linear) functor $F: \CA \to \CB$ yields  functors $\Ho\, F : \Ho\, \CA \to \Ho\, \CB$,  $\Ho^\cdot F : \Ho^\cdot \CA \to \Ho^\cdot \CB$. Our $F$ is a   {\it quasi-equivalence} if $\Ho^\cdot F$ is a fully faithful functor and $\Ho\, F$ is essentially surjective. 

{\it Examples.}  (i) An associative unital DG ($R$-)algebra $A$ is the same as a DG category with single object; we denote the latter also by $A$.

(ii) For any $R$-category $\CA$ the category $C\CA$ of complexes in $\CA$ is a DG category.

1.5.2. A {\it homotopy DG functor}  $\tilde{F}=(\tilde{F},\iota ):\CA\to\CB$ is a diagram $\CA \buildrel{\iota}\over\leftarrow \tilde{\CA} \buildrel{\tilde{F}}\over\lra \CB$ of DG functors such that $\iota$ is a quasi-equivalence. If $\CA =R$, then this referred to as {\it homotopy object} of $\CB$ (ordinary objects are the same as plain DG functors $R\to\CB$).
 
 {\it Remark.} $\tilde{F}$ yields a well-defined functor $\Ho\,\tilde{F}: \Ho\,\CA \to\Ho\,\CB$. E.g.~a homotopy object of $\CB$ yields a well-defined object of $\Ho\, \CB$. The converse is {\it not} true.
 
 Homotopy DG functors can be composed: the composition of $\CA \buildrel{\iota}\over\leftarrow \tilde{\CA} \buildrel{\tilde{F}}\over\lra \CB$ and $\CB \buildrel{\pi'}\over\leftarrow \tilde{\CB} \buildrel{\tilde{G}}\over\lra \CC$ is $\CA\leftarrow \tilde{\CA} \buildrel{\tilde{H}}\over\lra \CC$, where $\tilde{\CA}'$ is the ``homotopy fibered product" of $\tilde{\CA}$ and $\tilde{\CB}$ over $\CB$. Namely, objects of $\tilde{\CA}'$ are triples $(P,Q, \ell )$ where $P\in\tilde{\CA}$, $Q\in\tilde{\CB}$, $\ell :\tilde{F}(P)\to \pi' (Q)$ is a homotopy equivalence; $\Hom^n ((P,Q, \ell ),(P',Q', \ell' ))$ consists of triples $(\phi ,\psi ,\tau )$, $\phi \in\Hom^n (P,P')$, $\psi\in\Hom^n (Q,Q' )$, $\tau\in\Hom^{n-1}(\tilde{F}(P),\pi' (Q' ))$;  the differential on $\Hom$'s is determined by the property that natural  functors $\tilde{\CA}' \to \tilde{\CA},\tilde{\CB},\CB$,   $(P,Q, \ell )\mapsto P, Q,\CC one (\ell )$, are DG functors.

 1.5.3. Two DG categories $\CA$, $\CB$ yield DG categories $\CA\otimes \CB = \CA\mathop\otimes\limits_R \CB$ and (if $\CA$ is essentially small) $\CH om (\CA,\CB)$. Namely,   objects of  $\CA\otimes \CB$ are written as $M\otimes N$ where $M\in \CA$, $N\in \CB$, and $Hom (M\otimes N, M'\otimes N')= Hom (M,M') \mathop\otimes\limits_R Hom (N,N' )$.
Objects of $\CH om(\CA,\CB)$ are  DG functors $F: \CA\to \CB$;  the complex $Hom (F,F')$ has components $Hom^n (F,F')$ formed by degree $n$ morphisms of graded functors, the differential is  evident. There is an ``evaluation" DG functor $\CH om (\CA,\CB)\otimes \CA \to \CB$, $F\otimes M \mapsto F(M)$, which satisfies
a universality property: if $\CC$ is another essentially small DG category, then  $\CH om (\CC,\CH om (\CA,\CB))=\CH om (\CC\otimes \CA, \CB)$. We get:

\proclaim{\quad Lemma} Small DG categories form a tensor category $\CD \CG$ with inner Hom's.  \hfill$\square$
\endproclaim

{\it Remark.} While $\otimes$ preserves quasi-equivalences (when one of the arguments is homotopically $R$-flat), this is {\it not} true for $\CH om$, even if $\CA$ is semi-free in the sense of \cite{Dr} B4, for $\CH om (\CA ,\CB)$ has too few morphisms. As shown in \cite{T}, replacing strict morphisms between DG functors by lax ones, one gets  a natural 2-(DG)category structure  compatible with the homotopy localization.

1.5.4. Let $\CA$ be a DG category. For $M\in \CA$ its {\it $n$-translation}, $n\in\Bbb Z$, is an object $M[n]$ of $\CA$ representing the DG functor $N\mapsto \Hom (N,M)[n]$. Similarly, for a DG morphism $\phi : M\to M'$ its {\it cone} is an object $\CC one (\phi )$ of $\CA$ that represents the DG functor $N\mapsto \CC one (\Hom (N,M)\buildrel\phi\over\to \Hom (N,M'))$.
Our $\CA$ is  {\it strongly pretriangulated} if $\CA$ contains 0 and $M[n]$, $\CC one (\phi )$ exist for every $M$, $n$, and $\phi$.

For every  $\CA$ there is a fully faithful DG embedding $\CA\hra \CA^{\text{pretr}}$ such that $\CA^{\text{pretr}}$ is strongly pretriangulated and the smallest strongly pretriangulated DG subcategory that contains $\CA$ equals $\CA^{\text{pretr}}$. Such an embedding is essentially unique: more precisely, if $\CB$ is strongly pretriangulated, then every DG functor $F: \CA \to \CB$ extends in a unique, up to a unique DG isomorphism, way to a DG functor $F^{\text{pretr}} :\CA^{\text{pretr}}\to\CB$. For a concrete construction of $\CA^{\text{pretr}}$, see Remark (i) in 1.6.1.

A DG category $\CA$ is said to be {\it pretriangulated} if $\Ho\, \CA\to \Ho\,  \CA^{\text{pretr}}$ is an equivalence of categories.  If this happens, then $\Ho\, \CA$  is naturally a triangulated category\footnote{To see this, notice that the Yoneda embedding  $\CA\hra \CA^{\text{op}}$-dgm (see 1.6.3) identifies $\Ho\,\CA$ with the full subcategory of the homotopy category of right DG $\CA$-modules, which is triangulated by the standard argument. If $\CA$ is pretriangulated, then $\Ho\, \CA$ is  closed under cones and translations.}  (see \cite{BK}).
For an arbitrary $\CA$, we denote by $\CA^{\text{tri}} $ the triangulated category $ \Ho\, \CA^{\text{pretr}}$. Any DG functor $F:\CA \to\CB$ yields a triangulated functor $F=F^{\text{tri}}: \CA^{\text{tri}}\to\CB^{\text{tri}}$.

\proclaim{\quad Lemma} If $F$ is a quasi-equivalence, then such is $F^{\text{pretr}}$. \hfill$\square$
\endproclaim

{\it Example.} If $\CA$ is an additive category considered as a DG category, then $\CA^{\text{pretr}}$ is  the DG category $C^b \CA$ of bounded complexes in $\CA$, and $\CA^{\text{tri}}$ is the bounded homotopy category $K^b \CA$.

1.5.5. For a triangulated category $D$, a {\it DG structure} on $D$ is a pretriangulated DG category $\CA =\CA_D$ together with an equivalence of triangulated categories $ \CA^{\text{tri}} \iso D$. It is clear what a  {\it DG lifting} of a triangulated functor between triangulated categories equipped with DG structures is. 
A DG structure on $D$  automatically induces a DG structure on any (strictly) full triangulated subcategory $I \subset D$: take for $\CA_I$ the preimage of $I$ in $\CA$. 
If $I$ is right admissible, then $\CA_{I^\bot}$ is a DG structure on $D /I$ called the {\it   (right) DG quotient DG structure}\footnote{The construction is {\it not} self-dual: it is adapted for handling the right derived functors.}  and denoted by $\CA/\CA_{I}$. 

{\it Remark.}  The projection $D \to D/I$ lifts naturally to a homotopy DG functor $\CA \leftarrow \tilde{\CA} \to \CA /\CA_I$. Namely, an object of $\tilde{\CA}$ is a DG morphism $\phi : M\to N$ in $\CA$ such that $N\in \CA_{ I^\bot}$, $\CC one (\phi )\in \CA_I$, and  the complex $\Hom_{\tilde{\CA}}(\phi ,\phi' )$ is the subcomplex of $\Hom (\CC one (\phi ), \CC one (\phi' ))$ formed by those morphisms that map $N\subset \CC one (\phi )$ to $N' \subset \CC one (\phi' )$; the above DG functors send $\phi$ to  $M$ and $ N$.

 1.5.6. The proposition below is inspired by \S 6 of \cite{Bo1}.\footnote{For more general results in this vein see \cite{Bo2}. We are grateful to  M.~Bondarko for  sending us a draft of the article and relevant correspondence.}
  
 Let $D$ be a triangulated category and $B\subset D$ be a  full  subcategory closed under finite direct sums which  strongly generates $D$.  Suppose that for every $M,N\in B$ and $n>0$ one has $\Hom (M, N[n])=0$. Let $\CA^{\text{tri}}\iso D$ be a DG structure on $D$.
  
 \proclaim{\quad Proposition} (i) There is a natural homotopy DG functor $\epsilon :\CA \to C^b B$. 
 
 (ii) The functor $\epsilon^{\text{tri}}: D\to K^b B$ is conservative. For $D$ essentially small, it yields an isomorphism of the Grothendieck groups $K_0 (D)\iso K_0 (K^b B)=K_0 (B)$.
 
 (iii) If $B$ is Karoubian, then  $D$ is Karoubian. \endproclaim
 
 {\it Proof.} (i) Let $\CB $ be a DG subcategory of $\CA$ formed by those $M\in \CA$ that are isomorphic to objects of $B$ in $D$, and $\Hom_{\CB} (M,M'):=\tau_{\le 0}\Hom_{\CA}(M,M' )$. 
 Then $\CB^{\text{tri}}\iso\CA^{\text{tri}}\iso D$. So we can assume that $\CA =\CB^{\text{pretr}}$.
 
One has a DG functor $\epsilon_0 : \CB \to B$ which is the identity on objects and is the projection $\Hom_{\CB}(M,M' )=  \tau_{\le 0}\Hom_{\CB}(M,M' )   \to   H^0\Hom_{\CB}(M,M' )$ on morphisms. Set $\epsilon := \epsilon_0^{\text{pretr}}: \CB^{\text{pretr}}\to B^{\text{pretr}}=C^b B$.

(ii) By  Remark (i) in 1.6.1, every  $M\in\CA=\CB^{\text{pretr}}$ can be written as $(\oplus M^i [-i] ,\phi )$ where $M^i \in \CB$. The components $\phi_i^j \in \Hom^{i-j+1}(M^i ,M^j )$ vanish for $i\ge j$, so one has a finite decreasing filtration  $F^a M :=  \mathop\oplus_{i\ge a} M^i [-i]$ on $M$. Every  morphism $\theta : M\to M'$ of degree 0 is compatible with  $F$; we write $\theta^i := \gr^i_F \theta \in\Hom (M^i ,M^{\prime i})$. 
The complex $\epsilon (M)$ has components $M^i \in B$, its differential $d \in \Hom_B (M^i ,M^{i+1})$ is the class of $\phi^{i+1}_i $ in $H^0 \Hom (M^i ,M^{i+1})$.

$\epsilon^{\text{tri}} $ is conservative: 
Suppose $\epsilon (M )$ is contractible. This means that there are decompositions $M^i = L^i \oplus L^{i+1}$ in $B$ such that $d_i$ equals the composition of the projection and the embedding $M^i \twoheadrightarrow L^{i+1}\hra M^{i+1}$. A downward induction by $a$ shows that the composition
$F^a M\twoheadrightarrow M^a [-a] \twoheadrightarrow L^a [-a]$ is a homotopy equivalence.
For $a\ll 0$ we get $0=F^a M =M$, q.e.d.

$K_0 (\epsilon^{\text{tri}} )$ is an isomorphism: Since $B$ strongly generates $D$, $K_0 (D)$ is generated by classes $[M]$, $M\in B$. The relation $[M\oplus M' ]= [M]+[M']$ evidently holds in  $K_0 (D)$. Since $K_0 (B)$ is universal with respect to these generators and relations, we get a map  
$K_0 (K^b B)=K_0 (B)\to K_0 (D)$ which is evidently inverse to $K_0 (\epsilon^{\text{tri}} )$.

 (iii) Suppose $B$ is Karoubian; let us show that every $P\in D^\kappa$ lies in $D$.

Choose $M\in \CA$ and an idempotent $p\in H^0 \End M$ whose image equals $P$. Write $M=(\oplus M^i [-i],\phi )$ as above.
Shifting $M$, we can assume that $F^n M=0$, $F^1 M=M$ for some odd $n\ge 1$. 
By \cite{Th} (see Remark (i) in 1.1.1), $P[n]\oplus P$ lies in $D$ for every odd $n$.\footnote{A direct proof: If $Q\oplus P\in D$, then $ P[1]\oplus P =\CC one ( Q\oplus P \to Q\oplus P)\in D $, then $ Q[2]\oplus P =\CC one ((Q\oplus P)[1]\to (P[1]\oplus P))\in D$, then $P[3]\oplus P=\CC one ((Q\oplus P)[2]\to  (Q[2]\oplus P) )\in D$, etc. } Choose $N=  (\oplus N^i [-i ],\psi )   \in \CA$ that represents $P[n]\oplus P$. The projectors $\bar{p}, \bar{p}'\in H^0  \End\, N$  onto $P$, $P[n]$ factor through, respectively, $M$ and $M[n]$, so one can represent them by DG endomorphisms $p ,p'$ of $N$ such that $p$ factors through $M$ and $p'$ factors through $M [n]$.

 Consider the complex  $\epsilon (N)\in C^b B$. Its identity map  is homotopic to $\epsilon ( p )+ \epsilon ( p' )$; let $h \in \Hom^0 (\epsilon (N) ,\epsilon (N)[-1])$, be the homotopy. Set $q := h^1 d,\, q':= dh^0 \in H^0 \End (N^0 ) $. Since $ M^0 =0= (M[n])^0$, one has $ p^0 =0= p^{\prime 0}=0$. Therefore
  $\id_{ N^0} = q+ q' $ and $q^2 = h^1 d (h^1 d + dh^0 )=q$, $q^{\prime 2} = (h^1 d + dh^0 )dh^0 =q'$, $qq'=0$. The images of the idempotents $q,q'$ lie in $B$ since $B$ is Karoubian; let us represent them by $Q$, $Q'\in\CB$, so  we have an homotopy equivalence $Q\oplus Q' \iso  N^0$.
 
 Let $T$ be the cone of the composition $Q[-1] \to  N^0 [-1] \buildrel{\psi_0^{\cdot}}\over\lra  F^1 N$, so $F^0 T=T, F^1 T=F^1 N$. There is an evident DG  morphism $e : T \to N$ identical on $F^1$. We will show that  $e$  identifies $T\in D$ with $P$. Since $F^1 M =M$,  there is a DG morphism $t : N \to T$ such that $e t =p$. It remains to check that $te\in\End T$ is a homotopically invertible. Now $\epsilon (e )$ identifies $\epsilon (T )$ with a subcomplex of $\epsilon (N )$, and $h^{>0}$ is a homotopy between $\epsilon (t e )$ and $\id_{\epsilon (T )}$. Thus $\epsilon (te)=\id_{\epsilon (T)}$, and we are done by the conservativity of $\epsilon$.
\hfill$\square$ 
  
{\it Remark.} Probably, statement (iii) of the proposition is valid
without the assumption of existence of DG structure on $D$.

{\bf 1.6. Semi-free modules.} 1.6.1.  Let $\CA$ be a DG category. To perform homotopy constructions on the DG level, one often needs to enlarge $\CA$ to a DG category $\LA$  of  ``DG ind-objects" defined as follows.

Let $\CA^\natural$ be a DG category whose objects are formal direct sums $M=\oplus M_i [n_i ]$, where $i$ runs over some set $I$, 
 $M_i \in \CA$, $n_i \in \Bbb Z$; one has $\Hom (\oplus M_i [n_i ], \oplus M'_j [n_j ]):=\mathop\Pi\limits_{i\in I} \mathop\oplus\limits_{j\in J} \Hom (M_i ,M'_j )[n_j - n_i ]$, the composition law is the product of matrices (with finitely many non-zero entries in each column). 
An endomorphism $\phi$ of $M$ is said to be {\it f-nilpotent} if $I$ admits a filtration $ I_0 \subset I_1 \subset \ldots$ such that $\cup I_a = I$ and the corresponding filtration $M_0 \subset M_1 \subset \ldots$, $M_a := 
\mathop\oplus\limits_{i\in I_a} M_i [n_i ]$, of $M$ has property $\phi (M_a )\subset M_{a-1}$, $\phi (M_0 )=0$.
 Now $\LA$ is formed by pairs $(M, \phi )$ where $M\in \CA^\natural$ and $\phi $ is its f-nilpotent endomorphism of degree 1 which satisfies the Maurer-Cartan equation $d(\phi )+ \phi\cdot\phi =0$;
 the complex $\Hom ((M,\phi_M ), (M',\phi_{M'} ))$ equals $\Hom (M,M')$ as a $\Bbb Z$-graded module, its differential $d_\phi$ is $ d_\phi (f):= d(f)  +[\phi ,f]$. We usually denote $(M,\phi )$ simply by $M$, and write $\phi_M $ for $\phi$. 

The DG category $\LA$ is cocomplete and strongly pretriangulated. One has a fully faithful embedding $\CA \hra \LA$ which identifies $M\in \CA$ with  $(M[0],0)\in  \LA$.

{\it Remarks.} (i) The essential image of the corresponding fully faithful embedding  $\CA^{\text{pretr}}\hra \LA $ is formed by those $(M,\phi )$ that the sum $M=\oplus M_i [n_i ]$ is finite. Replacing $\CA$ by $\CA^{\text{pretr}}$ does not change $\LA $.

(ii)  The embedding $\CA\hra \LA$ does {\it not} commute with infinite direct sums.
 Namely, let  $\{ M_\alpha \}$ be an infinite family of non-zero objects of $\CA$ such that the direct sum of $ M_\alpha$'s in $\CA$ exists; denote it by $\oplus_\CA M_\alpha \in \CA$, and let $\oplus M_\alpha$  be the direct sum in  $\LA$. Then the  canonical morphism 
$\oplus M_\alpha \to \oplus_\CA M_\alpha$ is {\it not} an isomorphism. Same is true if we replace $\CA$ and $\LA$ by their homotopy categories.

{\it Exercise.} Suppose $\CA^{\text{tri}}$ has arbitrary direct sums and
is  generated by a set of compact objects. Then the essential image of
$\CA^{\text{tri}} \to \LA^{\text{tri}}$ consists of those  $M\in \LA$ that
for every set of objects $ N_i\in \CA$, the canonical morphism
$\Hom (\bigoplus_{\CA^{\text{tri}}} N_i, M) \to \prod Hom (N_i, M)$
is an isomorphism.\footnote{Hint: use (i) of 
Proposition in 1.4.2 for $I=\CA^{\text{tri}}$ and $D$ formed by $M\in \LA$ as above.}

\proclaim{\quad Lemma}   If a DG category  $\CB$ is cocomplete, strongly pretriangulated, and DG Karoubian, then every DG functor  $F : \CA\to \CB$  extends uniquely (up to a unique DG isomorphism) to a DG functor $ \LF : \LA \subset  (\LA )^\kappa \to \CB$ that commutes with arbitrary direct sums.
\endproclaim

{\it Proof.} The above $\CA^\natural$ is a full DG subcategory of $\LA $ formed by objects $(M,0)$. Our $F$ extends uniquely to a DG functor $A^\natural \to \CB$ which commutes with arbitrary direct sums, and then to $  \CA^{\natural\, \text{pretr}}\subset \LA $. 
 Let $A\,\tilde{}\subset \LA $ be the full DG subcategory formed by arbitrary direct sums of objects of $ \CA^{\natural \,\text{pretr}}$. Then $F$ extends uniquely to a DG functor ${A\,\tilde{}}\to\CB$ which commutes with arbitrary direct sums, hence to ${A\,\tilde{}}^{\,\,\text{pretr}}$. To finish the proof, it remains to show that $({\CA\,\tilde{}}^{\,\,\text{pretr}})^\kappa \iso (\LA )^\kappa$. Notice that the construction of $hocolim$ from  1.4.1 lifts
in the evident manner to DG level (and it becomes canonical, as cones are). Now for any $M\in \LA$ the terms of the filtration $M_0 \subset M_1 \subset \ldots$ on $M$ (which comes since $\phi_M$ is  f-nilpotent)  lie in $  \CA^{\natural\, \text{pretr}}$, hence $hocolim  M_a \in {\CA\,\tilde{}}^{\,\,\text{pretr}}$. Since the evident projection $hocolim  M_a \to M$  admits a DG splitting, one has
 $M \in ({\CA\,\tilde{}}^{\,\,\text{pretr}})^\kappa$,  q.e.d.  \hfill$\square$

1.6.2. A DG category $\CA$ is said to be {\it homotopically Karoubian} if $\Ho\, \CA$ is Karoubian. For any $\CA$ its {\it homotopy idempotent completion} is a fully faithful DG embedding $\CA \hra\CB$ such that $\Ho\,\CB$ is the idempotent completion of $\Ho\,\CA$. Here is a concrete construction: for $\CB$, take  a full DG subcategory of $\LA$ formed by those objects that represent in $\Ho\, \LA$ direct summands of objects in $\Ho\, \CA$; since $\LA^{\text{tri}}$ is Karoubian by 1.4.1, our $\CB$ is a homotopy idempotent completion of $\CA$. This construction also shows that all homotopy idempotent completions of $\CA$ are naturally quasi-equivalent.

1.6.3.   Suppose a DG category $\CA$ is essentially small.
Denote by $R\text{-dgm}$ the DG category  of complexes of $R$-modules. The objects of DG category $\CA\text{-dgm}:= \CH om (\CA,R\text{-dgm})$ are called {\it left (DG) $\CA$-modules}.  A {\it right $\CA$-module} is the same as an $\CA^{\text{op}}$-module. 
$\CA^{\text{op}} \text{-dgm}$ is  strongly pretriangulated, cocomplete, and DG Karoubian.
One has the Yoneda  fully faithful DG embeddings $\CA\hra \CA^{\text{op}} \text{-dgm}$, $\CA^{\text{op}} \hra \CA$-mod.

{\it Remarks.}  (i) If $\CA $ has single object,  $\CA =A$, then $\CA$-module =  DG $A$-module.

(ii) There are natural DG functors $\CA^{\text{op}} \text{-dgm}\otimes \CA\text{-dgm}\to (\CA^{\text{op}} \otimes \CA )\text{-dgm} \to R\text{-dgm}$, the first one assigns to $F\otimes G$ the $\CA^{\text{op}} \otimes \CA $-module $M\otimes N\mapsto F(M)\otimes_R G(N)$, the second one is the coend functor. We denote the composition as $F\otimes G\mapsto F\otimes_\CA G$. In the situation of (i) this is the usual tensor product of left and right $A$-modules.

The DG functor  $\LA \to \CA^{\text{op}} \text{-dgm}$ defined by the embedding $\CA\hra \CA^{\text{op}} \text{-dgm}$ (see the  lemma in 1.6.1) is itself a fully faithful embedding whose essential image consists of  {\it semi-free}  $\CA^{\text{op}}$-modules. Here an $\CA^{\text{op}}$-module $M$ is said to be {\it semi-free} if it admits an increasing exhaustive filtration $M_0 \subset M_1 \subset \ldots$ such that  each $\gr_a M$ is {\it free}, i.e., isomorphic to a direct sum of translations of objects in $\CA$. Thus $M_{a+1}= \CC one (C_a \to M_a )$ where $C_a$ is a free $\CA^{\text{op}}$-module.   Replacing   $\CA$ by $\CA^{\text{pretr}}$ does not change these categories.

\proclaim{\quad Lemma} $ \LA {\!}^{\text{tri}}= \LA {\!}^{\kappa\, \text{tri}}$ is a right-admissible subcategory of $ \CA^{\text{op}} \text{-dgm}^{\text{tri}}$. 
\endproclaim

{\it Proof.} See 1.4.1 and (i) of the proposition in 1.4.2.  \hfill$\square$

1.6.4. An $\CA^{\text{op}}$-module is said to be {\it acyclic} if the functor $\CA^{\text{op}} \to R\text{-dgm}$ takes values in acyclic complexes. Those objects form a full pretriangulated DG subcategory $\CA^{\text{op}} \text{-dgm}_{ac} \subset \CA^{\text{op}} \text{-dgm}$, and 
 $ \CA^{\text{op}} \text{-dgm}^{\text{tri}}_{ac}$ is the right orthogonal complement to $ \LA {\!}^{\text{tri}} $. By 1.2 and  the above lemma, $\CA^{\text{op}} \text{-dgm}^{\text{tri}}_{ac}$ is left admissible, its left orthogonal complement  equals $ \LA {\!}^{\text{tri}}$, and the restriction of the projection $\CA^{\text{op}} \text{-dgm}^{\text{tri}} \to D(\CA^{\text{op}} ):=   \CA^{\text{op}} \text{-dgm}^{\text{tri}}/ \CA^{\text{op}} \text{-dgm}^{\text{tri}}_{ac}$ (the {\it derived category} of right $\CA$-modules)  to $
 \LA {\!}^{\text{tri}}$ is an equivalence of categories:  $$   \LA {\!}^{\text{tri}} \iso D(\CA^{\text{op}} ).  \tag 1.6.1$$
Thus $ \LA$ is a DG structure on $D(\CA^{\text{op}} )$; we refer to it as the {\it canonical DG structure}.

{\it Remark.}  The procedure from (a) of the proof of the proposition in 1.4.2 (performed on DG level) provides for any $M\in \CA \text{-dgm}$  its semi-free left resolution $hocolim (M_a ,i_a )\to M$. In fact,   one can replace $hocolim$ by the plain direct limit $\limright M_a$ (which is semi-free):  the  morphism $hocolim (M_a ,i_a )\to M$ factors through  $\limright M_a$, and $hocolim (M_a ,i_a )\twoheadrightarrow \limright M_a$ is a homotopy equivalence.

  A DG functor $F:\CA_1 \to \CA_2$ yields an evident DG functor $F_* : \CA_2^{\text{op}} \text{-dgm}\to \CA_1^{\text{op}} \text{-dgm}$ and its left adjoint $F^* : \CA_1^{\text{op}} \text{-dgm} \to \CA_2^{\text{op}} \text{-dgm}$ that extends $F$ and sends $\mathop{\vtop{\ialign{#\crcr
  \hfil\rm $\CA_1$\hfil\crcr
  \noalign{\nointerlineskip}\rightarrowfill\crcr
  \noalign{\nointerlineskip}\crcr}}} $ to $\mathop{\vtop{\ialign{#\crcr
  \hfil\rm $\CA_2$\hfil\crcr
  \noalign{\nointerlineskip}\rightarrowfill\crcr
  \noalign{\nointerlineskip}\crcr}}}$;  
 they commute with arbitrary direct sums. Since $F_*$ preserves acyclic modules (or quasi-isomorphisms), it yields a triangulated functor $F_* : D(\CA_2^{\text{op}} ) \to D(\CA_1^{\text{op}} )$ whose left adjoint comes via (1.6.1) from $F^* :   \mathop{\vtop{\ialign{#\crcr
  \hfil\rm $\CA_1$\hfil\crcr
  \noalign{\nointerlineskip}\rightarrowfill\crcr
  \noalign{\nointerlineskip}\crcr}}} {\!}^{\text{tri}} \to \mathop{\vtop{\ialign{#\crcr
  \hfil\rm $\CA_2$\hfil\crcr
  \noalign{\nointerlineskip}\rightarrowfill\crcr
  \noalign{\nointerlineskip}\crcr}}}{\!}^{\text{tri}}$.

\proclaim{\quad Lemma}   If $F$ yields an equivalence $\CA^{\text{tri}\,\kappa}_1 \iso\CA^{\text{tri}\,\kappa}_2$, then  $F^* ,F_*  : D(\CA_1^{\text{op}} ) \leftrightarrows D(\CA_2^{\text{op}} ) $ are mutually inverse equivalences of triangulated categories. \hfill$\square$
\endproclaim

{\bf 1.7. Cocomplete vs.~small $\&$ Karoubian: DG setting.}  The next result is due to V.~Drinfeld.

Let $\CB$ be a  cocomplete pretriangulated DG category. An object of $\CB$ is said to be {\it perfect} if it is perfect (i.e., compact) as an object of  $ \CB^{\text{tri}}$. All perfect objects form a full DG subcategory $\CB^{\text{perf}} \subset\CB$. E.g.,
 $\CA \hra  \LA$ takes values in  $(\LA )^{\text{perf}}$.

Consider the following two classes of DG categories:

(a)  Pretriangulated DG categories $\CA$ with $\CA^{\text{tri}}$  essentially small and Karoubian.

(b) Cocomplete DG Karoubian strongly pretriangulated  DG categories  $\CB$ such that $ \CB^{\text{tri}}$ is   generated by a set of perfect objects. 

 They form naturally ``DG 2-categories" $\CD\CG^{(a)}$ and $\CD\CG^{(b)}$. Namely, 1-morphisms in $\CD\CG^{(a)}$ are DG functors, 1-morphisms in $\CD\CG^{(b)}$ are DG functors that commute with arbitrary direct sums and send perfect objects to perfect ones, and  
 the ``complex" of 2-morphisms between DG functors is defined as in 1.5.3 (the quotation marks mean that morphisms need not form sets).

For $\CA$ of type (a) the category $(\LA )^\kappa$ is of type (b); for $\CB$ of type (b) 
$\CB^{\text{perf}}$ is of type (a) (see (i) of the proposition in 1.4.2). By the lemma in 1.6.1, we have defined  ``adjoint DG functors" $$
\CD\CG^{(a)}\leftrightarrows \CD\CG^{(b)}. \tag 1.7.1$$

\proclaim{\quad Proposition} The adjunction 1-morphisms $\CA \to (\LA )^{\kappa\,\text{perf}}$,
$(\mathop{\vtop{\ialign{#\crcr
  \hfil\rm $\CB^{\text{perf}}$\hfil\crcr
  \noalign{\nointerlineskip}\rightarrowfill\crcr
  \noalign{\nointerlineskip}\crcr}}})^{\kappa}\to \CB$ are quasi-equivalences.
\endproclaim
 
 {\it Proof.} The assertion is implied by the next lemma, which follows from 1.4.2: 
 
 \proclaim{\quad Lemma} If $\CA \subset \CB^{\text{perf}}$ is  any small full DG subcategory that   generates $\CB^{\text{tri}}$, 
then  $\LA \to \CB$ is a quasi-equivalence, and
  $(\CB^{\text{perf}})^{ \text{tri}} = (\CA^{\text{tri}})^\kappa $.  \hfill$\square$
  \endproclaim
 
{\it Remarks.} (i) The proposition asserts that  $\CA \mapsto \LA$ and $\CB \mapsto \CB^{\text{perf}}$ are ``homotopically mutually inverse" correspondences between the categories of type (a) and of type (b). {\it Variant:} Let us enlarge (b) to the class (b)$'$ of cocomplete  pretriangulated  DG categories  $\CB$ such that $ \CB^{\text{tri}}$ is  generated by a set of perfect objects. For any such
 $\CB$   the embedding $\CB \hra \CB^{\text{pretr}\,\kappa}$ is a quasi-equivalence (see 1.4.1), and $\CB^{\text{pretr}\,\kappa}$ is of type (b). Thus  $\CA \mapsto \LA$ and $\CB \mapsto \CB^{\text{perf}}$ are ``homotopically mutually inverse" correspondences between the categories of type (a) and of type (b)$'$.

(ii) For any essentially small DG category $\CA$ the category  $ {\LA }^{\text{perf}}$  is a homotopy idempotent completion of $\CA^{\text{pretr}}$, see 1.6.2.

 {\bf 1.8. Homotopy DG quotients}.    Let $\CA$ be a  DG category such that $\CA^{\text{tri}}$ is essentially small, $\CI\subset\CA$ be a full DG subcategory. By 1.4.2, ${\LI }^{\text{tri}}\subset {\LA }^{\text{tri}}$ is a right-admissible subcategory. As in  1.5.5, let    $ \CB := \LA / \LI\subset \LA$ be the preimage of $({\LI }^{\text{tri}})^{\bot}\subset {\LA }^{\text{tri}}$, so one has  $\CB^{\text{tri}}\iso 
 \LA {}^{\text{tri}}/ \LI {}^{\text{tri}}$. Let $\CA /\CI \subset \CB$ be the preimage of the essential image of $\CA^{\text{tri}}$ in $
  \LA {}^{\text{tri}}/ \LI {}^{\text{tri}}$; this is  {\it the (right) Keller  DG quotient}.\footnote{See \cite{K} and \cite{Dr} \S 4. Notice that the construction is not self-dual.} As in Remark in 1.5.5, one has a natural homotopy DG functor $\CA \leftarrow \tilde{\CA}\to\CA /\CI$, where   $(\CA /\CI )^{\text{tri}}= \CA^{\text{tri}}/\CI^{\text{tri}}$.\footnote{So, by 1.5.2, one can assign to each object of $\CA$ a homotopy object of $\CA /\CI$.} We also set $\CA /^{\kappa }\CI := \CB^{\text{perf}}$, so  $(\CA /^{\kappa }\CI )^{\text{tri}} = \CA^{\text{tri}}/^\kappa \CI^{\text{tri}}$ (see 1.4.2). By the lemma in 1.7, one has natural quasi-equivalences $$\mathop{\vtop{\ialign{#\crcr
  \hfil\rm $\CA/\CI$\hfil\crcr
  \noalign{\nointerlineskip}\rightarrowfill\crcr
  \noalign{\nointerlineskip}\crcr}}}
  \iso \mathop{\vtop{\ialign{#\crcr
  \hfil\rm $\CA/^\kappa \CI$\hfil\crcr
  \noalign{\nointerlineskip}\rightarrowfill\crcr
  \noalign{\nointerlineskip}\crcr}}}
   \iso \LA /\LI  .\tag 1.8.1$$
 
{\it Remark.}    Small DG categories equipped with a marked zero object form naturally a pointed model category  (with quasi-equivalences as weak equivalences and surjective functors as fibrations),\footnote{It has a peculiar property:  suspensions of its objects are contractible.} see \cite{Dr} B4--6. This plain model structure is not  useful, for it does not consider morphisms between DG functors (cf.~Remark in 1.5.3). For $\CA$, $\CI$ as above, Drinfeld \cite{Dr} defines  a {\it (homotopy) DG quotient} as a homotopy DG functor $\CA \leftarrow \tilde{\CA}\to \CA /\CI$  where  the functor $\Ho\,\CA \to \Ho (\CA /\CI )$ is essentially surjective and yields an equivalence $\CA^{\text{tri}}/\CI^{\text{tri}}\to (\CA /\CI )^{\text{tri}}$. Thus $\CA /\CI$ is  a homotopy cone for the morphism $\CI \to \CA$. The Keller construction is a concrete way to construct a  DG quotient;\footnote{To retain smallness, one should replace  $\CA /\CI$ (which is not small) by  a small DG subcategory having the same homotopy category.} for another one see \S3 of \cite{Dr}. For the uniqueness property of homotopy DG quotients see \cite{T}; a weaker fact is in \cite{Dr}1.6.2.

  \medskip

{\bf 1.9. \bf Homotopy tensor DG categories.} For us, a {\it tensor triangulated category} is a triangulated category equipped with a tensor structure such that $\otimes$ is an exact bifunctor. Below we discuss DG liftings of this notion.

1.9.1. Let $\CA$ be a DG category. It is clear what a {\it DG tensor structure} on $\CA$ is.\footnote{This is a tensor structure such that  $\otimes : \CA \otimes\CA \to \CA$ is a DG functor and the associativity and commutativity constraints  are DG morphisms.} Unfortunately, DG  tensor structures are too rigid to pass readily to homotopy DG quotients. The next definition, due to V.~Drinfeld,  is  adapted for the {\it right} Keller localization:

A  {\it (right) homotopy  tensor  structure}  is a  pseudo-tensor structure on  $\CA$  (see \cite{BD} 1.1.3) such that the corresponding pseudo-tensor structure on $\Ho \CA$ is actually a tensor structure. Therefore for every finite non-empty set $I$ we have a DG functor $P_I : \CA^{{\text{op}} \otimes I}\otimes \CA \to R\text{-dgm}$, $(\otimes M_i )\otimes N\mapsto P_I (\{ M_i \} ,N)$ of $I$-polylinear operations,  the operations compose according to a usual format (and $P_I =\Hom$ if $|I|=1$), and $H^0  P_I (\{ M_i \} ,N) = \Hom_{\Ho\CA} (\otimes M_i ,N)$ for a tensor structure on $\Ho\,\CA$ uniquely defined by these operations.
  
Notice that any  DG tensor or homotopy  tensor structure extends in a unique way to $\CA^\kappa$ and $\CA^{\text{pretr}}$.  It yields  tensor structures on $\Ho\,  \CA$ and $\CA^{\text{tri}}$; the latter becomes a tensor triangulated category. 

{\it Exercise.} The tensor product   $\Ho\, \CA^{\otimes I}\to\Ho\, \CA$  lifts naturally to a homotopy  functor (see 1.5.2). So
for  $M_i \in\CA$ the tensor product  $\otimes M_i$ is well defined as a homotopy object of $\CA$.

{\it Remark.} The above definition of homotopy tensor structure is not self-dual:  algebras are pleasant to consider, while coalgebras are cumbersome objects. It would be  nice to find a definition of homotopy tensor structure  free of this nuisance.

1.9.2.  Any DG tensor structure on $\CA$ extends naturally to $\LA$ so that $\otimes$ commutes with arbitrary direct sums. Similarly,
 any homotopy  tensor structure on $\CA$ extends naturally to $\LA$, so that $P_I (\{ M_i \} ,N)$ transforms arbitrary direct sums of $M_i$'s to direct products, and if $M_i \in\CA$, then it 
 transforms arbitrary direct sums of $N$'s to direct sums. The lemma in 1.7 implies that   $(\CA^{\text{tri}})^\kappa = ({\LA}^{\text{perf}})^{\text{tri}}$ is a tensor subcategory of ${ \LA}^{\text{tri}}$, so ${\LA}^{\text{perf}}$ is a homotopy  tensor category.

Suppose $\CA$ is essentially small. Then any homotopy  tensor structure on $\CA$ yields, by means of (1.6.1), a tensor structure on $D(\CA^{\text{op}} )$ such that $\CA^{\text{tri}}\subset  (\CA^{\text{tri}})^\kappa$ are its  tensor subcategories. 

 If  $\CI\subset \CA$ is a full DG subcategory such that $\CI^{\text{tri}}$ is an $\otimes$-ideal in $\CA^{\text{tri}}$, then 
 $\CA /\CI \subset\CA /^\kappa \CI \subset\LA /  \LI $ are homotopy tensor DG subcategories of $\LA$. Thus the (right) Keller DG quotient inherits naturally a homotopy tensor structure. It lifts the evident tensor structure on $\CA^{\text{tri}}/\CI^{\text{tri}}$.
 
1.9.3. For  $M_i \in \CA$ with $i\in I$ as above, we have a left $\CA$-module 
 $P'_I (\{ M_i \})$, $N\mapsto P_I (\{ M_i \} ,N)$ (here $N\in\CA $). For any right $\CA$-module $F$ set $P_I ( \{ M_i \} ,F) := F\otimes_\CA P'_I (\{ M_i \})\in R$-dgm (see \cite{Dr} C3 for the notation); if $F\in\LA$, this complex coincides with its namesake from 1.9.2.

Denote by $\CP_I  (\{ M_i \},F)$ the right $\CA$-module  $Q\mapsto P_{\tilde{I}} ( \{ M_i ,Q \},F)$; here $Q$ is a test object in $\CA$ and $\tilde{I}$ is $I$ with one element added.
 Thus if $\CA$ is a tensor category, then $\CP_I (\{ M_i \} ,F)(Q)=F((\otimes M_i )\otimes Q)$. The DG endofunctor $\CP_I (\{ M_i \} ,\cdot ):
 F\mapsto \CP_I (\{ M_i \},F)$ of $\CA^{\text{op}}$-dgm
   preserves quasi-isomorphisms, so it yields a  triangulated (poly)functor $\CP_I  : (\CA^{\text{tri}})^{ {\text{op}} I} \times D(\CA^{\text{op}} )\to D(\CA^{\text{op}} )$. 
   
 The functors   $\{ M_i \} \mapsto P_I ( \{ M_i \} ,F), \CP_I ( \{ M_i \} ,F)$ extend naturally to $M_i\in \LA$ so that arbitrary direct sums of $M_i$'s are transformed to direct products. The corresponding DG endofunctor $ \CP_I (\{ M_i \},\cdot )$ of $\CA^{\text{op}}$-dgm preserves quasi-isomorphisms  if $M_i \in \LA^{\text{perf}}$ (see the lemma in 1.7). If $|I|=1$, we write $\CH om (M,F):=\CP_I (\{ M\} ,F)$.  
 
{\it Remark.} The endofunctor $ \CP_I (\{ M_i \} ,\cdot )$ of $D(\CA^{\text{op}} )$ is right adjoint to one $F\mapsto (\otimes M_i )\otimes F$. Thus $\CH om$ is an action of the tensor category $(\CA^{\text{tri}})^{\kappa\,{\text{op}}}$ on $D(\CA^{\text{op}} )$.

{\bf 1.10. Presheaves.}   Let $T$ be an essentially small category. Assume that for every finite family $\{ X_s \}$ of objects in $T$ the coproduct $\sqcup X_s$ exists, and  for every morphism $ Y\to X=  \sqcup X_i$ the fiber products $Y_s :=Y\times_X X_s $ are representable and $\sqcup Y_s \iso Y$.
For us, a presheaf  on $T$ always means a  presheaf  of $R$-modules  compatible with finite coproducts, i.e., a functor  $F: T^{\text{op}} \to R$-mod such that $F(\sqcup X_s )\iso \Pi_s F(X_s )$.\footnote{That is, a sheaf  for the stupid topology whose coverings $\{ X_s \to X\}$ are finite decompositions $\sqcup X_s \iso X$.} 
Let $R[T]$ be an $R$-category whose objects are denoted by $R[X]$, $X\in T$, and $\Hom (R[X], R[Y])$  is the direct limit of the directed family of free $R$-modules $\Pi_s R[\Hom (X_s ,Y)]$ labeled by finite decompositions $\sqcup X_s \iso X$. Then a presheaf on $T$ is the same as an $R [T ]^{\text{op}}$-module; we identify $R[X]$ with the corresponding representable presheaf. 

We have an
 abelian $R$-category $\CP\CS h$ of presheaves and a DG category $C\CP\CS h= R[T]^{\text{op}}$-dgm of complexes of presheaves. Let $\mathop{\vtop{\ialign{#\crcr
  \hfil\rm $R[T]$\hfil\crcr
  \noalign{\nointerlineskip}\rightarrowfill\crcr
  \noalign{\nointerlineskip}\crcr}}}$   be the category of semi-free DG $R [T ]^{\text{op}}$-modules  (see 1.6.3). Set
$\CP := \mathop{\vtop{\ialign{#\crcr
  \hfil\rm $R[T]$\hfil\crcr
  \noalign{\nointerlineskip}\rightarrowfill\crcr
  \noalign{\nointerlineskip}\crcr}}}^{\text{perf}}$, see 1.7. Then $\mathop{\vtop{\ialign{#\crcr
  \hfil\rm $R[T]$\hfil\crcr
  \noalign{\nointerlineskip}\rightarrowfill\crcr
  \noalign{\nointerlineskip}\crcr}}}\iso \LP$ is the canonical DG structure on the derived category $D\CP\CS h$, see (1.6.1).

Suppose $T$ is closed under finite direct products;  they define a tensor  structure on $T$, so $R [T ]$ is an tensor $R$-category. By 1.9.2,  $\mathop{\vtop{\ialign{#\crcr
  \hfil\rm $R[T]$\hfil\crcr
  \noalign{\nointerlineskip}\rightarrowfill\crcr
  \noalign{\nointerlineskip}\crcr}}}$ is a  tensor DG category, and $\CP$ its  tensor DG subcategory. We get a tensor structure on $\LP^{\text{tri}}=D\CP\CS h$. 
Notice that $\CP\CS h$ is a tensor abelian $R$-category in the 
usual way,\footnote{i.e., as the category of sheaves on the stupid topology from the previous footnote.} and the embedding $\LP\hookrightarrow C\CP\CS h$ is compatible with the tensor products, so the tensor structure on $D\CP\CS h$ is the standard one. The tensor DG category $\LP^{\text{op}}$ acts on $C\CP\CS h$ by the $\CH om$ action;  for $X, Y\in T$ and $F\in C\CP \CS h$ one has $\CH om (R [X],F)(Y)= F(X\times Y)$.

{\it Notation.} If $U_\cdot$ is a simplicial object of $T$, then $R [U_\cdot ]$ is a simplicial presheaf; by abuse of notation, we denote the  normalized complex also by $R [U_\cdot ]$. For a morphism $Y\to X$ we write $ R [Y/X]:= \CC one (R [Y]\to R [X])[-1]$, for an augmented simplicial object $U_\cdot /X$ we have $R[U_\cdot /X]$, etc.   So $R [X]\in\CP$, and $R [U_\cdot ]$, $R [U_\cdot /X] \in \LP$.

{\bf 1.11. Sheaves.} Suppose that our $T$ is equipped with  a Grothendieck topology  $\CT$ (see \cite{SGA4} II or \cite{D1} Arcata I  6) such that for every finite family $X_s \in T$ the morphisms $\{ X_s \to \sqcup X_s \}$ form a covering of $\sqcup X_s$. We have the abelian category $\CS  h^\CT$ of sheaves of $R$-modules, the exact sheafification functor $\CP \CS h \to\CS h^\CT$, $F\mapsto F_\CT $, and its right adjoint $\CS h^\CT \to \CP \CS h$ which is fully faithful and left exact.

A complex $F$ of presheaves is said to be {\it  locally acyclic}  if for every $X\in T$ and $a\in H^n F (X)$ those $X$-objects $V/X$  that   $a|_{V}\in H^n F (V )$ vanishes, form a 
 covering sieve of $X$. All such $F$ form a thick subcategory $I^\CT$ of $D\CS h^\CT$;
the sheafification functor yields an equivalence  $D\CP\CS h /I^\CT \iso D\CS h^\CT$.
 
For $X\in T$ its {\it hypercovering} (see \cite{SGA4} V 7.3)  is  an augmented simplicial object of $T$ which satisfies the next condition: Let $f: \partial \Delta^n \times Y /Y \to U_\cdot /X$ be any
  morphism of augmented simplicial objects; here $Y\in T$,  $n\ge 0$, $\Delta^n$ is the ``$n$-simplex" simplicial set, $\partial \Delta^n$ its boundary. Then those  $V/Y$ that   $f_V :   \partial \Delta^n \times V /V \to U_\cdot /X$  can be extended to  $\Delta^n \times V /V \to U_\cdot /X$, form a covering sieve of $Y$. 

Assume that every $X\in T$ is quasi-compact,\footnote{Which means (see \cite{SGA4} VI 1.1) that every covering  of  $X$ has a finite subcovering. }  has finite cohomological dimension,  and  each of its coverings  admits a refinement  $\{ Y_s \to X\}$ such that for every morphism $Z\to X$ the fiber products $Y_s \times_X Z$ are representable.

 \proclaim{\quad Lemma} (i) $I^\CT$ is a right-admissible subcategory of   $D\CP\CS h$.
 
(ii)  If   $U_\cdot /X$ is a  hypercovering of $X\in T$, then the complex $R [U_\cdot / X ]$ is  $\CT\!$ locally acyclic. The category   $I^\CT$ is  generated by all  such complexes  $R [U_\cdot / X ]$.  \hfill$\square$
\endproclaim

The objects of $I^{\CT \bot}$ are called {\it $\CT\!$-local} (or simply local) complexes;  by the lemma,  a complex $G$ is local if and only if for every hypercovering  $U_\cdot /X$ the map $G(X)\to G(U_\cdot ):= \Hom (R[U_\cdot ],G)$ is a quasi-isomorphism. Notice that if such $G$ is a single presheaf, then it is automatically a sheaf.
As in 1.2, one has the   
  {\it $\CT\!$-localization}  endofunctor of $D\CP\CS h $, $F \mapsto \CC^\CT (F):=F_{I^{\CT\bot}}$ which identifies $D\CS h^\CT$ with $I^{\CT\bot}$. 
 The t-structure picture here is from (i) of the lemma from loc.~cit.
 As  in 1.5.5, $\LP$ yields DG structures on the subcategories $I^\CT$ and $I^{\CT\bot}$.

  If $\CT$ has enough   points (see \cite{SGA4} IV 6), then $I^\CT$ is formed by  complexes of sheaves whose stalks at the points are acyclic complexes. The $\CT\!$-localization  can be lifted to a DG endofunctor  of  $C\CP\CS h$ which assigns to a complex of presheaves  its Godement resolution.

 \bigskip
 
\centerline{\bf \S 2. Correspondences and effective motives: first definitions }

\medskip

This section treats the material of lectures 1, 2,  8, 23 of \cite{MVW}, \S 6 of \cite{SV}, and 2.1.1--2.1.3, 3.4.1 of \cite{Vo2}. At variance with loc.~cit.,  the constructions are carried on DG level, and  the ``bounded above"  assumption on complexes is discarded.

{\bf 2.1. Presheaves  with transfer.} We fix a base field $k$; let $\CS ch =\CS ch_k $ be the category of all  $k$-schemes,  $\CS m =\CS m_k$ be the  subcategory of smooth  varieties. Notice that $\CS m$ is essentially small.

2.1.1. Using 1.10 for $T=\CS m$,  we get the abelian category $\CP\CS h$  of presheaves of $R$-modules on $\CS m$, the tensor DG categories  $\CP$ and $\LP$, the equivalence of triangulated tensor categories $$ { \LP}^{\text{tri}} \iso D\CP\CS h , \tag 2.1.1$$ etc.   Any scheme
 $X$ yields a presheaf $R [X] : Y\mapsto R [X](Y)=R [X(Y)]$ on $\CS m$. As in 1.10, for $Y\to X$ we have a complex $R [Y/X]:= \CC one ( R [Y]\to R [X])[-1]$,  a simplicial augmented scheme $U_\cdot / X$ yields  a complex $R [U_\cdot / X ]$, etc.

Any presheaf $F$ on $\CS m$ extends naturally to a presheaf on $\CS ch$: for a scheme $Y$ we define $F(Y)$ as $\limright F(X)$, the limit is taken over the category of pairs $(X,\nu )$ where $X\in\CS m$ and $\nu $ is a morphism $ Y\to X$.\footnote{A morphism $(X,\nu )\to (X',\nu' )$ is a morphism $f: X\to X'$ such that $f\nu =\nu'$.} If $Y$ can be represented as the projective limit 
 of a directed family of smooth varieties $Y_\alpha$ connected by affine  morphisms, then   
 $F(Y)= \limright F( Y_\alpha )$ (see \cite{EGA IV} 8.8); in particular, the functor $F\mapsto F(Y)$ is exact.  Examples of such $Y$'s:  (i)  For a point $x$ of a smooth variety $X$ such is  the local scheme $X_x := \Spec\, \CO_{Xx}$ and its  Henselization $X^h_x$;  
(ii) If $K/k$ is a separable extension, then such is any $Y\in \CS m_K$.

2.1.2. If $X, Y$ are schemes and $Y$ is normal, then we denote by $\CC or (Y,X)=\CC or_R (Y,X)$ the free $R$-module generated by the set of  reduced irreducible subschemes $\Gamma \subset X\times Y$ such that the projection $\Gamma \to Y$ is finite and its image is an irreducible component of $Y$ (the latter condition can be replaced by that $\dim \Gamma =\dim Y$). Its elements $\gamma =\Sigma a_i \Gamma_i$, $a_i \in R$, are referred to as {\it finite correspondences}. For  $Z$  normal, $Y$  smooth,  we have the composition map $\CC or (Y,X)\otimes\CC or (Z,Y)\to\CC or (Z,X)$, $\gamma\otimes\gamma' \mapsto 
\gamma\gamma'$ := the push-forward by the projection $X\times Y\times Z \to X\times Z$ of the cycle $(X\times \gamma'  )\cap( \gamma \times Z)$ on $X\times Y\times Z$. Here $\cap$ is the cycle-theoretic intersection, which is well defined due to smoothness of $Y$ (our  cycles intersect properly).  The composition is associative.

{\it Remarks.}  (a) Informally, finite correspondences are multi-valued maps; the composition is the usual composition of multi-valued maps.
(b) Any  $\gamma\in \CC or (Y,X)$ is uniquely determined by its value at the generic point(s) $\gamma (\eta_Y )$ of $Y$, which is a zero cycle on $X_{\eta_Y}$. (c)  By Noether's normalization lemma, every (reduced, irreducible) affine variety $\Gamma$ can be realized as a finite correspondence between some vector spaces $X$, $Y$. (d) The projection $\Gamma\to Y$ need not be flat: indeed, the flatness implies that $\Gamma$ is a  Cohen-Macaulay scheme.

Let $\CC or (\CS m ) =\CC or (\CS m_k )$ be the $R$-category whose objects 
are smooth varieties and morphisms are finite correspondences.
We denote $X\in \CS m$ considered as an object of $\CC or (\CS m)$ by  $\Z [X]$.  There is a faithful embedding $R [\CS m ] \hra\CC or (\CS m)$, $R [X] \mapsto \Z [X]$, which sends a morphism $f: Y\to X$ to its graph considered as a cycle on $X\times Y$. As in 1.10, $f$ yields a complex $\Z [Y/X]$, etc.

We write $\Z$ for $\Z [\Spec\, k ]$. For
   $X\in \CS m$  the structure projection $X\to \Spec\, k$ yields the {\it augmentation} morphism $aug_X :\Z [X]\to \Z$ and the 2-term complex
$\Z [X/\Spec \, k ]$; any 0-cycle of degree 1 on $X$  yields a splitting of the augmentation morphism and an isomorphism $\Z [X/\Spec \, k] \oplus  \Z  \iso \Z [X]$. If $Y$ is connected, then $\Hom (\Z [Y],\Z )=R$, so we have a morphism of $R$-modules $deg : \CC or (Y,X)\to R$, $\gamma \mapsto aug_X \gamma$. One has
$deg (\gamma\gamma' )=deg (\gamma )deg (\gamma' )$.

One refers to $\CC or (\CS m)^{\text{op}}$-modules as {\it presheaves with transfers} (on $\CS m_k$); they form an abelian category $\CP\CS h_{\tr}$. 
The evident forgetful functor $o:\CP \CS h_{\tr} \to \CP\CS h$ is exact and faithful, so presheaves with transfers can be considered as presheaves with an extra structure (that of transfer). Notice that any scheme $X\in\CS ch$ 
yields a presheaf with transfers $Y\mapsto \CC or (Y,X)$ which we denote by 
 $\Z [X]$. We have defined a functor $ \CS ch \to \CP\CS h_{\tr}$.

{\it Exercises.}
(i) Let $Y/X$ be a Galois \'etale covering of smooth varieties with the Galois group $G$, and 
$F$ be a presheaf with transfers such that the multiplication by $|G|$ on $F$ is invertible. Then $F(X)\iso F(Y)^G$.

(ii) If $K/k$ be a separable extension, then any presheaf with transfers on $\CS m_k$ defines naturally a 
presheaf with transfers on $\CS m_K$ (as a mere presheaf it was defined in 2.1.1).

2.1.3. The next approach to $\CC or$ is sometimes convenient. 
 Suppose $R=\Bbb Z$.  A correspondence $\gamma =\Sigma a_i \Gamma_i \in\CC or (Y,X)$ is said to be {\it effective} if all $a_i \ge 0$, so $\CC or (Y,X)$ is the group completion of the commutative monoid $\CC or(Y,X)^{\text{eff}} $ of effective correspondences. The composition of effective correspondences is effective.
 
 Suppose $X$ is quasi-projective; then $\Sym^\cdot X = \mathop\sqcup\limits_{n\ge 0} \Sym^n X$ is a universal commutative monoid scheme  generated by $X$. 
 Notice that  if   $Y$ is another quasi-projective scheme, then  $(\Sym^\cdot X)(Y)$ is the same as the set of morphisms of monoid schemes   $\Sym^\cdot Y \to \Sym^\cdot X$.
 
 \proclaim{\quad Proposition} There is a natural morphism of  monoids $$\sigma : \,\CC or(Y,X)^{\text{eff}}\to (\Sym^\cdot X)(Y). \tag 2.1.2$$  If $p=$ char $k$ is zero, this is an isomorphism; if $p$ is finite, then $\sigma$ becomes an isomorphism after inverting $p$ (so $\sigma $ is injective, for $\CC or$ has no torsion). The map $\sigma$ is compatible   with composition (for quasi-projective $Y$).
 \endproclaim
 
{\it Proof.} We can assume that $Y$ is irreducible.

 (a) We define $\sigma$ on the base of $\CC or(Y,X)^{\text{eff}} $ formed by irreducible correspondences $\Gamma \subset X\times Y$. Let $n$ be the degree of $\Gamma$. Write $\Gamma =\Spec \, \CA_Y$ where $\CA_Y$ is a coherent $\CO_Y$-algebra; 
let $U\subset Y$ be the open subset such that $\CA_U$ is flat. The $Y$-scheme $\Sym^n (\Gamma /Y )=\Spec\, (\CA_Y^{\otimes n})^{\Sigma_n}$ has a natural section $s_{\Gamma }$ over $U$: for $f\in \CA_U$ the function $s_\Gamma^* f^{\otimes n}\in \CO_U$  is the determinant of the multiplication by $f$ action on the $\CO_U$-module $\CA_U$. Since $Y$ is normal, $s_\Gamma$ extends uniquely to a section  of $ \Sym^n (\Gamma /Y )$ over $Y$, and (2.1.2) assigns to $\Gamma$ the composition  $Y\buildrel{s_\Gamma}\over\lra   \Sym^n (\Gamma /Y )\to \Sym^n X$.  

{\it Remark.} For $Y$ regular, we have an invertible $\CO_Y$-module $\det_{\CO_Y} \CA_Y$; the group  Aut$_{\CO_Y}(\CA_Y )$ acts on it by transport of structure. For  $f\in \CA^\times_Y$ the function $s_\Gamma^* f^{\otimes n}$ is the action of the multiplication by $f$ automorphism of $\CA_Y$ on $\det_{\CO_Y} \CA_Y$.

(b) Let us show that $\gamma_V \in\CC or(V,X)^{\text{eff}}$, $V\subset Y$  an open subset, comes from $ \CC or(Y,X)^{\text{eff}} $ if (and only if) the map $\sigma (\gamma_V ): V\to \Sym^n X$, $n= deg(\gamma_V )$, extends to $Y$. Let $Z_V \subset X\times V$ be the support of the cycle $\gamma_V$, $Z$ its closure in  $X\times Y$; then $\gamma_V$ extends to $Y$ if (and only if) $Z$ is finite over $Y$. Suppose $\sigma (\gamma_V )$ extends to  $f: Y\to \Sym^n X$. Consider  $X^n $ as a finite $ \Sym^n X$-scheme; let $\tilde{Z}:= X^n \mathop\times\limits_{\Sym^n X} Y \subset X^n \times Y$ be its $f$-pull-back, and $Z' \subset X\times Y$ be the image of $\tilde{Z}$ by any of the projections $X^n \to X$. Then $Z'$ is finite over $Y$ (for such is  $\tilde{Z}$) and $Z' \supset Z_V$, so $Z$ is finite over $Y$, q.e.d.

(c)  By (b), it suffices to check the isomorphism assertion
  at the generic point $\eta$ of $Y$. Replacing $k$ by $\eta$, $X$ by $X_\eta$, we can assume that $Y=\Spec\, k$ (so $\CC or (Y,X)$ is the group $Z_0 (X)$ of zero cycles on $X$). It is clear  that $\sigma$ is an isomorphism if
 $k$ is algebraically closed; by elementary Galois theory same is true if $k$ is perfect. The case of arbitrary $k$ reduces to  its perfect closure, since for any purely inseparable finite extension $k'/k$ both compositions of the usual maps $Z_0 (X)\rightleftarrows  Z_0 (X_{k'})$ are multiplications by $[k':k]$ (and $(\Sym^\cdot X)(k)\to    (\Sym^\cdot X)(k') $ is injective), q.e.d.
 
 (d)  Suppose $Z$ is normal, $Y$ is quasi-projective and smooth; let us check that for $\gamma \in\CC or(Y,X)^{\text{eff}} $, $\gamma' \in \CC or(Z,Y)^{\text{eff}} $ one has $\sigma (\gamma\gamma' )= \sigma (\gamma )\sigma (\gamma' )$, where $\sigma (\gamma )$ is understood as a morphism of monoid schemes $\Sym^\cdot Y \to \Sym^\cdot X$. It suffices to check the equality at a generic point  of $Z$; by the base change to the generic point(s) of $\gamma'$, we are reduced to the case $Z=\Spec\, k$, $\gamma' =y\in Y(k)$.
 One can assume that $\gamma$ is as in (a), and, replacing $X$ by an affine neighborhood of $\gamma (y)$,  that $X$ is affine, $X=\Spec\, B$.  Then $z:=\sigma (\gamma\gamma' )$ and $z':= \sigma (\gamma )\sigma (\gamma' )$ are $k$-points of $\Sym^n X =\Spec\, (B^{\otimes n})^{\Sigma_n}$. The latter algebra is generated by elements $b^{\otimes n}$ with $b\in B$  invertible at $\gamma (y)$. It remains to show that $b^{\otimes n} (z)=b^{\otimes n}(z')$. If $\gamma (y)= \Sigma n_i x_i$, then $b^{\otimes n} (z)= \prod  \text{Nm}_{k(x_i )/k} b(x_i )^{n_i }$ where $ \text{Nm}_{k(x_i )/k}:  k(x_i )\to k$ is the norm map. By Remark in (a), $b(z')$ is the determinant of the multiplication by $b$ action on the derived  fiber $Li^*_y \CA_Y$ of $ \CA_Y $ at $y$.  
 The latter is a complex of $B$-modules supported at $\{ x_i \}$ with the Euler characteristic at $x_i$ equal to $n_i$, and we are done.
\hfill$\square$

As an application, let us show that for any 
 commutative algebraic group $G$ over $k$ the presheaf $X\mapsto G(X)$ admits a natural transfer structure. Let $\pi : \Sym^\cdot G \to G$ be the morphism of monoid schemes
 which equals id$_G$ on $\Sym^1 G $. Now for 
  an effective $\gamma\in\CC or (Y,X)$  the  map $\gamma^* :G(X)\to G(Y)$   sends $a\in G(X)$ to the composition $Y\buildrel{a\gamma}\over\lra \Sym^\cdot G\buildrel\pi\over\to G$. In particular,
the sheaf $\CO^\times$ of invertible functions admits a natural transfer structure (given by the norm).

{\it Remark.} Consider an additive  category whose objects are  quasi-projective schemes and morphisms are group completions of the morphisms between the monoid schemes $\Sym^\cdot X$. The proposition shows that the category of smooth quasi-projective varieties and finite correspondences embeds into this category (and the embedding is fully faithful if $char\, k =0$).

2.1.4. Consider the DG category of complexes $C\CP\CS h_{\tr}$,  which is the same as the DG category $\CC or (\CS m)^{\text{op}} \text{-dgm}$ of DG $\CC or (\CS m)^{\text{op}}$-modules; it
contains the full DG category $\mathop{\vtop{\ialign{#\crcr
  \hfil\rm $\CP_{\tr}$\hfil\crcr
  \noalign{\nointerlineskip}\rightarrowfill\crcr
  \noalign{\nointerlineskip}\crcr}}}$
 of semi-free $\CC or (\CS m)^{\text{op}}$-modules (see 1.6). As in Remark (ii) in 1.7, let $\CP_{\tr} := (
 \mathop{\vtop{\ialign{#\crcr
  \hfil\rm $\CP_{\tr}$\hfil\crcr
  \noalign{\nointerlineskip}\rightarrowfill\crcr
  \noalign{\nointerlineskip}\crcr}}}
  )^{\text{perf}}$ be the homotopy idempotent completion of $\CC or (\CS m)^{\text{pretr}}$, so  $\CP^{\text{tri}}_{\tr}=(\CC or(\CS m)^{\text{tri}})^{\kappa} $. By (1.6.1),  $\mathop{\vtop{\ialign{#\crcr
  \hfil\rm $\CP_{\tr}$\hfil\crcr
  \noalign{\nointerlineskip}\rightarrowfill\crcr
  \noalign{\nointerlineskip}\crcr}}}$ is a DG structure on $D\CP\CS h_{\tr}$:  $$
  \mathop{\vtop{\ialign{#\crcr
  \hfil\rm $\CP^{\text{tri}}_{\tr}$\hfil\crcr
  \noalign{\nointerlineskip}\rightarrowfill\crcr
  \noalign{\nointerlineskip}\crcr}}} \iso D\CP\CS h_{\tr}. \tag 2.1.3$$

{\it Remarks.} \! (i)
 The   embedding $R [\CS m ] \hra\CC or (\CS m)$, $R[X]\mapsto R_{\tr}[X]$,  yields DG functors $\CP \to \CP_{\tr}$, $\LP \to \mathop{\vtop{\ialign{#\crcr
  \hfil\rm $\CP_{\tr}$\hfil\crcr
  \noalign{\nointerlineskip}\rightarrowfill\crcr
  \noalign{\nointerlineskip}\crcr}}}$. The corresponding  functor ${  \LP}^{\text{tri}} \to  \mathop{\vtop{\ialign{#\crcr
  \hfil\rm $\CP^{\text{tri}}_{\tr}$\hfil\crcr
  \noalign{\nointerlineskip}\rightarrowfill\crcr
  \noalign{\nointerlineskip}\crcr}}}$ is left adjoint, via identifications (2.1.1) and (2.1.3), to the forgetful functor $o: D\CP\CS h_{\tr} \to D\CP\CS h$, so we have the adjoint pair  $$D\CP\CS h \rightleftarrows D\CP\CS h_{\tr} . \tag 2.1.4$$

(ii) If $K/k$ is any extension, then the base change functor $\CS m_k \to \CS m_K$, $X\mapsto X_K$, extends naturally to  $\CC or (\CS m)$ and the DG categories $\mathop{\vtop{\ialign{#\crcr
  \hfil\rm $\CP_{\tr}$\hfil\crcr
  \noalign{\nointerlineskip}\rightarrowfill\crcr
  \noalign{\nointerlineskip}\crcr}}}$,  $\CP_{\tr}$.  If $K/k$ is separable, then we have a similar exact functor for the categories of presheaves with transfers (see Exercise (ii) in 2.1.2), and   (2.1.3) is compatible with these functors.

Voevodsky's notation for $\CP\CS h_{\tr}$ is $Pre Shv (Sm Cor (k))$ in \cite{Vo2}, and {\bf PST} in \cite{MVW}.
 
\medskip {\bf 2.2. The tensor product and inner $\CH om$; the Tate twist.} 
The direct product of smooth varieties defines on  $\CC or(\CS m)$ a  tensor category structure,   $\Z [X]\otimes\Z [Y] =\Z [X\times Y]$, with the unit object $\Z $. Thus $\CP_{\tr}$ and $\mathop{\vtop{\ialign{#\crcr
  \hfil\rm $\CP_{\tr}$\hfil\crcr
  \noalign{\nointerlineskip}\rightarrowfill\crcr
  \noalign{\nointerlineskip}\crcr}}}$ are DG tensor  categories (see 1.9.2), and $ D\CP\CS h_{\tr} $ is a tensor triangulated category (via (2.1.3)). The tensor product is right t-exact for
 the evident t-structure on $ D\CP\CS h_{\tr} $.

As in 1.9.3, the DG tensor  category $\CP_{\tr}^{\text{op}}$ acts naturally on $C\CP\CS h_{\tr}$. The action of $P\in \CP_{\tr}$ is denoted by $F\mapsto \CH om (P,F)$, so for $P=R_{\tr}[X]$ one has $\CH om (R_{\tr}[X], F)(Y)= F(X\times Y)$. Passing to the homotopy categories, we get the action of the  tensor triangulated category $\CP_{\tr}^{\text{op\, tri}}$ on $D\CP\CS h_{\tr}$. The endofunctors  
 $F\mapsto P\otimes F$ and  $F\mapsto \CH om (P,F)$ of  $D\CP\CS h_{\tr}
= \mathop{\vtop{\ialign{#\crcr
  \hfil\rm $\CP_{\tr}$\hfil\crcr
  \noalign{\nointerlineskip}\rightarrowfill\crcr
  \noalign{\nointerlineskip}\crcr}}}^{\text{tri}} $ are  adjoint.

 By 1.10, the above picture remains valid if we replace   $\CC or (\CS m)$  by $R [\CS m ]$. 
 Since $R [\CS m ] \hra \CC or (\CS m) $ is a tensor functor,    $\CP \to \CP_{\tr}$,  $\LP \to \mathop{\vtop{\ialign{#\crcr
  \hfil\rm $\CP_{\tr}$\hfil\crcr
  \noalign{\nointerlineskip}\rightarrowfill\crcr
  \noalign{\nointerlineskip}\crcr}}}$, and $D\CP\CS h \to D\CP\CS h_{\tr}$ are also tensor functors. The forgetful functor $o:  C\CP\CS h_{\tr} \to C\CP\CS h$ is compatible with the $\CP^{\text{op}}$-action. 
   
{\it Remark.}  Since $o: D\CP\CS h_{\tr} \to D\CP\CS h$ is
  right adjoint to a tensor functor ${ \LP}^{\text{tri}} \to  \mathop{\vtop{\ialign{#\crcr
  \hfil\rm $\CP_{\tr}$\hfil\crcr
  \noalign{\nointerlineskip}\rightarrowfill\crcr
  \noalign{\nointerlineskip}\crcr}}}^{\text{tri}} $, it acts naturally on polylinear operations (i.e., in the terminology of \cite{BD} 1.1.5, $o$ is a {\it pseudo-tensor} functor). Explicitly, for 
a finite collection of objects  $F_i , G\in D\CP\CS h_{\tr}$ one has a canonical map $\Hom (\otimes F_i ,G)\to \Hom (\otimes o(F_i ), o(G))$ which assigns to $\phi : \otimes F_i \to G$ the composition $\otimes o (F_i  )\to o(\otimes F_i )\buildrel{o(\phi )}\over\lra o(G)$, where the first arrow comes from the evident maps $\otimes\, \CC or (Z,X_i )\to \CC or (Z, \prod X_i )$ (the latter are not isomorphisms in general, so $o$ is {\it not} a tensor functor). Therefore $o$ transforms any kind of  algebras in $D\CP\CS h_{\tr} $ to algebras of similar kind in $ D\CP\CS h$.

Set $R (1):= R [\Bbb G_m /\Spec \, k][-1]\in \CP $ where $\Bbb G_m := \Bbb A^1 \smallsetminus \{ 0\}$.\footnote{At the moment, we do not care about the group structure on $\Bbb G_m$. Our notation is {\it not} consistent with that of \cite{MVW} and \cite{Vo2}: in loc.~cit., $\Bbb G_m$ denotes the {\it pointed} variety $(\Bbb A^1 \smallsetminus \{ 0\}, 1)$.} We denote the corresponding object of $\CP_{\tr}$ by $\Z (1)$. As in 2.1.1,  every $a\in \Bbb G_m (k)$ yields a morphism $R\to R [\Bbb G_m ]$ so that 
 $R \oplus R (1)[1] \iso R [\Bbb G_m ]$; the standard choice is $a=1$.
 The  tensor product with $R (1)$ or $\Z (1)$ in $\CP$, $\CP_{\tr}$ and related tensor categories is denoted by $F\mapsto F(1)$; this is the {\it Tate twist} endofunctor.

For  $F\in\CP\CS h$ we define a presheaf $F_{-1}$ by formula $F_{-1}(X):= F((X\times\Bbb G_{m })/X)= \Coker (F(X)\to F(X\times \Bbb G_{m }))$ where the arrow comes from the projection $X\times\Bbb G_{m }\to X$. Same construction works  in the case of presheaves with transfers. We can also pass to complexes;  one has 
 $F_{-1}= \CH om (R (1)[1],F)$. The endofunctor  $F\mapsto F_{-1}[1]$ of $D\CP\CS h_{\tr}  $ is right adjoint to the Tate twist.

The  iterations of the above functors are denoted by $F\mapsto F(n), F_{-n}$, $n\in \Bbb Z_{\ge 0}$.

\medskip

{\bf 2.3. The  DG category of effective motives.}  Consider the following two types of complexes in $\CC or (\CS m)$:

($\Delta$) The 2-term complexes $\Z [ \Bbb A^1 \times X]\to \Z [X]$, the differential comes from the projection $\Bbb A^1 \times X \to X$, $X\in\CS m$.

(MV$_{\text{Zar}}$)  The 3-term Mayer-Vietoris ones $\Z [U\cap V]\to \Z [U]\oplus \Z [V]\to\Z [X]$. Here $U,V\subset X$ are Zariski open subsets such that $U\cup V=X$, 
and the differentials are, respectively, the difference and the sum of the  embedding morphisms.

Let $\CI^{\Delta}_{\tr},\CI^{\text{Zar}}_{\tr} ,\CI^{\Delta \text{Zar}}_{\tr} \subset \CP_{\tr}$ be the 
homotopy idempotent completions of the full pretriangulated 
 DG subcategories strongly generated by the objects of the corresponding types (see 1.6.2). Set (see 1.8) $$\CD_\CM^{\text{eff}}= \CD_{\CM /k}^{\text{eff}}  := \CP_{\tr} /^\kappa \, \CI^{\Delta \text{Zar}}_{\tr}.\tag 2.3.1$$  This is the  DG category of {\it effective geometric motives}. The homotopy object of $\CD_\CM^{\text{eff}}$ that corresponds to  $\Z [X]$  is denoted by $M(X)$; this is the {\it motive} of $X$.

Passing to the  ``infinite" categories, we see, by (1.8.1), that $$
\mathop{\vtop{\ialign{#\crcr
  \hfil\rm $\CD^{\text{eff}}_{\CM}$\hfil\crcr
  \noalign{\nointerlineskip}\rightarrowfill\crcr
  \noalign{\nointerlineskip}\crcr}}} 
 \iso \mathop{\vtop{\ialign{#\crcr
  \hfil\rm $\CP_{\tr}$\hfil\crcr
  \noalign{\nointerlineskip}\rightarrowfill\crcr
  \noalign{\nointerlineskip}\crcr}}} /\mathop{\vtop{\ialign{#\crcr
  \hfil\rm $\CI^{\Delta}_{\tr}$\hfil\crcr
  \noalign{\nointerlineskip}\rightarrowfill\crcr
  \noalign{\nointerlineskip}\crcr}}}{\!}^{ \text{Zar}} . \tag 2.3.2$$ This is the  DG category of {\it effective  motives}.

Consider the corresponding homotopy categories. We use identification (2.1.3). Set 
 $I^{\text{Zar}}_{\tr}:= \mathop{\vtop{\ialign{#\crcr
  \hfil\rm $\CI^{\text{Zar}}_{\tr}$\hfil\crcr
  \noalign{\nointerlineskip}\rightarrowfill\crcr
  \noalign{\nointerlineskip}\crcr}}}{}^{\text{tri}}$, etc; these are right-admissible subcategories of $D\CP\CS h_{\tr}$, and $$\mathop{\vtop{\ialign{#\crcr
  \hfil\rm $\CD^{\text{eff}}_{\CM}$\hfil\crcr
  \noalign{\nointerlineskip}\rightarrowfill\crcr
  \noalign{\nointerlineskip}\crcr}}}{\! }^{\text{tri}} \iso D\CP\CS h_{\tr}  /I^{ \Delta  \text{Zar}}_{\tr} \buildrel\sim\over\leftarrow (I^{ \Delta \text{Zar}}_{\tr} )^\bot  , \tag 2.3.3$$   $$  \CD_\CM^{\text{eff}}{}^{\text{tri}}\iso (D\CP\CS h_{\tr}  /I^{ \Delta \text{Zar}}_{\tr})^{\text{perf}}
. \tag 2.3.4$$ 

The right $I^{ \Delta \text{Zar}}_{\tr}$-localization endofunctor $\CC^\CM : D\CP\CS h_{\tr} \to D\CP\CS h_{\tr}$ is called the {\it motivic localization}; it will be studied in the next two sections.

Notice that  the tensor product by any object of $\CC or (\CS m)$ preserves the complexes of type
 ($\Delta$) and of type  (Zar).  As in 1.9,  the above $\CI^?_{\tr}$'s are ideals in $\CP_{\tr}$, and $\CD_\CM^{\text{eff}}$ is a homotopy tensor DG category with the unit object $R_\CM := M(\Spec \, k )$. Same is true for the  ``infinite" categories.  

The image in $ \CD_\CM^{\text{eff}}$ of $\Z (1) \in \CP_{\tr}$  is  the {\it Tate motive}  $R_\CM (1)$.

Voevodsky's notation for $\CD_\CM^{\text{eff}}{}^{\text{tri}}=  \CD_{\CM /k}^{\text{eff\, tri}}$ is   $DM^{\text{eff}}_{gm} (k)$ in \cite{Vo2}. 

{\it Remarks.} (i) In \S 4 we will see how to perform the motivic localization explicitly if
 the base field $k$ is perfect.  It is not clear if the above definition is reasonable for non perfect $k$ (see e.g.~Remark (b) in 4.4). 

(ii) The embedding $ \mathop{\vtop{\ialign{#\crcr
  \hfil\rm $\CD^{\text{eff}}_{\CM}$\hfil\crcr
  \noalign{\nointerlineskip}\rightarrowfill\crcr
  \noalign{\nointerlineskip}\crcr}}}{}^{\text{tri}} \iso (I^{ \Delta \text{Zar}}_{\tr} )^\bot \subset D\CP\CS h_{\tr} $ is a pseudo-tensor functor
(as the right adjoint to a tensor functor, cf.~Remark in 2.2). So $\CC^\CM :D\CP\CS h_{\tr} \to D\CP\CS h_{\tr} $ is a pseudo-tensor functor as well.

\medskip

{\bf 2.4. The Artin motives.} Let $\CS m_0 $ be the category of smooth zero-dimensional varieties.  The embedding  $\CS m_0 \hra \CS m$ admits a left adjoint functor $\pi_0 : \CS m\to \CS m_0$, where $\pi_0 (X)$ is the spectrum of the integral closure of $k$ in $\CO (X)$. The corresponding full subcategory $\CC or(\CS m_0 )\subset \CC or(\CS m)$ is an additive $R$-category, and $\pi_0$ extends naturally to an $R$-linear functor $\pi_0 : \CC or(\CS m)\to \CC or(\CS m_0 )$ left adjoint to the embedding.

Pick a separable closure $\bar{k}/k$ of $k$; let $G$ be the Galois group. By Galois theory, the functor $X\mapsto X(\bar{k})$ identifies $\CS m_0 $ with the category of finite $G$-sets. This implies easily that the functor $R_{\tr} [X] \mapsto R[X(\bar{k})]$ is a fully faithful embedding of $\CC or(\CS m_0 )$ into the category of $R[G]$-modules; its essential image is the category $R[G]$-prm of  permutative  $R[G]$-modules.

The DG category $\CA r =\CA r_k$ of {\it Artin motives} is defined as the homotopy idempotent completion of   $\CC or(\CS m_0 )$ inside $\CP_{\text{tr}}$ (see 1.6.2, 1.7). So the triangulated category $\CA r^{\text{tri}}$ is the idempotent completion of the homotopy category of finite complexes of objects of $\CC or(\CS m_0 )$, i.e.,  $\CA r^{\text{tri}}=(K^b R[G]$-prm$)^\kappa$. The  Artin motives form
 a full subcategory in $\CD_{\CM}^{\text{eff\, tri}} $:

\proclaim{\quad Lemma}  The
functor $\pi_0^{\text{tri}}$ factors through a triangulated functor $\CD_\CM^{\text{eff\, tri}} \to \CA r^{\text{tri}}$, which is left adjoint to the composition $\CA r^{\text{tri}} \hra (\CC or(\CS m )^{\text{tri}})^\kappa \to \CD_{\CM}^{\text{eff\, tri}} $. The latter composition is a fully faithful embedding $\CA r^{\text{tri}} \hra  \CD_{\CM}^{\text{eff\, tri}} $.
\endproclaim

{\it Proof.} To see that $\pi_0^{\text{tri}} : \CC or(\CS m )^{\text{tri}} \to \CA r^{\text{tri}}$ factors through $\CD_\CM^{\text{eff\, tri}}$, we need to check that $\pi_0^{\text{tri}}$ kills complexes ($\Delta$) and (MV$_{\text{Zar}}$) from  2.3. This is evident since  $\pi_0 ( \Bbb A^1 \times X )= \pi_0 (X)$ and $\pi_0 (U_\cdot )\to \pi_0 (X)$ is a homotopy equivalence. The rest is immediate. \hfill$\square$

{\it Remark.} Let $G$ be any (pro)finite group; suppose that $R$ is noetherian regular. 
Then, according to \cite{R}, 
 the evident functor $K^b R[G]$-prm $\to D^b R[G]$-mod (the bounded derived category of finitely generated $R[G]$-modules) yields an equivalence $(K^b R[G]$-prm$/I)^\kappa \iso D^b R[G]$-mod;
here $I$ is the thick subcategory of acyclic complexes of permutative modules.

 \bigskip
 
\centerline{\bf \S 3. Premotives and the $\Bbb A^1$-homotopy localization }

\medskip

In this section  we  consider the  $\Bbb A^1$-homotopy localization procedure used in the definition of the category of motives more thoroughly. The material corresponds to lectures 4 and 9 of \cite{MVW} and \cite{Vo3}. The constructions (if not the terminology) of 3.1  go back to \cite{S}. The principal result is theorem 3.3 (Voevodsky's cancellation theorem);  its proof is a variant of Voevodsky's  original one \cite{Vo3}. The construction  in 3.4.1  is used to define the residue maps in 4.6; the theme of 3.4.2 is developed in 5.1.4.

{\bf 3.1. The $\Bbb A^1$-homotopy localization.} 3.1.1.  Let $I^\Delta \subset D\CP\CS h$  be the  triangulated subcategory generated by all complexes of type $R [ \Bbb A^1 \times X/X] $  where $X\in\CS m$ and $ \Bbb A^1 \times X \to X$ is the projection. Its objects are {\it $\Bbb A^1$-homotopy contractible} complexes. One has the corresponding subcategory $I^{\Delta\bot} \subset D\CP\CS h$ 
 of {\it $\Bbb A^1$-homotopy invariant}, or  {\it $\Bbb A^1$-homotopy local}, complexes. 
  A complex of presheaves $F$ lies in $I^{\Delta\bot}$ if and only if  for every $X$ the morphism $F(X)\to F( \Bbb A^1 \times X)$ is a quasi-isomorphism, i.e., the morphism $F\to \CH om (R[\Bbb A^1 ], F)$ is a quasi-isomorphism in $D\CP\CS h$. $F$ is $\Bbb A^1$-homotopy invariant if and only if 
 all the cohomology $H^a F$ are   $\Bbb A^1$-homotopy invariant. Thus $I^{\Delta\bot}$ carries a  non-degenerate t-structure whose core is the abelian  category $\CP\CS h^\Delta$ of $\Bbb A^1$-homotopy invariant presheaves. 
  
  {\it Exercise.}  $F\in D\CP\CS h$ is $\Bbb A^1$-homotopy invariant iff  the   pull-back maps $i_0^* , i_1^* : \CH om (R[\Bbb A^1] ,F)\rightrightarrows F$ for the embeddings $i_a :\Spec\, k\hra \Bbb A^1 $, $\cdot \mapsto a\in\Bbb A^1 $, coincide.

Our $I^\Delta$ is right admissible by 1.4.2, so $I^{\Delta\bot}\iso D\CP\CS h /I^\Delta$, and the  t-structure picture fits into the format of part (ii) of the lemma in 1.2.

In fact, the $I^\Delta$-localization functor lifts to a DG 
endofunctor $\CC^\Delta$  of $ C\CP\CS h$ defined as follows.
Let $\Delta^\cdot$ be the standard  cosimplicial scheme of the ``simplices"   $\Delta^n := \{ (t_0 ,\ldots ,t_n )\in \Bbb A^{n+1}: \sum t_i =1\}$. A presheaf $F$ yields then a simplicial presheaf $\CH om (R[\Delta^\cdot ],F): X\mapsto F(X\times \Delta^\cdot )$; let $C^\Delta (F)$ be its normalized complex. Explicitly,   $C^\Delta (F)^n (X):= \{ f\in F(X\times\Delta^n ): \partial_1 (f)=\ldots =\partial_n (f)=0\}$,  the differential is $f\mapsto \partial_0 (f)$.    For $F\in C\CP\CS h$ we denote by $C^\Delta (F)$  the total complex of the corresponding bicomplex (so   $\CC^\Delta (F)^n (X)= \mathop\sum\limits_{a-b =n} F^a (X\times\Delta^b )$).
 In the notation of 1.9 and 2.2, one has $\CC^\Delta (F)= \limright \CH om (\sigma_{\le n} R [\Delta^\cdot ]^{\text{nor}}, F)$; here $R [\Delta^\cdot ] $ is a cosimplicial object in  $R [\CS m ]$,  $R[\Delta^\cdot ]^{\text{nor}}$ its normalized complex, and $\sigma_{\le n} R[\Delta^\cdot ]^{\text{nor}} $ are its stupid truncations (which are finite complexes). A
natural embedding $F\hra \CC^\Delta (F)$ identifies $F$ with the component of the $b$-degree 0.
Since $\CC^\Delta$ preserves quasi-isomorphisms, it yields an endofunctor on $D\CP\CS h$. The latter  equals the $\Bbb A^1$-homotopy localization:

\proclaim{\quad Proposition} One has
 $ \CC^\Delta (F)\in I^{\Delta\bot},$ $\CC^\Delta (F)/F \in I^\Delta$.
 \endproclaim
 
{\it Proof.} We show that for a presheaf $F$ the complex $\CC^\Delta (F)$ is $\Bbb A^1$-homotopy invariant, leaving the rest for the reader. It suffices to construct a  simplicial homotopy $h_\cdot$ between the morphisms of simplicial $R$-modules $i_0^* ,i_1^* : F(\Bbb A^1 \times X\times\Delta^\cdot )\rightrightarrows F(X\times\Delta^\cdot )$.
Let $T_\cdot$ be the simplicial interval, i.e., the simplicial set with $T_n$ = the set of monotone maps $a: [0,n]\to [0,1]$. Our $h_\cdot$ is comprised by maps $h_n : F(\Bbb A^1 \times X\times\Delta^n)\otimes R[T_n ]\to F(X\times\Delta^n )$, i.e., by  those $h_n^a : F(\Bbb A^1 \times X\times\Delta^n)\to F(X\times\Delta^n )$, $a\in T_n$, $h_n (f\otimes a)=h_n^a (f).$ We define $h_n^a$ as the pull-back  for a map $X\times\Delta^n \to\Bbb A^1 \times X\times\Delta^n$, $(x, (t_0 ,\ldots ,t_n ))\mapsto (\sum a(i) t_i , x ,( t_0 ,\ldots ,t_n ))$. 
\hfill$\square$

{\it Exercise.} Suppose a presheaf $F$ admits a natural map $\kappa : F(X)\to F(\Bbb A^1 \times X)$ which preserves the zero section ($\kappa$ need not be additive) and satisfies $i_0^* \kappa =\id_F$, $i_1^* \kappa =0$. Then $\CC^\Delta (F)$ is acyclic, i.e., $F\in I^\Delta$.

3.1.2. The  story of 3.1.1 remains literally valid if we replace presheaves by presheaves with transfers. So  we have  subcategories
 $I^\Delta_{\tr}, I^{\Delta\bot}_{\tr}$ $ \subset D\CP\CS h_{\tr}$ of  $\Bbb A^1$-homotopy contractible and $\Bbb A^1$-homotopy invariant presheaves with transfers, a DG endofunctor $F\mapsto \CC^\Delta (F):= \limright \CH om (\sigma_{\le n} R_{\tr} [\Delta^\cdot ]^{\text{nor}}, F)$ of $C\CP\CS h_{\tr}$ which descends to the  $I^\Delta_{\tr}$-localization endofunctor of $D\CP\CS h_{\tr}$, a non-degenerate t-structure on  $ I^{ \Delta\bot}_{\tr}\iso D\CP\CS h_{\tr}  /I^{\Delta}_{\tr}$ whose core is the abelian category $\CP\CS h^\Delta_{\tr}$ of $\Bbb A^1$-homotopy invariant presheaves with transfers, etc. 

 {\it Exercise.} Show that $ I^{ \Delta\bot}_{\tr}$ is {\it not} equivalent to the derived category of
 $\CP\CS h^\Delta_{\tr}$.
 \footnote{Hint: Since $R_{\tr}[X]$ is a projective object of  $\CP\CS h_{\tr}$ and the embedding   $ I^{ \Delta\bot}_{\tr}\hra   D\CP\CS h_{\tr}$ is t-exact,  $H^0 \CC^\Delta (R_{\tr}[X] )$ is a projective object in  $\CP\CS h^\Delta_{\tr}$.   Were $D \CP\CS h^\Delta_{\tr} \to I^{ \Delta\bot}_{\tr}$ be an equivalence, this would imply that    $  \CC^\Delta (R_{\tr}[X] )\iso H^0 \CC^\Delta (R_{\tr}[X] ) $.}

  The $\Bbb A^1$-homotopy localization is compatible with functors (2.1.4). Namely, the functor $D\CP\CS h \to D\CP\CS h_{\tr}$ sends $I^\Delta$ to $I^\Delta_{\tr}$;  its right adjoint (the forgetful functor) $D\CP\CS h_{\tr} \to D\CP\CS h$  commutes with $\CC^\Delta$ (hence sends  $I^\Delta_{\tr}$ to $I^\Delta$,  $I^{\Delta\bot}_{\tr}$ to $I^{\Delta\bot}$).

 As in 2.3, our triangulated categories carry natural DG structures. Namely, set $ \CD_{p\CM}^{\text{eff}}:=  \CP_{\tr} /^\kappa \CI^{\Delta}_{\tr}$. Then, by (1.8.1), one has  $\mathop{\vtop{\ialign{#\crcr
  \hfil\rm $\CD^{\text{eff}}_{p\CM}$\hfil\crcr
  \noalign{\nointerlineskip}\rightarrowfill\crcr
  \noalign{\nointerlineskip}\crcr}}}  \iso \mathop{\vtop{\ialign{#\crcr
  \hfil\rm $\CP_{\tr}$\hfil\crcr
  \noalign{\nointerlineskip}\rightarrowfill\crcr
  \noalign{\nointerlineskip}\crcr}}}  /\mathop{\vtop{\ialign{#\crcr
  \hfil\rm $\CI^\Delta_{\tr}$\hfil\crcr
  \noalign{\nointerlineskip}\rightarrowfill\crcr
  \noalign{\nointerlineskip}\crcr}}}  $ and $$\mathop{\vtop{\ialign{#\crcr
  \hfil\rm $\CD^{\text{eff}}_{p\CM}$\hfil\crcr
  \noalign{\nointerlineskip}\rightarrowfill\crcr
  \noalign{\nointerlineskip}\crcr}}}{\!}^{\text{tri}} \iso  D\CP\CS h_{\tr}  /I^{\Delta}_{\tr}  \buildrel\sim\over\leftarrow I^{ \Delta\bot}_{\tr} .\tag 3.1.1 $$

 Since $\CI^\Delta_{\tr}$ is a tensor ideal in $\CP_{\tr}$, our $\CD_{p\CM}^{\text{eff}}$ and $\mathop{\vtop{\ialign{#\crcr
  \hfil\rm $\CD^{\text{eff}}_{p\CM}$\hfil\crcr
  \noalign{\nointerlineskip}\rightarrowfill\crcr
  \noalign{\nointerlineskip}\crcr}}} $ are homotopy tensor DG categories.
The endofunctor $\CC^\Delta$ of $D\CP\CS h_{\tr}$ is right t-exact. Hence the tensor product on $\mathop{\vtop{\ialign{#\crcr
  \hfil\rm $\CD^{\text{eff}}_{p\CM}$\hfil\crcr
  \noalign{\nointerlineskip}\rightarrowfill\crcr
  \noalign{\nointerlineskip}\crcr}}}{\!}^{\text{tri}}$ is right t-exact.

Objects of $\mathop{\vtop{\ialign{#\crcr
  \hfil\rm $\CD^{\text{eff}}_{p\CM}$\hfil\crcr
  \noalign{\nointerlineskip}\rightarrowfill\crcr
  \noalign{\nointerlineskip}\crcr}}} $  (resp.~$\CD_{ p \CM}^{\text{eff}}$) are called {\it  effective (geometric) premotives}. The image $\Z^\Delta [X]$ of $\Z [X]$ in $\CD_{p\CM}^{\text{eff}}$ is called the {\it premotive} of $X$. We   have the Tate premotive $\Z^\Delta (1)$, etc.

{\it Remark.}   Let $\Ho\, \CC or^\Delta (\CS m )$ be the full subcategory of $ \CD_{p \CM}^{\text{eff\, tri}}$ formed by objects $\Z^\Delta [X]$, $X\in \CS m$, i.e., its morphisms $\CC or^\Delta (Y,X) :=H^0 \Hom (\Z^\Delta [Y],\Z^\Delta [X])$ are $\Bbb A^1$-homotopy classes of correspondences. Then $\Bbb A^1$-homotopy invariant presheaves with transfers are the same as right $\Ho\Z^\Delta [\CS m]$-modules, i.e., $R$-linear functors $F: \Ho\Z^\Delta [\CS m]^{{\text{op}}} \to R$-mod.

3.1.3. Let $X$, $Y$ be affine smooth varieties.

{\it Question.} Is it true that the complex  $ \Hom (\Z^\Delta [X],\Z^\Delta [Y])= \CC or (X\times\Delta^\cdot ,Y)$  is acyclic in degrees $<-\dim Y$?

\medskip

{\bf 3.2. The Tate premotive.}  Here are some typical computations in $\CD_{p\CM}^{\text{eff}}$:

 \proclaim{\quad Proposition} (i) The $\Bbb Z$-presheaf with transfers $\CO^\times $ (see 2.1.3)  is $\Bbb A^1$-homotopy invariant and there is a canonical homotopy equivalence in $\CD_{p\CM}^{\text{eff}}$ $$ \epsilon_R  : \Z^\Delta (1)[1]\iso \CO^\times \buildrel{L}\over\otimes R. \tag 3.2.1$$

(ii) The transposition of factors action of the symmetric group $\Sigma_n$ on  the premotive $\Z^\Delta (1)^{\otimes n} =\Z^\Delta (n)$ is homotopically trivial (i.e., is trivial in $\CD_{p\CM}^{\text{eff\, tri}}$).
\endproclaim

{\it Proof.} (i) The $\Bbb A^1$-homotopy invariance of $\CO^\times $ is clear (our varieties are reduced). 

Let $\epsilon :   \Bbb Z_{\tr}(1)[1]\to \CO^\times  $ be a morphism of $\Bbb Z$-presheaves with transfers
defined as the composition    $  \Bbb Z_{\tr}(1)[1]\hra  \Bbb Z_{\tr} [\Bbb G_m ]\to \CO^\times $, the latter arrow sends $\gamma \in \CC or (X,\Bbb G_m )$ to $\gamma^* (t)\in  \CO^\times (X)$ (see 2.1.3); here $t\in \CO^\times (\Bbb G_m )$ is the standard parameter. Our $\epsilon$ is evidently surjective;  we will check in a moment that its kernel $K$ is in $I^\Delta$.  Since $\Bbb Z$-presheaves $K\subset \Bbb Z_{\tr }(1)[1] \subset \Bbb Z_{\tr} [\Bbb G_m ]$ have no $\Bbb Z$-torsion,  $\epsilon$ yields a quasi-isomorphism $\CC one (K\otimes R \to R_{\tr} (1)[1]) \to  \CO^\times \buildrel{L}\over\otimes R$. Then $\epsilon_R$ is its composition with the embedding $\Z (1)[1] \to \CC one (K\otimes R \to \Z (1)[1])$.
 
To show that $K\in I^\Delta$, we present $\epsilon$ as follows. For $X\in \CS m$ a morphism $ \Bbb Z_{\tr} [X] \to \Bbb Z_{\tr} (1)[1]$ is  a divisor  $D$ on $\Bbb G_{mX}$ finite over $X$ of relative degree 0. 
Let $\Phi (X)$ be the multiplicative group of meromorphic functions $\phi$ on $\Bbb P^1_X$ such that $\phi$ is defined on a Zariski neighborhood of $\{ 0,\infty \} \times X$, $\phi (\infty ,x)\equiv 1$, and $\phi (0,x)$ is invertible. One has an isomorphism $\Phi (X) \iso \Hom ( \Bbb Z_{\tr} [X] , \Bbb Z_{\tr }(1)[1])$, $\phi \mapsto div\,\phi$, and $\epsilon (div\,\phi )(x)=\phi (0,x)$. Therefore $K$ equals the subgroup of $\phi \in \Phi$ with $\phi (0,x)\equiv 1$. One has a natural map $\kappa :   K(X)\to K(\Bbb A^1 \times X )$, $ \kappa (\phi )(t, u,x):= u + (1-u )\phi (t,x)$ (here $u$ is the parameter on $\Bbb A^1$), which satisfies the conditions of the second exercise in 3.1.1, and we are done by that exercise.

(ii) It suffices to consider the case $n=2$. Let $\sigma$ be the transposition of factors of $\Bbb G_m \times \Bbb G_m $. We want to show that $\sigma +\id$ kills the direct summand $(R_{\tr}^\Delta (1)[1])^{\otimes 2}\subset \Z^\Delta [\Bbb G_m \times \Bbb G_m ]$ in the homotopy category, i.e.,
 that $\sigma + \id \in \CC or (\Bbb G_m \times\Bbb G_m , \Bbb G_m \times\Bbb G_m )$ is $\Bbb A^1$-homotopic to the sum of correspondences that take values in $1\times\Bbb G_m $ and in $\Bbb G_m \times 1$. Let $\alpha \in \End (\Bbb G_m \times\Bbb G_m )$ be the map $(x,x')\mapsto (xx',1)$; 
we check that $\sigma + \id$ is $\Bbb A^1$-homotopic to $\sigma\alpha +\alpha$. One has  $\sigma + \id =  p^\vee p$, where $p: \Bbb G_m \times\Bbb G_m \to \Sym^2 \Bbb G_m $ is the projection, $p^\vee$ is the transposed correspondence. Since $ \Sym^2 \Bbb G_m \iso \Bbb G_m \times \Bbb A^1$, $p(x,x')\mapsto (xx', x+x')$, our $p$ is $\Bbb A^1$-homotopic to $p\alpha$. Thus $\sigma + \id =  p^\vee p$ is $\Bbb A^1$-homotopic to $p^\vee \! p\alpha = \sigma\alpha +\alpha$, q.e.d. \hfill$\square$

\medskip 

{\it Question.} By (ii) above, $R^\Delta_{\tr} (* ):= \mathop\oplus\limits_{n\ge 0} R^\Delta_{\tr}(n)$ is  a commutative unital algebra in the tensor category $\mathop{\vtop{\ialign{#\crcr
  \hfil\rm $\CD^{\text{eff}}_{p\CM}$\hfil\crcr  \noalign{\nointerlineskip}\rightarrowfill\crcr
  \noalign{\nointerlineskip}\crcr}}}{\!}^{\text{tri}} $. Can one realize it naturally as an $E_\infty$-algebra in the homotopy tensor DG category   $\mathop{\vtop{\ialign{#\crcr
  \hfil\rm $\CD^{\text{eff}}_{p\CM}$\hfil\crcr  \noalign{\nointerlineskip}\rightarrowfill\crcr
  \noalign{\nointerlineskip}\crcr}}} $?

{\it Remark.} The next generalization of the proposition can be proved by a simple modification of the argument. Let $C$ be a connected smooth projective curve over $k$, $A\subset C$ be a  finite non-empty subscheme \'etale over $k$. 

(i) One has a canonical homotopy equivalence  $R^\Delta_{\tr } [C\smallsetminus A] \iso Pic (C,A)$  in $\CD_{p\CM}^{\text{eff}}$. Here  $Pic (C,A)$ is the generalized Jacobian (the group scheme of pairs $(\CL ,\iota )$, $\CL$ is a line bundle on $C$, $\iota$ is a trivialization of $\CL|_A$) considered as a presheaf with transfers (see 2.1.3); notice that  this presheaf is $\Bbb A^1$-homotopy invariant. 

(ii) (M.~Nori, cf.~\cite{M} 5.8). Let  $B \subset C\smallsetminus A$ be another non-empty finite subscheme \'etale over $k$. Set $M (C; A,B):= \CC one ( \Bbb Z^\Delta_{\tr} [B] \to \Bbb Z^\Delta_{\tr} [C\smallsetminus A]) \in  \CD_{\CM}^{\text{eff}}$. Consider the action of the symmetric group $\Sigma_n $ on $M(C;A,B)^{\otimes n}$; let $\tau_n \in \End (M(C;A,B)^{\otimes n})$ be the  $\Sigma_n$-averaging (the sum of operators $s$, $s\in\Sigma_n$). Then $\tau_n$ is naturally homotopic to 0 if $n> 2g-2 +|A|+|B|$; here $g$ is the genus of $C$.

\proclaim{\quad 3.3. Theorem} The Tate premotive is homotopically quasi-invertible:  
 on $\CD_{p\CM}^{\text{eff}}$ (or $\mathop{\vtop{\ialign{#\crcr
  \hfil\rm $\CD^{\text{eff}}_{p\CM}$\hfil\crcr
  \noalign{\nointerlineskip}\rightarrowfill\crcr
  \noalign{\nointerlineskip}\crcr}}} $) the Tate twist  endofunctor $M\mapsto M(1)$  is homotopically fully faithful.
 \endproclaim
 
{\it Remark.} For $F\in D\CP\CS h_{\tr}$ let $\nu_F : F\to (\CC^\Delta (F(1)))_{-1}[-1]$ be obtained   by adjunction (see 2.2) from the canonical morphism $F(1)\to \CC^\Delta (F(1))$. The theorem means that $\nu_F$ is a quasi-isomorphism if $F$ is $\Bbb A^1$-homotopically invariant.

 {\it Proof of the theorem.} We want to show that the Tate twist map
$\Hom ( M',M)\to \Hom (M'(1),M(1))$ is a quasi-isomorphism  for every premotives $M$, $M'$.
It suffices to consider the case when  $M$, $M'$ are the premotives of some $Z,Z' \in\CS m$. We can assume that $R=\Bbb Z$.  Notice that for $Z=\Spec\, k$  the statement follows from 3.2(i).

 As in 3.1,  we represent $ \Hom (\Bbb Z_{\text{tr}}^\Delta [Z'] ,\Bbb Z_{\text{tr}}^\Delta [Z] )$  by the complex $C=C(Z',Z)$ with components  $C^a = \CC or (Z'\times\Delta^a ,Z)$. Set $C^\flat=C^\flat(Z',Z):= C(Z'\times\Bbb G_m ,Z\times\Bbb G_m )$, $E:= \CC or (\Bbb G_m ,\Bbb G_m )=\End (\Bbb Z_{\tr} [\Bbb G_m ])$; our $C^\flat$ is an $E$-bimodule. Let $\tau\in E$ be the correspondence $x\mapsto (x)-(1)$; this is a projector onto $\Bbb Z_{\tr} (1)[1]\subset \Bbb Z_{\tr} [\Bbb G_m ]$. Let $C^\flat_\tau $ be the direct summand $\tau C^\flat \tau $ of $C^\flat$,
 $E_\tau := \tau E\tau = \End (\Bbb Z_{\tr} (1))$, so $C^\flat_\tau $ is a $E_\tau$-bimodule.
  The Tate twist map is represented as $\rho : C\to C^\flat_\tau$, $\gamma\mapsto \gamma\otimes \tau$. Let us show that $\rho$ is a quasi-isomorphism.

3.3.1. {\it $\rho$ is injective on cohomology:} 
For $n\ge 1$ let $\iota_n \in E$ be the correspondence $x\mapsto x^{1/n}$ (so $deg ( \iota_n )=n$); one has $\iota_{mn}=\iota_m \iota_n$. 
Notice that the image of $\tau \iota_n \tau =\iota_n \tau\in E_\tau$ in $H^0 \End (\Bbb Z_{\tr}^\Delta (1))=\Bbb Z$ equals 1,  so the left multiplication by $\iota_n$ preserves  $C^\flat_\tau $ and is homotopic to the identity on it.

We will define  a subcomplex $C^\natural \subset C^\flat $ preserved by the $\iota_n$,
and a morphism of complexes $\kappa : C^\natural \to C$ with the next properties: (a) For any $\chi \in C^\flat $ one has $\iota_n \chi \in C^\natural$ for $n\gg 1$; (b) For any  $\gamma \in C$ and $n>1$ one has $\iota_n \rho (\gamma)\in C^\natural$ and   $\kappa (\iota_n \rho (\gamma))=\gamma$. This yields the injectivity of $H^\cdot \rho$.

For a smooth $Y$ (we will need $Y=Z'\times \Delta^\cdot$) consider $D\in  \CC or (Y \times \Bbb G_m ,\Bbb G_m )$ given by an effective divisor.  Let $t, t' \in \CO^\times (\Bbb G_m \times (Y\times \Bbb G_m ))$ be the projections to the two $\Bbb G_m$-factors. Locally on $Y$ we can write  $D = (f =0)$ where $f= f(t,t')\in \CO_Y (t,t')$ is a Laurent polynomial. 
The condition on $D$ says that its closure  in $\Bbb P^1 \times (Y\times\Bbb G_m )$ does not meet the divisor $t=\{ 0,\infty\} $. This means that, after multiplication by an invertible function on $Y$ and  powers of $t$, $t'$, our $f$ can be written as $1 + p_0 (t,t')$ and as $1 + p_\infty (t^{-1},t^{\prime -1})$ where $p_0 , p_\infty \in t\CO_{Y}[t,t', t^{\prime -1}]$. We say that $D$ is {\it nice} if  $ p_0 , p_\infty $ are sums of monomials of degree $>0$. This condition assures that  $D\cap (t=t')\subset \Bbb G_{m} \times Y$ (here $\Bbb G_m$ is the diagonal $t=t'$) is closed in $\Bbb P^1 \times Y$, hence it is finite over $Y$. 

Now our  $C^\natural \subset C^\flat $  is generated by those effective correspondences $\chi \in C^\flat$ that
 $\pi\chi \in \CC or (Z' \times\Delta^a \times \Bbb G_m ,\Bbb G_m )$ is  nice; here $\pi$ is the projection $Z\times\Bbb G_m \to \Bbb G_m$. By the above,  any $\chi \in C^\natural$ meets the divisor $t=t' $ on $(Z\times\Bbb G_m )\times(Z'\times\Bbb G_m )$ properly, and the cycle-theoretic intersection $\chi_{\text{diag}}:=\chi \cap (t=t')$ is finite over $Z'$. We define  $\kappa (\chi )\in C$ as the image of $-\chi_{\text{diag}}$ by the projection to $Z\times Z' \times\Delta^a$. Properties (a), (b) are evident.

3.3.2. {\it   $\rho   $ is surjective on cohomology:} Let $\chi \in C^\flat_\tau$ be any cycle.  Let us show that 
for sufficiently large $n$ our $\chi $ is homotopic to $ \rho \kappa (\iota _n \chi )$.

Consider another copy $\Bbb G^c_m$ of $\Bbb G_m$; let $\tau^c$ be the projector $\tau$ on it, etc. Let us perform the above constructions replacing $Z^{(\prime )}$ by $Z^{(\prime )}\times\Bbb G_m$ and $\Bbb G_m$ by $\Bbb G^c_m$. We get a subcomplex $C^\natural  (Z'\times\Bbb G_m ,Z\times \Bbb G_m )  \subset C^\flat (Z'\times\Bbb G_m ,Z\times \Bbb G_m )$, a morphism $\kappa^c : 
C^\natural  (Z'\times\Bbb G_m ,Z\times \Bbb G_m )\to C  (Z'\times\Bbb G_m ,Z\times \Bbb G_m )$, etc.

Let $\sigma$ be the transposition symmetry of $\Bbb G_m \times\Bbb G^c_m$. By 3.2(ii), one has 
$\chi \otimes\tau^c -\sigma (\chi \otimes\tau^c )\sigma  = d \phi$ for some  
$\phi \in C(Z' \times \Bbb G_m \times \Bbb G_m^c ,Z\times\Bbb G_m \times\Bbb G_m^c )$ which is
invariant with respect to the left and right multiplications by $\tau$, $\tau^c$. Take any $n\ge 2$ 
such that $\iota^c_n \phi \in C^\natural  (Z'\times\Bbb G_m ,Z\times \Bbb G_m )$; then $
{\iota}^c_n (\chi \otimes{\tau}^c )$, ${\iota}^c_n \sigma (\chi\otimes{\tau}^c )\sigma $ lie in this subcomplex as well. Since $\chi =
{\kappa}^c ({\iota}^c_n (\chi \otimes{\tau}^c ))$ and $\rho \kappa (\iota_n \chi ):=\kappa (\iota_n \chi )\otimes\tau = {\kappa}^c ( {\iota}^c_n \sigma (\chi \otimes{\tau}^c )\sigma )$, one has 
$\chi - \rho \kappa (\iota_n \chi )=  d{\kappa}^c ({\iota}^c_n \phi )$, q.e.d.  \hfill$\square$

\medskip

{\bf 3.4. Specialization and Tate untwisting.} 3.4.1. Let $X$ be a smooth variety, $U\subset X$ its dense  open subset. We call such a datum simply a {\it pair}, and denote it by $(X,U)$. If $(X',U')$ is another pair, then a {\it correspondence} $\theta : (X',U') \to (X,U)$ is a correspondence $\theta \in \CC or (X',X)$ such that the restriction $\theta^o$  of $\theta$ to $U'$ lies in $\CC or (U',U)\subset \CC or (U',X)$. We have an $R$-category $\CC or (\CP a)$ whose objects are pairs and morphisms are correspondences $\CC or ((X',U'),(X,U))$. One has two evident faithful functors $\CC or (\CP a )\to \CC or (\CS m)$, $(X,U)\mapsto X,U$, and a functor $\CC or (\CP a)\to \CP_{\tr}$, $(X,U) \mapsto \Z [U/X ]:=\CC one (\Z [U]\to \Z [X])[-1]$. For $\theta$ as above we denote the morphism $\Z[U'/X']\to \Z [U/X]$ by $\theta^{rel}$.

A pair $(X,U)$ is said to be {\it smooth} if $Z:=X\smallsetminus U$ is a smooth divisor. Let $\CC or (\CP a^{\text{sm}})$ be the subcategory of smooth pairs. If $\theta$ as above is a correspondence between smooth pairs, then the divisor $Z\times X' $ in $X\times X'$ meets $\theta$ properly.   Set  $$Sp (\theta) := \theta \cap (Z \times X') \tag 3.4.1$$ where the intersection is computed in the cycle-theoretic sense. This is a cycle supported at $Z\times Z' \subset X\times X'$ and finite over $Z'$. Therefore $Sp(\theta )\in \CC or (Z' ,Z )$. 

We leave it to the reader to check that the specialization is compatible with the composition, i.e., $Sp (\theta\theta')=Sp (\theta )Sp (\theta')$. Thus we have a functor $Sp : \CC or (\CP a^{\text{sm}})\to \CC or (\CS m)$, $(X,U) \mapsto Z$, $\theta \mapsto Sp (\theta )$.

{\it Remark.} If $\theta \in \CC or (X' ,U)\subset \CC or ((X',U'),(X,U))$, then $\theta^{rel}$ is homotopic to 0 and $Sp (\theta )=0$.

3.4.2. We call $\Bbb T :=(\Bbb A^1 ,\Bbb G_m ) $ the {\it Tate pair}. So for every smooth $Z$ we have a smooth pair $Z\times \Bbb T := (Z\times\Bbb A^1 /Z\times \Bbb G_m )$. One has a canonical homotopy equivalence
$\Z^\Delta [Z\times \Bbb G_m/Z\times\Bbb A^1  ]\iso \Z^\Delta [Z](1)[1]$. 

Take any $\theta \in \CC or (Z'\times \Bbb T ,Z\times \Bbb T)$. By the above homotopy equivalence, we can consider $\theta^{rel}$ as a homotopy morphism 
$ \Z^\Delta [Z'](1)\to \Z^\Delta [Z](1)$. We also have $Sp (\theta ):   \Z [Z']\to \Z [Z]$.
 
\proclaim{\quad Proposition }  $\theta^{rel}$ is $\Bbb A^1$-homotopic to $Sp (\theta )(1)$.
\endproclaim

{\it Proof.} We use the notation from 3.3. Notice that $\iota_n ,\tau \in E$ from loc.~cit.~come from correspondences in $\CC or (\Bbb T,\Bbb T )$ which we denote  by $\tilde{\iota}_n $, $\tilde{\tau}$. 
By the remark  above, replacing $\theta$ by $\tilde{\tau }\theta\tilde{\tau}$ does not change $\theta^{rel}$ and $Sp (\theta )$, so we can assume that $\theta =\tilde{\tau}\theta\tilde{\tau}$.  Then $\theta^{rel}$ is homotopic to $\theta^o \in C^\flat_\tau$. By 3.3.2 applied to $\chi=\theta^o$,
 it remains to show that  $Sp (\theta )$ is homotopic to $\kappa (\iota_n \theta^o )$ for sufficiently large $n$.

Set $\tilde{C}^\flat := C (Z' \times \Bbb A^1 , Z\times \Bbb A^1 )$. We define a subcomplex $\tilde{C}^\natural \subset \tilde{C}^\flat$ and the morphism $\tilde{\kappa} : \tilde{C}^\natural \to C(Z',Z)$ in the same way as we have defined $C^\natural$, $\kappa$ in 3.3.1 replacing $\Bbb G_m$ by $\Bbb A^1$ (there is a single difference: to define the nice correspondences in $\CC or (Y\times \Bbb A^1 ,\Bbb A^1 )$, we look at $p_\infty$ only).  As in loc.~cit., for any $\theta'\in \tilde{C}^\flat$ one has $\tilde{\iota}_n \theta' \in\tilde{C}^\natural$ for sufficiently large $n$. For any $\gamma \in C$ one has  $\tilde{\iota}_n (\gamma \otimes\tilde{\tau})\in \tilde{C}^\natural$, but, at variance with $\kappa$, one has $\tilde{\kappa}\tilde{\iota}_n (\gamma \otimes\tilde{\tau})=0$ (for $\tilde{\kappa}$ includes the intersection with the diagonal $t=t'$ at $t=0$).

Since $\tilde{\tau}$ is 
$\Bbb A^1$-homotopic to zero in $\CC or (\Bbb A^1 ,\Bbb A^1 )$ (for its degree equals 0), so is 
 $\theta = \tilde{\tau}\theta$, i.e., $ \theta =d\psi$ in the complex $\tilde{C}^\flat$. Thus $\tilde{\kappa} (\tilde{\iota}_n \theta )=d (\tilde{\kappa}(\tilde{\iota}_n \psi ))$ for large $n$. Now $\tilde{\kappa} (\tilde{\iota}_n \theta )= \kappa (\iota_n \theta^o )- Sp' (\tilde{\iota}_n \theta )$ where $Sp' (\tilde{\iota}_n \theta )\in \CC or (Z',Z)$ is the intersection of the cycle $\tilde{\iota}_n \theta$ with the diagonal $t=t'$ at $t=0$. The latter equals the intersection of $\theta$ with the divisor $t^n = t'$ at $t=0$, which is $Sp (\theta )$ for large $n$, and we are done.
 \hfill$\square$

 \bigskip
 
\centerline{\bf \S 4. The motivic localization }

\medskip
This section deals with the subject of lectures 12, 13, 24 of \cite{MVW} and 3.1, 3.2 of \cite{Vo2}. 
The principal result  is the theorem in 4.4 (cf.~\cite{Vo2}  3.2.6,   3.1.12,  \cite{MVW} 13.8, 13.9,  \cite{Deg2}  5.1). We deduce it from the theory of residues (the Cousin resolution) for $\Bbb A^1$-homotopy invariant presheaves with transfers. The latter is formulated  in 4.6, the proofs are in \S 5. 
The original argument of Voevodsky is rather different (the Cousin resolution emerges in the aftermath on the  last page of \cite{MVW}), though both proofs are based on the same  key input (see 5.1, 5.3).

  \medskip
  
  {\bf 4.1. The Zariski  localization: a complaint.}   The Zariski localization of  plain presheaves is a particular case of the format from 1.11 for categories of schemes such as $\CS m$ equipped with the Zariski topology. 
  The local objects for the Zariski topology are  local schemes (which are localizations of the schemes from our category). A complex of presheaves is Zariski locally acyclic (see 1.11) if and only if it is stalkwise acyclic. The  subcategory $I^{\text{Zar}}$ in $D\CP\CS h$ of Zariski locally acyclic complexes  is  generated by all Mayer-Vietoris complexes $R[U\cap V]\to R[U]\oplus R[V]\to R[X]$ (here $U,V\subset X$ are  opens  such that $U\cup V=X$), see \cite{BG} (or modify the proof of the proposition in 4.2.1 for the Zariski site).

  Passing to presheaves with transfers, we see that the category $I_{\tr}^{\text{Zar} \bot}\subset D\CP\CS h_{\tr}$ from 2.3 consists of all  $F$'s that are Zariski local as complexes of mere presheaves.  But the category $I_{\tr}^{\text{Zar}}$ itself contains not only all complexes which are  Zariski stalkwise acyclic, but many other objects as well: e.g.~the Mayer-Vietoris complex  $R_{\tr}[U\cap V]\to R_{\tr}[U]\oplus R_{\tr}[V]\to R_{\tr}[X]$
  is {\it not}  Zariski stalkwise acyclic if $X$ is irreducible and $U,V\neq X$. So 
    the forgetful functor $o: D\CP\CS h_{\tr} \to D\CP\CS h$ does {\it not} send $I_{\tr}^{\text{Zar}}$ to  $I^{\text{Zar} }$, hence it
     does {\it not} commute with the Zariski localization, which makes the Zariski localization of arbitrary presheaves with transfers unwieldy. 
   
     The cause of the nuisance is that the image of a local scheme by a correspondence is merely a {\it semi-local} scheme which need not be  the disjoint union of local schemes. Thus a correspondence may {\it not} be continuous for the Zariski topology:
   the Zariski coverings do {\it not} define a Grothendieck topology on $\CC or(\CS m )$. 
 
The problem disappears if one replaces the Zariski topology by the Nisnevich one (see the lemma in 4.3). We will see that for perfect $k$ this  does not change the motivic localization.

\medskip

{\bf 4.2. A Nisnevich localization reminder.}  For details, see \cite{Ni} or \cite{TT} App.~E. To fix the setting, we play with the category $\CS m$;  it can be replaced by the category of all $k$-varieties (or all  Noetherian schemes of finite Krull dimension). 

{\it Remark.} The constructions  below work with obvious modification (replacing \'etale maps by local isomorphisms) for the Zariski topology.

4.2.1. The Nisnevich topology lies in between the Zariski and \'etale topologies. A Nisnevich covering of $X$ is an \'etale  covering $ U\to X$ such that for every field $K$ the map $U(K)\to X(K)$ is surjective, or, equivalently, $U/X$ admits a constructible section (see 4.2.2).  The local objects for the Nisnevich topology are  Henselian local schemes. 

 $\CS m$ equipped with the Nisnevich topology satisfies the conditions from 1.11, so the 
general format from loc.~cit.~readily applies. 
We have  the abelian category of Nisnevich sheaves $\CS h^{\text{Nis}}$,
the exact sheafification functor $\CP\CS h \to\CS h^{\text{Nis}}$, $F\mapsto F_{\text{Nis}}$, the right-admissible subcategory $I^{\text{Nis}}\subset D\CP\CS h$ of Nisnevich locally acyclic  complexes,\footnote{A complex $F$ of presheaves is Nisnevich locally acyclic iff it is Nisnevich stalkwise acyclic, i.e., for any $X\in\CS$, $x\in X$ the complex $F(X_x^h )$ is acyclic; here $X^h_x$ is the Henselization of $X$ at $x$ considered as a proobject of $\CS m$.}
  the category $I^{\text{Nis} \bot}$ of Nisnevich local complexes, and  the equivalences of the triangulated categories   $$D\CP\CS h/I^{\text{Nis}}\iso D\CS h^{\text{Nis}} \buildrel\sim\over\leftarrow I^{\text{Nis} \bot}. \tag 4.2.1$$ 
 
 {\it Exercise.} Suppose $\ell$ is prime to $char\, k$. Show that  the presheaf $X\mapsto H^1_{\acute e t}(X,\mu_\ell )$ is a Nisnevich sheaf; if $k$ is separably closed, then the stalk   $H^1_{\acute e t}(X^h_x ,\mu_\ell )$  vanishes if and only if $x$ is a closed point of $X$. 
  \footnote{ This sheaf was denoted by $\CO^\times /\ell$ in  \cite{MVW} 12.9. The exercise shows that the formula from \cite{MVW} 12.10 should be corrected: the product must take into account all points $x\in X$, not only the closed ones as in loc.~cit. }
 
 Suppose we have an \'etale covering $\{ U,V\}$ of $X\in\CS m$ such that $V$ is  an open subset of $X$, and for $Z:= X\smallsetminus V$ one has $U_Z := U\times_X Z \iso Z$. We call such $\{U,V\}$ a {\it Nisnevich Mayer-Vietoris covering}, and the corresponding complex of presheaves $R[ U_V ] \to R[U]\oplus R[V] \to R[X]$ (the differentials are the difference and the sum of  evident maps)  a {\it Nisnevich Mayer-Vietoris complex.} It lies in $I^{\text{Nis}}$.

\proclaim{\quad Proposition}  $I^{\text{Nis}}$  is generated by the Nisnevich Mayer-Vietoris complexes.\footnote{See \cite{Vo4} for a more general story.} 
\endproclaim

{\it Proof.}  By (i) of the proposition in  1.4.2 (for $D= I^{\text{Nis}}$, $S$ the set of all Nisnevich Mayer-Vietoris complexes), the assertion amounts to the next statement: Suppose $F\in I^{\text{Nis}}$ is such that for
each Nisnevich Mayer-Vietoris covering as above  the total complex of $F(X)\to F(U)\oplus F(V)\to F(U_V )$ is acyclic; then $F$ is acyclic. 
To show this, we check by induction by $n$ the next claim:  For every $X$ and $h\in H^a F (X)$ there is an open $V\subset X$ with $X\smallsetminus V$ of codimension $\ge n$ such that $h|_V \in H^a F (V)$ vanishes.

By the induction assumption, there is an open $V' \subset X$ with $Z:= X\smallsetminus V'$ of codimension $\ge n-1$ such that $h|_{V'}=0$. Let $z_i$ be the generic points of $Z$ of codimension $n-1$. Since $F\in I^{\text{Nis}}$, there is a Nisnevich covering of $X$ on which $h$ vanishes. Thus  there is an open $V''\subset X$  containing $V'$ and $z_i$'s,
and an \'etale $U'/V''$ such that $U'_{V'' \cap Z} \iso V'' \cap Z$ and $h|_{U'} =0$. Then $\{ U',V'\}$ is a Nisnevich Mayer-Vietoris covering of $V''$, so the exact sequence $H^{a-1}F (U'_{V'})\to H^a F (V'')\to H^a F (U')\oplus H^a F(V')$ shows that $h|_{V''}$ comes from $g\in H^{a-1}F (U'_{V'})$. By the induction assumption, there is a closed $Y' \subset U'_{V'}$ such that the image of $g$ in $H^{a-1}F (U'_{V'}\smallsetminus Y' )$ vanishes. Let $Y$ be the closure of the image of $Y'$ in $X$. Then $V:=V'' \smallsetminus Z\cap Y$ is the promised open subset of $X$. Indeed, $X\smallsetminus V$ has codimension $\ge n$ since $z_i \notin Z\cap Y$. Set $U:= U'_{V'' \smallsetminus Y}$; then $\{ U,V' \}$ is a Nisnevich Mayer-Vietoris covering of $V$, so $h|_V$ comes from $g|_{U_{V'}} \in H^{a-1}F (U_{V'})$. Since $U_{V'}\subset U'_{V'}\smallsetminus Y' $, one has $g|_{U_{V'}} =0$, hence $h|_V =0$, q.e.d. \hfill$\square$

In particular, we see that $I^{\text{Nis}}$ is  generated by a set of compact objects, so the Nisnevich localization fits into the picture of the proposition from 1.4.2.

 {\it Exercises.} (i) A presheaf $F$ is a Nisnevich sheaf iff for every Nisnevich Mayer-Vietoris covering  the sequence $0\to F(X)\to F(U)\oplus F(V)\to F(U_V )$ is exact.
 
 (ii) For any Nisnevich sheaf $F$ and $h\in H^n_{\text{Nis}}(X,F)$ there is a closed $Y\subset X$ of codimension $\ge n$ such that the image of $h$ in $H^n_{\text{Nis}}(X\smallsetminus Y,F)$ vanishes.
 
4.2.2.  As was mentioned in 1.11, the Nisnevich localization functor $D\CP\CS h \to D\CS h^{\text{Nis}} \iso I^{\text{Nis}\bot}\subset D\CP\CS h$ can be lifted to the level of complexes by means of the Godement resolution. A characteristic feature of the Nisnevich (as well as the Zarisky) topology is that  the Godement resolution admits here an economic adele-style version $\CC^{\text{Nis}}$.\footnote{Here ``economic" means, in particular,  that for a sheaf $F$
on a variety $X$ of dimension $n$ the resolution $\CC^{\text{Nis}} (F)$ has length $\le n$.} In the rest of the section we construct $\CC^{\text{Nis}}$.
We will not use it any heavily, so the reader who is not inclined to follow the details of the construction may pass directly to 4.3.

 For $X\in\CS m$ a {\it partition} $\{ Z_\alpha \}$ of $X$ always means a (finite) partition by locally closed reduced subschemes. For an $X$-scheme $Y/X$ a {\it constructible section}  is a
set-theoretic section $s: X\to Y$ such that  for some partition $\{ Z_\alpha \}$ of $X$ each $s|_{Z_\alpha} : Z_\alpha \to Y$ is a morphism of schemes.\footnote{I.e., the closure of $s(Z_\alpha )$ in $Y_{Z_\alpha}$, considered as a closed reduced subscheme in $Y_{Z_\alpha}$, projects isomorphically to $Z_\alpha$.}

For $X\in\CS m$ consider the category $\Xi (X)$ of  pairs $\xi =(  U^\xi /X ,  s^\xi )$, where   $ U^\xi /X $ is an \'etale $X$-scheme and $s^\xi$ is a constructible section of   $ U^\xi /X $; a morphism $\xi \to \xi'$ is a morphism of $X$-schemes $U^\xi \to U^{\xi'}$ which sends $s^\xi$ to $s^{\xi'}$. 
Let $\Xi_o (X)$ be the subcategory of those $\xi_o$ that each connected component of $U^{\xi_o}$ intersects  $s^{\xi_o} (X )$. Our $\Xi (X)$ admits finite products, for  any $\xi $, $\xi_o$ as above there is at most one morphism $\xi_o \to \xi$, and for each $\xi$ one can find some $\xi_o \to \xi$. Thus $\Xi (X)^{\text{op}}$ is directed,  $\Xi_o (X)$ is its cofinal subset.  Set $X_0^{\text{Nis}} := \limleft_{\Xi (X)} U^\xi = \limleft_{\Xi_{o} (X)} U^\xi$; this is a proobject of $\CS m / X$. 

 It fits into a simplicial proobject $X_\cdot^{\text{Nis}}/X$   defined as follows. For $n\ge 0$  all collections $(\xi_0 ,\ldots ,\xi_n )$ where  $\xi_i = (U^{\xi_{i}}/U^{\xi_{i-1}}, s^{\xi_i})\in \Xi (U^{\xi_{i-1}})$ (we set $U^{\xi_{-1}}:=X$)  form  a directed category $\Xi^n (X)$;  $\xi$'s with every $\xi_i $ in $\Xi_o$ comprise a  directed set $\Xi^n_o (X)$ cofinal in $\Xi^n (X)$. Our   $\Xi^\cdot (X)$ is a  cosimplicial system of categories: for a monotone map  $\phi : [0,n]\to [0,m]$ the corresponding functor $ \Xi^n (X)\to\Xi^m (X)$ is $\xi \mapsto \xi^\phi$,  $\xi^\phi_j := (U^{\xi_{\phi^{-1}[j]}}/U^{\xi_{\phi^{-1}[j-1]}}, s^{\xi_{\phi^{-1}[j],\phi^{-1}[j-1]}})$; here for $j\in [0,m]$ we set $\phi^{-1}[j ]:= \max \{ i : \phi (i)\le j\} \in [0,n]$ and $s^{\xi_{a ,b}}$ is the composition $s^{\xi_a}\ldots s^{\xi_{b+1}}$ if $a> b$ and the identity map if $a=b$. 
Now  $X^{\text{Nis}}_n := \limleft_{\Xi^n (X)} U^{\xi_n}$; the simplicial structure maps on   $X_\cdot^{\text{Nis}}/X$ are $\phi^{\text{Nis}}:  
X^{\text{Nis}}_m \to X^{\text{Nis}}_n$ defined as the $\Xi^n (X)$-projective limit of the maps $X^{\text{Nis}}_m \to U^{\xi^\phi_m} = U^{\xi_n}$.

Any morphism of varieties $f : Y\to X$ lifts naturally to a morphism of simplicial proobjects $f^{\text{Nis}}: Y^{\text{Nis}}_\cdot \to X^{\text{Nis}}_\cdot$ (for one has an evident pull-back functor $f^* : \Xi^\cdot (X)\to \Xi^\cdot (Y)$). Thus any presheaf $F$ yields a cosimplicial presheaf $F^\cdot_{\text{Nis}}$, $F^\cdot_{\text{Nis}}(X):=
 F(X_\cdot^{\text{Nis}})$ (i.e.,   $F^n_{\text{Nis}}(X) = \limright_{\Xi^n (X)}\! F(U^{\xi_n} )$), augmented by $F$.   Let $\CC^{\text{Nis}}(F)$ be its normalized complex and $\epsilon : F\to\CC^{\text{Nis}}(F)$ be  the  augmentation.

\proclaim{\quad Proposition} (i) The functor $F\mapsto \CC^{\text{Nis}}(F)$ is exact.  (ii) Each presheaf $\CC^{\text{Nis}}(F)^n$ is Nisnevich local (see 1.11), hence is a Nisnevich sheaf. (iii) The augmentation $\epsilon$ makes  $\CC^{\text{Nis}}(F)$  a resolution of the Nisnevich sheafification  of $F$.  (iv) Every 
 section of  $\CC^{\text{Nis}}(F)^n$  is supported in codimension $n$, so $\CC^{\text{Nis}}(F)(X)$ sits in degrees $[0,\dim X]$.
\endproclaim

{\it Proof.} (i) Since  $\Xi^n (X)$ is directed, the functor $F\mapsto F^n_{\text{Nis}}$ is exact, and  $\CC^{\text{Nis}}(F)^n$ is a natural direct summand in $F^n_{\text{Nis}}$. 

(ii)  By construction, the functor $F\mapsto F_{\text{Nis}}^n$ is $n+1$-iteration of one $F\mapsto F_{\text{Nis}}^0$ (i.e., $F_{\text{Nis}}^n =(F_{\text{Nis}}^{n-1})_{\text{Nis}}^0$). Replacing $F_{\text{Nis}}^{n-1}$ by $F$, we see that  it suffices to show that $F_{\text{Nis}}^0$ is Nisnevich local. Indeed, for every Nisnevich hypercovering $V_\cdot $ of $X$ the complex $F_{\text{Nis}}^0 (V_\cdot /X)$ is homotopically trivial, which follows easily  from the existence of a constructible (left) contraction of $V_\cdot /X$;\footnote{ Recall that a {\it left contraction} of an augmented simplicial object $V_\cdot /X$  is a datum $c_\cdot$ of sections $c_n : V_n \to V_{n+1}$ of $\partial_{0}: V_{n+1}\to V_n$, $n\ge -1$ (here $V_{-1}:=X$),  compatible with the simplicial structure  maps, in the sense that  for every monotone $\phi : [0,n] \to [0,m]$ one has  $c_n \phi_V = \tilde{\phi}_V c_m : V_m \to V_{n+1}$; here $\tilde{\phi} : [0,n+1]\to [0,m+1]$ is defined as $\tilde{\phi}(0)=0$, $\tilde{\phi}(a)= \phi (a-1)+1$ for $a\in [1,n+1]$.}
the details are  left to the reader.

(iii) Take any $X\in\CS m$ and a point $x\in X$; let $Y:=X^h_x$ be the  Henselization of $X$ at $x$ considered as a proobject of $\CS m$. To check that $F(Y )\to C^{\text{Nis}}(F)(Y )$ is a quasi-isomorphism, one shows that  the augmented simplicial proobject $Y^{\text{Nis}}_\cdot /Y$ is contractible. The contraction $c$ is unique; its construction is left to the reader.

(iv) Let us define a cofinal subset $\Xi^n_c $ of $\Xi_o^n$ which clarifies the structure of $F^n_{\text{Nis}}$.
 
 Let $\{ Z_\alpha \}$ be a stratification of $X$ with equidimensional strata. Its set $A$ of indices $\alpha$ is naturally ordered ($\alpha' \ge \alpha$ means $\bar{Z}_{\alpha'} \supset Z_{\alpha}$), so we have the  simplicial set $S(A)_\cdot $ where $S(A)_n$ is the set of ordered $n+1$-tuples $(\alpha_0 \ldots \alpha_n )$, $\alpha_0 \le \ldots\le \alpha_n$. 

 For  $\xi \in \Xi (X)$ we denote by $\{ Z_\alpha^\xi \}$ the induced stratification of $U^\xi$, and by 
 $U^\xi_\alpha$  the union of connected components of $U^\xi$ that contains $s^\xi (Z_\alpha )$.
We say that $\xi $ is  {\it $\{ Z_\alpha \} $-controlled} if  $U^\xi $ equals the disjoint union of $ U^\xi_\alpha$, for each $\alpha$ the intersection  $Z^\xi_{\alpha'}\cap U^\xi_\alpha $ is empty unless $\alpha'\ge \alpha$, and $Z^\xi_{\alpha}\cap U^\xi_\alpha \iso Z_\alpha$. Notice that $s^\xi |_{Z_\alpha}$ is the inverse to the latter isomorphism. 
 
We say that $\xi 
\in\Xi^n (X)$ is {\it $\{ Z_\alpha \} $-controlled} if
 each $\xi_i \in \Xi (U^{\xi_{i-1}})$ is $\{ Z_\alpha^{\xi_{i-1}}\}$-controlled; here  $\{ Z_\alpha^{\xi_{i-1}}\}$ is the stratification of $U^{\xi_{i-1}}$ induced by $\{ Z_\alpha \}$.  Then every $U^{\xi_m}$, $m\in [0, n]$, is the disjoint union of components $U^{\xi_m}_{\alpha_0 \ldots \alpha_m}$ labeled by $(\alpha_0 \ldots  \alpha_m )$ $\in S(A)_m$, defined as the intersection of the preimages of $U^{\xi_i}_{\alpha_i}\subset U^{\xi_i}$ for $i\le m$.

Suppose  $\xi 
\in\Xi^n (X)$, $n>0$, is $\{ Z_\alpha \} $-controlled. For $m\in [0,n-1]$ let $\sigma_m : [0,n]\to [0,n-1]$ be the monotone surjection with $\sigma_m (m)=\sigma_m (m+1)$. Then $\xi^{\sigma_m}\in\Xi^{n-1}(X)$
need not be $\{ Z_\alpha \} $-controlled, but there is a natural morphism $\bar{\xi}^{\sigma_m}\to \xi^{\sigma_m}$ in    $\Xi^{n-1}(X)$ with $\{ Z_\alpha \} $-controlled $\bar{\xi}^{\sigma_m}$. Namely, 
$U^{\bar{\xi}^{\sigma_m}_i}$ equals $U^{\xi^{\sigma_m}_i}$ if $i\le m$; otherwise $U^{\bar{\xi}^{\sigma_m}_i}$ is an open subset of $ U^{\xi^{\sigma_m}_i} = U^{\xi_{i+1}}$ equal to the union of all components $U^{\xi_{i+1}}_{\alpha_0 \ldots \alpha_{i+1}}$ with $\alpha_m =\alpha_{m+1}$.

  \proclaim{\quad Lemma} $\xi$'s in  $\Xi^n (X)$ that are  controlled by some stratification form a cofinal subset $\Xi^n_c (X)$ in $\Xi^n_o (X)$. If $n>0$, then for every $m\in [0,n]$ all $\bar{\xi}^{\sigma_m}$, $\xi\in \Xi^n_c (X)$, form a cofinal subset in $\Xi^{n-1}_c (X)$.  \hfill$\square$
\endproclaim

By the lemma, $F_{\text{Nis}}^n (X) = \limright_{\Xi^n_c (X)}\! F(U^{\xi_n} )= \limright_{\Xi^n_c (X)}\! \oplus F(U^{\xi_n}_{\alpha_0 \ldots\alpha_n } )$, and its submodule  $\CC^{\text{Nis}}(F)^n (X)$ is the inductive limit of the sum of components $F(U^{\xi_n}_{\alpha_0 \ldots\alpha_n } )$ with $\alpha_0 <\ldots <\alpha_n$. Since the image of $F(U^{\xi_n}_{\alpha_0 \ldots\alpha_n } )$ in $F_{\text{Nis}}^n (X)$ vanishes on the complement of $\bar{Z}_{\alpha_0}$, we are done with the proposition.
\hfill$\square$

For a complex of presheaves $F$ we denote by $\CC^{\text{Nis}}(F)$ the corresponding total complex. We see  that  $\CC^{\text{Nis}}$ is a DG endofunctor that lifts the Nisnevich localization. 
 
{\bf 4.3. Compatibility with transfers.} The Nisnevich topology, as opposed to the Zariski one, is correspondence friendly:

 \proclaim{\quad Proposition}  The Nisnevich coverings define a Grothendieck topology on $R_{\tr}[\CS m]$. 
 \endproclaim

{\it Proof.}   For a smooth $X$ each  
Nisnevich cover $U/X$ yields a sieve on $R_{\tr}[X]$ formed by those correspondences $\CC or (\cdot ,X)$ that factor through $U$. The axioms of Grothendieck topology (see e.g.~\cite{D1} Arcata I  6) are immediate, except the one which is the next lemma for $X$, $Y$ smooth:

\proclaim{\quad Lemma}  Let $X,Y$ be schemes, and  $\pi : U \to  X$ be a Nisnevich cover; suppose that $Y$ is normal. Then for any
 $\gamma \in \CC or (Y,X)$ one can find  a Nisnevich cover $\pi' :  V \to Y$  and $\tilde{\gamma}\in \CC or ( V,U)$ such that $\pi \tilde{\gamma}= \gamma \pi'$. 
 \endproclaim
 
 {\it Proof of Lemma.}
 Those $\gamma$ for which the assertion is true form an $R$-submodule of $\CC or (Y,X)$. Indeed, if $\gamma =\Sigma a_i \gamma_i$, $a_i \in R$, and we have found $V_i /Y$ and $\tilde{\gamma}_i$ for each $\gamma_i$, then we take for $V/Y$  the fiber product of $V_i$'s over $Y$, and define $\tilde{\gamma} $ as $ \Sigma a_i \tilde{\gamma}'_i$, where $\tilde{\gamma}'_i$ is the composition of $\tilde{\gamma}_i$ and the projection $V\to V_i$.

Thus  it suffices to consider the case of  $\gamma$ given by a reduced irreducible $\Gamma \subset X\times Y$. Consider the Nisnevich cover $U_\Gamma := U\mathop\times\limits_X \Gamma$ of $\Gamma$. Since $\Gamma$ is finite over $Y$, one can find a Nisnevich cover $\pi' : V\to Y$ together with a morphism of $\Gamma$-schemes $\theta : V_\Gamma := V\mathop\times\limits_Y \Gamma \to U_\Gamma$. 
Take for $\tilde{\gamma}$ the composition $V\to V_\Gamma \buildrel\theta\over\lra U_\Gamma \to U$ where the first arrow is the finite correspondence transposed to the projection $V_\Gamma \to V$.
\hfill$\square$

A {\it Nisnevich sheaf with transfers} is a presheaf with transfers which is a Nisnevich sheaf. Equivalently, this is an $R$-sheaf $F$ on the topology from the lemma such that for every  $X,Y\in R_{\text{tr}}[\CS m]$ the pull-back map $\CC or (Y,X)\times F(X)\to F(Y)$ is $R$-bilinear. These objects form an abelian category $\CS h^{\text{Nis}}_{\tr}$, and we have the exact sheafification functor $\CP\CS h_{\tr} \to \CS h^{\text{Nis}}_{\tr}$, $F\mapsto F_{\text{Nis}}$, left adjoint to the fully faithful embedding $\CS h^{\text{Nis}}_{\tr}\hra \CP\CS h_{\tr}$. The forgetful functors $o: \CP\CS h_{\tr}\to\CP\CS h$, $ \CS h^{\text{Nis}}_{\tr} \to \CS h^{\text{Nis}}$ commute with the sheafification and the embedding.

Let $I^{\text{Nis}}_{\tr}$ be the triangulated subcategory of $D\CP\CS h_{\tr}$  generated by all Nisnevich Mayer-Vietoris complexes $\Z [U_V ]\to \Z [U]\oplus\Z [V]\to \Z [X]$ (see 4.2.1). It is right admissible (by (i) of the proposition from 1.4.2). The forgetful functor $o: D\CP\CS h_{\tr}\to D\CP\CS h$ sends $I^{\text{Nis}}_{\tr}$ to $I^{\text{Nis}}$ (by the lemma) and $I^{\text{Nis}\bot}_{\tr}$ to $I^{\text{Nis}\bot}$ (by the proposition in 4.2.1), i.e., $o$ commutes with the Nisnevich localization $\CC^{\text{Nis}}$. Since $o$ is conservative, $I^{\text{Nis}}_{\tr}$, $ I^{\text{Nis}\bot}_{\tr} $ are formed exactly by those complexes  of presheaves  with transfers that are
 Nisnevich locally acyclic, resp.~Nisnevich local,  as  mere complexes of presheaves. So the Nisnevich sheafification  yields  equivalences  $$D\CP\CS h_{\tr}/I^{\text{Nis}}_{\tr}\iso D\CS h^{\text{Nis}}_{\tr} \buildrel\sim\over\leftarrow I^{\text{Nis}\bot}_{\tr}. \tag 4.3.1$$ 
  The t-structure picture here is that of (i) from the lemma in 1.2.
 
The canonical resolution $\CC^{\text{Nis}}$ from 4.2.2 is compatible with transfers. Namely,
any $X\in \CS m$ yields a simplicial proobject $\Z [X^{\text{Nis}}_\cdot ]$ in $\CC or(\CS m )$ (see 4.2.2), and, as in the lemma above, any $\gamma \in \CC or (Y,X)$ lifts naturally to a morphism $\Z [Y^{\text{Nis}}_\cdot ]\to \Z [X^{\text{Nis}}_\cdot ]$.
Therefore for any presheaf with transfers  $F$  its  resolution $\CC^{\text{Nis}}(F)$   has a canonical transfer structure, and the morphism $F\to \CC^{\text{Nis}}(F)$ is compatible with transfers. So the Nisnevich sheafification $F_{\text{Nis}} =H^0 \CC^{\text{Nis}}(F)$  has a natural transfer structure as well.  
We see that $\CC^{\text{Nis}}$ is a DG endofunctor of $C\CP\CS h_{\tr}$ that lifts the Nisnevich localization.

{\it Remark.} The usual Godement resolution \cite{Go} for the Nisnevich topology is compatible with transfers as well.

Recall that $D\CP\CS h_{\tr}$ has a natural DG structure  $
\mathop{\vtop{\ialign{#\crcr
  \hfil\rm $\CP_{\tr}$\hfil\crcr
  \noalign{\nointerlineskip}\rightarrowfill\crcr
  \noalign{\nointerlineskip}\crcr}}} $ (see 2.1.4). Therefore, as in 1.5.5,  the DG subcategories
$\mathop{\vtop{\ialign{#\crcr
  \hfil\rm $\CI^{\text{Nis}}_{\tr}$\hfil\crcr
  \noalign{\nointerlineskip}\rightarrowfill\crcr
  \noalign{\nointerlineskip}\crcr}}},    \mathop{\vtop{\ialign{#\crcr
  \hfil\rm $\CI^{\text{Nis}}_{\tr}$\hfil\crcr
  \noalign{\nointerlineskip}\rightarrowfill\crcr
  \noalign{\nointerlineskip}\crcr}}}{\!}^\bot    \subset \mathop{\vtop{\ialign{#\crcr
  \hfil\rm $\CP_{\tr}$\hfil\crcr
  \noalign{\nointerlineskip}\rightarrowfill\crcr
  \noalign{\nointerlineskip}\crcr}}} $  of objects that belong to $I^{\text{Nis}}_{\tr}, I^{\text{Nis}\bot }_{\tr} \subset D\CP\CS h_{\tr}= \mathop{\vtop{\ialign{#\crcr
  \hfil\rm $\CP^{\text{tri}}_{\tr}$\hfil\crcr
  \noalign{\nointerlineskip}\rightarrowfill\crcr
  \noalign{\nointerlineskip}\crcr}}}$   are DG structures  on $I^{\text{Nis}}_{\tr}, I^{\text{Nis}\bot }_{\tr}$, and  one has a quasi-equivalence of DG categories $
\mathop{\vtop{\ialign{#\crcr
  \hfil\rm $\CI^{\text{Nis}}_{\tr}$\hfil\crcr
  \noalign{\nointerlineskip}\rightarrowfill\crcr
  \noalign{\nointerlineskip}\crcr}}}{\!}^{\bot}
  \iso \mathop{\vtop{\ialign{#\crcr
  \hfil\rm $\CP_{\tr}$\hfil\crcr
  \noalign{\nointerlineskip}\rightarrowfill\crcr
  \noalign{\nointerlineskip}\crcr}}} /\mathop{\vtop{\ialign{#\crcr
  \hfil\rm $\CI^{\text{Nis}}_{\tr}$\hfil\crcr
  \noalign{\nointerlineskip}\rightarrowfill\crcr
  \noalign{\nointerlineskip}\crcr}}}$.

{\bf 4.4. The motivic localization.} Set 
$\CI_{\tr}^{\text{Nis}} := (\mathop{\vtop{\ialign{#\crcr
  \hfil\rm $\CI^{\text{Nis}}_{\tr}$\hfil\crcr
  \noalign{\nointerlineskip}\rightarrowfill\crcr
  \noalign{\nointerlineskip}\crcr}}} )^{\text{perf}}    \subset \CP_{\tr}$, i.e., $\CI_{\tr}^{\text{Nis}}$ is  the homotopy idempotent completion of the full pretriangulated DG subcategory strongly  generated by  the Nisnevich Mayer-Vietoris complexes. As in 2.3, let 
 $  \CI_{\tr}^{\Delta\text{Nis}} \subset \CP_{\tr}$ be the corresponding subcategory strongly generated    by $  \CI_{\tr}^{\Delta} $ and $\CI_{\tr}^{\text{Nis}}$. 
 Then $I^{\Delta  \text{Nis}}_{\tr} :=( \mathop{\vtop{\ialign{#\crcr
  \hfil\rm $\CI^{\Delta\text{Nis}}_{\tr}$\hfil\crcr
  \noalign{\nointerlineskip}\rightarrowfill\crcr
  \noalign{\nointerlineskip}\crcr}}} )^{\text{tri}}$ is the full triangulated subcategory  of $ D\CP\CS h_{\tr}$  generated by $I^{\text{Nis}}_{\tr}$ and $ I^\Delta_{\tr}$.

\proclaim{\quad Theorem} Suppose that the base field $k$ is perfect.

 (i)  The subcategories $  I^\Delta_{\tr}, I^{\text{Nis}}_{\tr}\subset D\CP\CS h_{\tr}$ are compatible (see 1.3).

(ii) One has 
$I^{\Delta \text{Zar}}_{\tr}= I^{\Delta \text{Nis}}_{\tr}$ (see 2.3), i.e.,  $$
\mathop{\vtop{\ialign{#\crcr
  \hfil\rm $\CD^{\text{eff}}_{\CM}$\hfil\crcr
  \noalign{\nointerlineskip}\rightarrowfill\crcr
  \noalign{\nointerlineskip}\crcr}}}{\!}^{\text{tri}} \iso D\CP\CS h_{\tr}  /I^{\Delta \text{Nis}}_{\tr}= I^{\Delta\bot}_{\tr}\cap I^{\text{Nis}\bot}_{\tr} , \tag 4.4.1$$   $$\CD_\CM^{\text{eff}}{}^{\text{tri}}\iso (D\CP\CS h_{\tr}  /I^{\Delta \text{Nis}}_{\tr})^{\text{perf}}. \tag 4.4.2$$ 
Therefore one has a quasi-equivalence of the DG categories $\CD^{\text{eff}}_\CM \iso \CP_{\tr}/^\kappa \CI_{\tr}^{\Delta\text{Nis}}$, and
the motivic localization  endofunctor $\CC^\CM  $ of $D\CP\CS h_{\tr}$ can be written as $$ \CC^\CM = \CC^\Delta \CC^{\text{Nis}} =\CC^{\text{Nis}} \CC^\Delta . \tag 4.4.3$$
\endproclaim

{\it Remarks.} (a) For an interpretation of (i) on finite level (inside  $\CP_{\tr}$), see 1.4.3.

(b) The theorem implies that in the definition of $\CD_\CM^{\text{eff}}$ (see 2.3) one can replace $\CC or(\CS m )$ by its full subcategory of affine varieties:
 the category of motives does not change. We do not know if this remains true for non perfect $k$.

(c) We see that the  identification $\mathop{\vtop{\ialign{#\crcr
  \hfil\rm $\CD^{\text{eff}}_{\CM}$\hfil\crcr
  \noalign{\nointerlineskip}\rightarrowfill\crcr
  \noalign{\nointerlineskip}\crcr}}} {\!}^{\text{tri}} = I_{\tr}^{\Delta\bot}/I^{\Delta\bot}_{\tr}\cap I^{\text{Nis}}_{\tr}$ yields a  t-structure on $\mathop{\vtop{\ialign{#\crcr
  \hfil\rm $\CD^{\text{eff}}_{\CM}$\hfil\crcr
  \noalign{\nointerlineskip}\rightarrowfill\crcr
  \noalign{\nointerlineskip}\crcr}}} {\!}^{\text{tri}}$ whose core is the abelian category of $\Bbb A^1$-homotopy invariant Nisnevich sheaves with transfers $\CS h^{\Delta \text{Nis}}_{\tr}$.
 To sum the picture up, we have a tensor DG category $\mathop{\vtop{\ialign{#\crcr
  \hfil\rm $\CP_{\tr}$\hfil\crcr
  \noalign{\nointerlineskip}\rightarrowfill\crcr
  \noalign{\nointerlineskip}\crcr}}} $ and its three quotients. 
The corresponding homotopy categories carry non-degenerate t-structures with cores $\CP\CS h_{\tr}$, $\CS h^{\text{Nis}}_{\tr}$, $\CP\CS h_{\tr}^\Delta$,  $\CS h^{\Delta \text{Nis}}_{\tr}$. The projections  $\mathop{\vtop{\ialign{#\crcr
  \hfil\rm $\CP_{\tr}$\hfil\crcr
  \noalign{\nointerlineskip}\rightarrowfill\crcr
  \noalign{\nointerlineskip}\crcr}}} {\!}^{\text{tri}} \twoheadrightarrow  \mathop{\vtop{\ialign{#\crcr
  \hfil\rm $\CP_{\tr}$\hfil\crcr
  \noalign{\nointerlineskip}\rightarrowfill\crcr
  \noalign{\nointerlineskip}\crcr}}}  {\!}^{\text{tri}}/I^{\text{Nis}}_{\tr}$ and $\mathop{\vtop{\ialign{#\crcr
  \hfil\rm $\CD^{\text{eff}}_{p\CM}$\hfil\crcr
  \noalign{\nointerlineskip}\rightarrowfill\crcr
  \noalign{\nointerlineskip}\crcr}}}  {\!}^{\text{tri}} \twoheadrightarrow   \mathop{\vtop{\ialign{#\crcr
  \hfil\rm $\CD^{\text{eff}}_{\CM}$\hfil\crcr
  \noalign{\nointerlineskip}\rightarrowfill\crcr
  \noalign{\nointerlineskip}\crcr}}} {\!}^{\text{tri}}$ are t-exact, while  $\mathop{\vtop{\ialign{#\crcr
  \hfil\rm $\CP_{\tr}$\hfil\crcr
  \noalign{\nointerlineskip}\rightarrowfill\crcr
  \noalign{\nointerlineskip}\crcr}}} {\!}^{\text{tri}} \twoheadrightarrow   \mathop{\vtop{\ialign{#\crcr
  \hfil\rm $\CD^{\text{eff}}_{p\CM}$\hfil\crcr
  \noalign{\nointerlineskip}\rightarrowfill\crcr
  \noalign{\nointerlineskip}\crcr}}} {\!}^{\text{tri}}$  is right t-exact (its right adjoint is exact). The tensor products are right t-exact.

(d) We do not know if $\mathop{\vtop{\ialign{#\crcr
  \hfil\rm $\CD^{\text{eff}}_{\CM}$\hfil\crcr
  \noalign{\nointerlineskip}\rightarrowfill\crcr
  \noalign{\nointerlineskip}\crcr}}} {\!}^{\text{tri}} $ is  equivalent to the derived category of  $\CS h_{\tr}^{\Delta\text{Nis}}$ (cf.~the exercise in 3.1.2).

{\it Plan of the proof.}  First, we reformulate the theorem in concrete terms, reducing it to the theorem in 4.5. The latter is deduced in 4.6.4  from the 
 (itself   important) theory of residues for $\Bbb A^1$-homotopy invariant presheaves with transfers, which is the subject of 4.6. The proofs of the principal facts about  residues (the proposition in 4.6.1 and  the theorem in 4.6.3) and of the lemma in 4.6.4  are given in \S 5.

\proclaim{\quad  4.5.  Theorem} If  $k$ is perfect and  $F\in \CP\CS h_{\tr}$ is $\Bbb A^1$-homotopy invariant, then for every $X\in\CS m$  one has  $$H^\cdot_{\text{Nis}}(X,F)\iso H^\cdot_{\text{Nis}}(X\times\Bbb A^1 ,F), \quad H^\cdot_{\text{Zar}}(X,F)\iso H^\cdot_{\text{Nis}}(X,F). \tag 4.5.1$$ In particular, the Zariski and Nisnevich sheafifications of $F$ coincide: $F_{\text{Zar}} \iso F_{\text{Nis}}$.
\endproclaim

Here  $H^a_{\text{Nis}}(X,F):= H^a (X_{\text{Nis}}, F_{\text{Nis}}) =H^a  \CC^{\text{Nis}}(F)(X) $, etc. 

{\it Remark. } Since the small Zariski and Nisnevich topologies have finite cohomological dimension, the  theorem automatically implies that (4.5.1) remains true for any complex $F\in I^{\Delta\bot}_{\tr}$.

Let us show that the above theorem implies the theorem in 4.4:

(a) {\it 4.4(i) is equivalent to the first isomorphism in (4.5.1).}
By  1.3(ii),  4.4(i) amounts to the fact that   the Nisnevich localization $\CC^{\text{Nis}}$ preserves $I^\Delta_{\tr}$ and $I^{\Delta\bot}_{\tr}$. 

The statement about $I^\Delta_{\tr}$ is easy (we need not assume 4.5, and make no use of transfers). Indeed, for any presheaf $F$ the canonical embedding $\CC^{\text{Nis}}(F)\hra \CC^\Delta (\CC^{\text{Nis}}(F))$ is  the composition of natural embeddings $\CC^{\text{Nis}}(F)\hra \CC^{\text{Nis}}(\CC^\Delta (F)) \hra \CC^\Delta (\CC^{\text{Nis}}(F))$. Now if $F\in I_{\tr}^\Delta$, then $\CC^\Delta (F)$ is acyclic, hence the morphism $\CC^{\text{Nis}}(F)\hra \CC^\Delta (\CC^{\text{Nis}}(F))$ vanishes in $D\CP\CS h$, which implies that $\CC^{\text{Nis}}(F)\in I^\Delta_{\tr}$.

The statement about $I^{\Delta\bot}_{\tr}$ says that for every 
$\Bbb A^1$-homotopy invariant complex of presheaves with transfers $F$ its Nisnevich localization $\CC^{\text{Nis}}(F)$ is again $\Bbb A^1$-homotopy invariant, i.e., for every $X\in \CS m$  one has $H^\cdot_{\text{Nis}}(X,F)\iso H^\cdot_{\text{Nis}}(X\times\Bbb A^1 ,F)$, which by the above Remark is the first isomorphism in (4.5.1), q.e.d.

(b) {\it 4.4(ii) follows from the second isomorphism in (4.5.1).} Indeed,  4.4(ii) amounts to the fact that 
$I^{\text{Nis}\bot}_{\tr}\cap I^{\Delta\bot}_{\tr} = I^{\text{Zar}\bot}_{\tr}\cap I^{\Delta\bot}_{\tr}$. 
Since $I_{\tr}^{\text{Nis}\bot}\subset I^{\text{Zar}\bot}_{\tr}$, this means that 
every $\Bbb A^1$-homotopy invariant complex of presheaves with transfers
 which is Zariski local (as a mere complex of presheaves) is automatically Nisnevich local. It suffices to check that for every 
$F \in  I^{\Delta\bot}_{\tr}$ its Nisnevich localization coincides, as an object of $D\CP\CS h$, with the Zariski localization $F^{\text{Zar}}$ of $F$ (in the sense of mere presheaves),\footnote{Hence, in particular,  $F^{\text{Zar}}$ acquires a transfer structure (cf.~4.1).} which by the above Remark  is the  second isomorphism in (4.5.1), q.e.d.

{\it Exercise} (Gaitsgory). The second isomorphism in (4.5.1) is equivalent to  4.4(ii) combined with the fact that  $o: D\CP\CS h_{\tr} \to D\CP\CS h$ sends $I_{\tr}^{\Delta\bot} \cap I_{\tr}^{\text{Zar}}$ to $I^{\text{Zar}}$.

The next corollary of the  theorems in 4.4, 4.5 captures their essence:

\proclaim{\quad Corollary} If  $Y$ is a smooth local scheme,\footnote{I.e., $Y$ is the localization of a smooth variety at a point. In fact, by    (ii) of the theorem in 4.6.3, the statement is valid whenever $Y$ is a smooth {\it semi-local} scheme.} e.g.~$Y=\Spec \, k$, then the functor $X$ $\mapsto \CC^\Delta (\Z [X] )(Y)$  has Nisnevich local nature, i.e., for any Nisnevich hypercovering $U_\cdot $ of $X$ the map $\CC^\Delta (\Z [U_\cdot ] )(Y)\to \CC^\Delta (\Z [X] )(Y)$ is a quasi-isomorphism. \hfill$\square$
\endproclaim

 {\bf 4.6. Residues and the Cousin resolution.} 4.6.1.  We use the notation of 3.4.   Let $X$ be a smooth variety, $Z\subset X$ a smooth divisor, i.e., we have a smooth pair $(X,X\smallsetminus Z)$. 

{\it Definition.}  $\theta\in \CC or (Z\times\Bbb T , (X,X\smallsetminus Z))$ is a {\it link correspondence} (for $(X,X\smallsetminus Z)$, or simply at $Z$) if  $Sp (\theta ) \in \CC or (Z,Z)$  equals $ \id_Z$. 

Using the homotopy equivalence $\Z^\Delta [ Z\times\Bbb G_m /Z\times \Bbb A^1 ] \iso \Z^\Delta [Z](1)[1] $, we can rewrite  $\theta^{rel}$ as a  morphism in the homotopy category of premotives $$\theta^{rel}: \,\Z^\Delta [Z](1)[1] \to\Z^\Delta [ (X\smallsetminus Z)/X]. \tag 4.6.1$$

\proclaim{\quad Proposition} Let $X$, $Z$ be as above, and $P\subset Z$ be a finite subset of  points. Suppose that the base field $k$ is  perfect and $X$ is quasi-projective.

(i) There is a  Zariski neighborhood $U\subset X$ of $P$ such that $(U,U\smallsetminus Z) $ admits  a  link correspondence $\theta$.\footnote{For this statement $k$ need not be perfect: it suffices to demand that $P$ consists of smooth points (i.e., the generic points of smooth locally closed subschemes).}

(ii) For every such $U$ and link correspondences  $\theta$, $\theta'$  for $(U,U\smallsetminus Z)$ one can find an  open  $Z' \subset U\cap Z$ which contains $ P$ such that the restrictions $\Z^\Delta [Z'](1)[1]\to \Z^\Delta [(U\smallsetminus Z)/U]$ of 
 $\theta^{rel}$, $\theta^{\prime rel} $ to $ Z'$ are $\Bbb A^1$-homotopic.  
\endproclaim
 
For a proof, see 5.4. 

{\it Remarks.} (a) It follows from 5.3 and  5.2 that  $U$ and $\theta$ from (i) can be chosen so that $\theta^{rel}$ is a homotopy equivalence.

(b) In 5.4 we prove the proposition under the additional assumption that $k$ is infinite. The case of finite $k$ follows then by the next trick. A standard argument (localizing $R$ at two different primes, etc.) shows that one can assume that certain prime $\ell$ is invertible in $R$.  
 Let $k'$ be a maximal pro-$\ell$-extension of $k$,  $G:=\text{Gal} (k'/k)\simeq \Bbb Z_\ell$. For any
 $Y_1 , Y_2  \in \CS m_k$ we have the $R$-module $\CC or (Y_{1 k'}, Y_{2 k'})_{k'} $  of correspondences between the  $k'$-varieties $Y_{1 k'}, Y_{2 k'}$;
  the group $G$ acts naturally on it, and $\CC or (Y_1 ,Y_2 )= \CC or (Y_{1 k'}, Y_{2 k'})_{k'}^G$. Due to the condition on  $R$, the functor of  $G$-invariants is exact on $R[G]$-modules, 
 so the existence of a link correspondence over $k'$ implies that over $k$. Same is true if we consider $\Bbb A^1$-localized correspondences, so  (ii) of the proposition for $k$ follows from that for $k'$, q.e.d. 
 
 {\it Question.} Is it true that  morphisms $\theta^{rel}$ for all
 link correspondences at $Z$ localized at $P$ form naturally a {\it contractible} space? That is, can one present, after the localization at $P$, a natural complex $C$ of $R$-modules  with $H^0 C=R$, $H^{\neq 0} C=0$, and a morphism $C\to \Hom (\Z^\Delta [Z](1)[1],\Z^\Delta [ (X\smallsetminus Z)/X])$ whose image contains the $\theta^{rel}$'s?

\medskip

4.6.2. Assuming the above proposition, let us define the residue maps.

 {\it Suppose that $k$ is perfect}. This means that any reduced scheme $X$ of finite type has a smooth open dense subscheme. As in 2.1.1, for  $F\in \CP\CS h$ we have a well-defined abelian group $F(\eta_X )$; here $\eta_X$ is the union of the generic points  of $X$.

Let $F$ be  an $\Bbb A^1$-homotopy invariant  presheaf with transfers.  As in 2.2,  we have  $F_{-n}\in \CP\CS h_{\tr}$, $n\ge 0$; they are all   $\Bbb A^1$-homotopy invariant.\footnote{In fact, all functors $\CH om (P, \cdot )$, $P\in \CP_{\tr}$, from  2.2 preserve $I_{\tr}^{\Delta\bot} \subset D\CP\CS h_{\tr}$.} 

 Let $T$ be a semi-local reduced  $k$-scheme of Krull dimension 1 which is  localization of a $k$-scheme of finite type. Let $s=s_T$ be the (disjoint) union of closed points, $\eta =\eta_T := T\smallsetminus s$ the union of the generic points.  Let us define the residue map  $$\Res = \Res^T =\Res_{s }: F(\eta )\to F_{-1}(s ). \tag 4.6.2$$  
 
 If $T$ is regular, then  $(T, s)$ is the localization of a pair $(X,Z)$,  where $X$ is a smooth affine variety, $Z\subset X$ a smooth divisor, at the generic points of $Z$. Recall  that $F_{-1}(Z)= \Hom (\Z^\Delta [Z](1)[1],F)$ (see 2.2). 
 Replace $X$ by a base of Zariski neighborhoods of $s$;  then the  pull-back maps for   morphisms (4.6.1) provide, by the proposition, a well-defined map $\Coker (F(T)\to F(\eta ))\to F_{-1}(s)$, which yields $\Res_s$. 

{\it Remark.} Let $X$, $Z$ be as above. By the proposition,  we have a natural morphism  $\Res_Z : F(X\smallsetminus Z)\to (F_{-1})_{\text{Zar} } (Z)$ defined as the pull-back map for morphisms (4.6.1) on all possible opens of $ Z$. For  $\phi \in F(X\smallsetminus Z)$ one has $\Res_s (\phi_\eta )=(\Res_Z (\phi ))_s$.

 For arbitrary $T$, let $\tilde{T}$ be the normalization of $T$; then $\Res$ is the composition $F(\eta )\to F_{-1}(s_{\tilde{T}})\to F_{-1}(s )$, where the first arrow is the residue map for $\tilde{T}$ and the second one is the transfer map for the finite projection $s_{\tilde{T}}\to s$.

\medskip

4.6.3. As above, our $k$ is perfect. For $X\in\CS m$  denote by $X^{(n)}$ the set of its points  of codimension $n$. For $x\in X^{(n)}$ let $\eta_x \subset X$ be the corresponding generic point.  For a presheaf  $F$  set $Cous(F)^n (X):= \mathop\oplus\limits_{x\in X^{(n)}} F_{-n}(\eta_x )$.

A smooth morphism $f : Y\to X$ yields a pull-back map $f^* : Cous(F)^n (X)\to Cous(F)^n (Y)$ which sends $F_{-n}(\eta_x )$ to $\oplus \, F_{-n}(\eta_y )$ where the $y$ are the generic points of $f^{-1}(x)$. If $F$ has a transfer structure and   $g: Z\to X$ is any morphism in $\CS m$ of codimension $m$, then for every $a$ such that $a, a+m \ge 0$ we have  the  trace map $g_* : 
Cous(F_{-a-m})^{ n-m}(Z)\to Cous(F_{-a})^n (X)$ whose non-zero components are  the transfers $F_{-a-n}(\eta_z )\to F_{-a-n}(\eta_{g(z)} )$ for $\eta_z $ finite over $ \eta_{g(z)}$.
The pull-back and  the trace  are  compatible with morphisms of $F$'s and with compositions of $f$'s and $g$'s;  if $f$, $g$ are as above and $Y \buildrel{g'}\over\leftarrow Y\mathop\times\limits_X Z \buildrel{f'}\over\rightarrow Z$ are the projections, then $f^* g_* = g'_* {f'}^*$.

In particular, the pull-back maps for smooth morphisms make each $Cous (F)^n$ a presheaf on the small smooth topology of any $X\in\CS m$. 

\proclaim{\quad Lemma}  Each $Cous (F)^n$ is a  local presheaf  for the small Nisnevich and Zariski topologies,\footnote{See 1.11; in fact,  $Cous (F)^n$  are flabby sheaves (see \cite{SGA4} V 4.1 for the definition). } hence $H^{>0}_{\text{Nis}}(X,  Cous (F)^n )= H^{>0}_{\text{Zar}}(X,  Cous (F)^n )= 0$.  \hfill$\square$
\endproclaim

There is an evident morphism $ F(X)\to Cous (F)^0 (X)$ compatible with the pull-back maps for smooth morphisms of $X$'s. By the lemma, it extends by adjunction to the morphisms of the sheafifications $F_{\text{Zar}}\to F_{\text{Nis}}\to Cous (F)^0$.

Now suppose that  $F$ is an $\Bbb A^1$-homotopy invariant presheaf with transfers. We define the {\it   Cousin differential } $$d=d_X : \, \CC ous (F)^n (X) \to \CC ous(F)^{n+1}(X) \tag 4.6.3$$  as follows. Our $d$ is an $X^{(n)}\times X^{(n+1)}$-matrix whose  entries are  maps $d_{xy} : F_{-n} (\eta_x )\to F_{-n-1}(\eta_y )$.  If $y$  lies in the closure of $x$, then $d_{xy}:= \Res_{s}$ where $s=\eta_Y$ and $T$ is the localization of the closure of $x$ at $y$; otherwise $d_{xy}=0$. Since every  $\nu\in F_{-n}(\eta_x )$ comes from an open subscheme of the closure of $x$, there is only finitely many $y$'s such that  $d_{xy}(\nu )\neq 0$, so $d$ is correctly defined.

\proclaim{\quad  Theorem}  (i) One has $d^2 =0$, so $d$  makes $Cous (F)$ a complex. For smooth $f$ and proper $g$ as above the corresponding maps $f^* : Cous (F)(X)\to Cous (F)(Y)$ and $g_* : Cous (F_{-a-m})(Z)[-m]\to Cous (F_{-a})(X)$  are morphisms of complexes.

(ii) The morphism $F_{\text{Zar}}\to Cous (F)^0$  makes  $Cous (F)$ a resolution of $F_{\text{Zar}}$ as sheaves on the small Zariski topology. In fact, $\CC ous (F)(\hat{X})$ is a resolution of $F(\hat{X})$ if $\hat{X}$ is any smooth semi-local scheme.\footnote{Which means that $\hat{X}$ is localization of an affine smooth scheme at a finite subset of points.}
\endproclaim

 For a proof, see 5.6.

{\it Remark.}
The last statement in (ii)  explains  why $F_{\text{Zar}}$ carries a transfer structure (avoiding the reference to the Nisnevich localization). Indeed, a correspondence sends a local scheme to a {\it semi-local} one.  Since  sections of $F_{\text{Zar}}$ on the latter are the same as those of $F$,  we know how the correspondence acts on them.

4.6.4. We need the following auxiliary lemma (to be proved in 5.1.2):

\proclaim{\quad Lemma} 
Let $F$ be an $\Bbb A^1$-homotopy invariant presheaf with transfers, $\eta$ be  the generic point of a smooth variety. Then $F(\Bbb A^1_\eta ) \iso H^0_{\text{Zar}} (\Bbb A^1_\eta ,F)$, $H^1_{\text{Zar}} (\Bbb A^1_\eta ,F)=0.$
\endproclaim

Let us deduce  theorem  4.5 from the above theorem and  the lemma.

 Of course, 4.6.3(ii) implies that $\CC ous (F)$ is a resolution of $F_{\text{Nis}}$ as sheaves on small Nisnevich topology as well, hence $F_{\text{Zar}}=F_{\text{Nis}}$. By 4.6.3(ii) and the lemma in 4.6.3, one has $R\Gamma (X_{\text{Zar}}, F_{\text{Zar}})= \CC ous (F)(X)= R\Gamma (X_{\text{Nis}}, F_{\text{Nis}})$, so we have  the second isomorphism in (4.5.1). 
 
 It remains  to show that the pull-back map $p_X^* : \CC ous (F)(X)\to \CC ous (F) (X\times\Bbb A^1 )$ for the projection $p_X : X\times \Bbb A^1 \to X$ is a quasi-isomorphism.  Our complexes carry natural  finite filtrations: on $Cous (F) (X)$ this is the stupid filtration, and on $Cous (F) (X\times \Bbb A^1 )$  the filtration whose $m$th term is formed by all  $F_{-n}(\eta_z )$'s such that  $ p_X (z )$ has codimension $\ge m$ in $X$. Then $p^*_X$ is compatible with the filtrations. We will check that
 $p^*_X$ is, in fact, a {\it filtered} quasi-isomorphism. 
 
 Indeed,  by   4.6.3(i), $\gr\, p_X^*$ is the direct sum of morphisms $p_{\eta_x}^* : F_{-n}(\eta_x )\to Cous (F_{-n}) (\Bbb A^1_{\eta_x} )$,  $x\in X^{(n)}$. Since $F_{-n}$ is $\Bbb A^1$-homotopy invariant,  our assertion follows from 4.6.3(ii) and the above lemma.

4.6.5. {\it Question.} Let $\{ Z_i \}$ be a smooth stratification of $X$. Can one, upon  suitable refinement of $\{ Z_i \}$, recover $M(X)$ by gluing $M(Z_i )(\codim Z_i )$ together using link correspondences and homotopies between them? Cf.~Question in 4.6.1.

For example, let $Z\subset X$ be a smooth divisor which admits a link correspondence $\theta$. Set $U:= X\smallsetminus Z$ and $\Z [U\cup_\theta Z] := \CC one (\Z [Z\times \Bbb G_m ]\to \Z [U] \oplus \Z [Z\times\Bbb A^1 ])$, the arrow is the difference of $\theta$ and the embedding $Z\times \Bbb G_m \hra Z\times\Bbb A^1$. One has a morphism $\alpha : \Z [U\cup_\theta Z]\to \Z [X]$ which is  $U\hra X$ on $\Z [U]$ and $\tilde{\theta}$ on $\Z [Z\times\Bbb A^1 ]$. 

{\it Exercise.} Show that
 $\alpha$ is a homotopy equivalence in $\CD_{\CM}^{\text{eff}}$.\footnote{Hint: use the proposition in 4.6.1 and Remark  (a) in loc.~cit.}

 \bigskip
 
\centerline{\bf \S 5. The residual proofs }

\medskip

In this section we prove the assertions from 4.6. The material of 5.1 is a version of the contents of lectures 11, 21 in \cite{MVW} and  \S\S  2, 4.1--4.3 of  \cite{Vo1}, the proposition in 4.3 is M.~Walker's theorem from  \cite{MVW} 11.17. 
The proof in 5.6 of the theorem from 4.6.3 is based on
 a variant of Gabber's argument \cite{G}.

{\bf 5.1. Splitting off a puncture on a  curve.}  Let $X$ be a smooth affine curve, $z$ a $k$-point of $X$, $U:= X\smallsetminus \{ z\}$. Suppose  $z$ is a principal divisor, so we have $\varphi\in \CO (X)$ with $div\,\varphi =z $. Consider the morphism $R [U] \to  R [X] \oplus R (1)[1]$ whose components come from the embedding $U\hra X$ and the map $U \buildrel{\varphi}\over\to \Bbb G_m$ followed by a projector to $R (1)[1]$ (see 2.2). The gist of this section is the next assertion:

{\it  $\Z^\Delta [U] \to \Z^\Delta [X]\oplus R_{\tr}^\Delta (1)[1]$ is a homotopy equivalence of premotives.
The projection $\Z^\Delta [U]\to R_{\tr}^\Delta (1)[1]$ is  left inverse to the composition $R_{\tr}^\Delta (1)[1] \to  \Z^\Delta [U/X] \to\Z^\Delta [U]$, where the first arrow is the map (4.6.1) coming from any link correspondence for $(X,z)$. }

Below we consider a  more general setting when we have a finite subset of points (instead of a single $z$),  and the picture depends on parameters.

5.1.1. For us, a {\it triple} is a datum $(X/S,U) $ formed by smooth affine schemes $U$, $X$, and $S$, an open embedding $U\hra X$, and a smooth morphism $X\to S$ of relative dimension 1, such that  $Z:= X\smallsetminus U$ is finite over $S$ (we consider $Z$ as a reduced scheme). A triple with properties (a), (b) below is called {\it Voevodsky triple}:

 (a) The diagonal divisor  $\Delta_Z : Z\hra  Z\mathop\times\limits_S X'$ is  principal. 

(b) There exists an open embedding of $S$-schemes $X\hra \bar{X}$ such that $\bar{X}$ is normal and  proper over $S$, $X^\infty := \bar{X}\smallsetminus X$ is finite over $S$, and $T:= X^\infty \sqcup Z$ admits a quasi-affine Zariski neighborhood in $\bar{X}$.
\footnote{ The datum of (b) is essentially the same as {\it  standard triple} from \cite{Vo1} 4.1, \cite{MVW} 11.5. } 

Notice that $\bar{X}$ is uniquely defined,  $Z$ is a closed subscheme of $\bar{X}$, and $T$ is finite over $S$. 

\proclaim{\quad Proposition} (i) For a Voevodsky triple $(X/S,U)$ the morphism
$\Z^\Delta [U]\to\Z^\Delta [X]$ in  is  a homotopy split surjection,  i.e., it admits a right inverse $s$ in $\CD_{p\CM}^{\text{eff\, tri}}$.\footnote{Notice that $\CC^\Delta ( \Z [U] )\to \CC^\Delta (\Z [X])$  is {\it injective} as a morphism of complexes, yet  it is a {\it split surjection} as a morphism in the homotopy category.} Therefore $\Z^\Delta [U/X]\to \Z^\Delta [U]$ is the kernel of this surjection.

(ii) For any $S$-scheme $Z/S$ the premotives   $\Z^\Delta [U /X]$ for all Voevodsky triples $(X/S,U)$ with $Z=X\smallsetminus U$ are canonically identified in the homotopy category
$ \CD_{p\CM}^{\text{eff\, tri}}$. 
\endproclaim

{\it Proof.} The promised splitting and identifications are produced by special correspondences defined as follows:

 Let $(X/S,U)$, $(X'/S,U')$ be Voevodsky triples with the same $Z/S$. A meromorphic
 function  $f$ on $X\mathop\times\limits_{S} X'$ is said to be {\it special} if for some open neighborhood $W$ of $T$ in $ \bar{X}$ (see (b) above) our  $f$ is a regular function on $W \mathop\times\limits_S X'$ whose restriction to $Z \mathop\times\limits_{S} X'$ is an
 equation of the divisor $\Delta_Z$ (see (a)), and $f$ equals 1 on $X^\infty \mathop\times\limits_{S}X'$. Such an $f$ exists due to  (a): indeed, if $W\subset \bar{X}$ is any quasi-affine neighborhood of $T$, then every function on $T\mathop\times\limits_S X'$ extends to  $W\mathop\times\limits_S X'$.\footnote{Let us check the last assertion. Since $S$ is affine, $\Spec\, \CO (W)$ is an $S$-scheme. Consider the embeddings $T \hra     W\hra \Spec\,\CO (W)$. 
    Since $T$ is finite over $S$, it is a closed subscheme of $ \Spec\,\CO (W)$, hence  $T\mathop\times\limits_S X' \hra    \Spec\,\CO (W)\mathop\times\limits_S X' $ is a closed embedding of affine schemes, q.e.d.}
Then $div (f)\subset X\mathop\times\limits_{S} X $ is finite over $X'$, i.e., $\theta_f :=div (f)\in\CC or (X',X)$. We call such a correspondence {\it $S$-special}, or simply {\it special}. 

If the two triples coincide, then one can demand that $f$, in addition, vanishes on the diagonal $\Delta (W \cap X)\subset W\mathop\times\limits_S X$. Such an $f$ exists due to (a), since every function on $T\mathop\times\limits_S X$ that vanishes on $\Delta_Z $ extends to a function on $W\mathop\times\limits_S X$ that vanishes on $\Delta (W \cap X)$.\footnote{For the closure of 
    $\Delta (W \cap X)$ in $\Spec\,\CO (W)\mathop\times\limits_S X$ intersects the closed subscheme
     $T\mathop\times\limits_S X$ (see the previous footnote) by $\Delta_Z$.}  We refer then to $\theta_f$  as a {\it very special} correspondence.

Notice that for a special correspondence $\theta$ its restriction $\theta^o$  to $U'$ lies in $\CC or (U',U)$, i.e., $\theta \in \CC or ((X',U'),(X,U))$ (see 3.4.1), so we have $\theta^{rel}: \Z [U'/X'] \to \Z[U/X]$.

\proclaim{\quad Lemma} (i) Any special $\theta  $  is $\Bbb A^1$-homotopic to zero in $\CC or (X',X)$.

(ii) If $\theta_1, \theta_2 \in\CC or (X',X)$ are special, then $\theta_1 -\theta_2 \in \CC or (X',U)$.

(iii) The composition of special correspondences is special.
\endproclaim

{\it Proof of Lemma.} (i) Use a homotopy $div (1-t + tf)\in \CC or (X'\times\Bbb A^1 ,X)$.

(ii) Notice that $f_1/f_2$ is invertible on $Z\times X'$.

 (iii) Let $(X''/S,U'')$ be another Voevodsky triple with same $Z$, $\theta_{f'} \in \CC or (X'',X')$ a special correspondence. Then $\theta_f\theta_{f'} = \theta_g$ where $g(x,x'')=f(x,\theta_{f'}(x''))$ (here $(x,x'')\in X\mathop\times\limits_S X$). Let us check that $g$ is special.
   Our $f$ is invertible  on $T\mathop\times\limits_S  (X'\smallsetminus W')$  and $X'\smallsetminus W'$ is finite over $S$, so there is an open neighborhood $W_0 \subset W$ of $T$ such that   $f$ is invertible  on $W\mathop\times\limits_S  (X'\smallsetminus W')$. Then $g $ is regular on  $W_0\mathop\times\limits_S X''$. The rest is left to the reader. \hfill$\square$

We return to the proof of the proposition. (i) If $\theta\in\CC or (X,X)$ is very special, then $s:=\id_X - \theta$ lies in $\CC or (X,U)$. By (i) of the lemma, it is right inverse to $\Z^\Delta [U]\to\Z^\Delta [ X]$ in the homotopy category, and we are done. Notice that the corresponding projector of $\Z^\Delta [U]$ onto the kernel of $\Z^\Delta [U]\to \Z^\Delta [X]$ equals the restriction $\theta^o$ of $\theta$ to $U$. Thus $\theta^{rel}\in \End (  \Z^\Delta [U/X] )$ is homotopic to the identity.

 (ii) By (ii) of the lemma, the composition $\Z^\Delta [U'/X'] \buildrel{\theta_1 -\theta_2}\over\lra \Z^\Delta [U/X]\to\Z^\Delta [U]$ vanishes. Since, by (i), $\Z^\Delta [U/X]\to\Z^\Delta [U]$ is injective in the homotopy category,  the morphism $\theta^{rel}: \Z^\Delta [U'/X'] \to \Z^\Delta [U/X]$ in the homotopy category does not depend on the choice of special $\theta$. This system of canonical morphisms $$ \Z [U'/X'] \to \Z[U/X] \tag 5.1.1$$ in the homotopy category is transitive by (iii) of the lemma. Since  $\theta^{rel}$ for $(X,U)=(X',U')$ equals identity,  (5.1.1) is a transitive system of isomorphisms, q.e.d.  \hfill$\square$

\medskip

{\it Remarks.} (i)  Voevodsky triples
 are stable with respect to the base change by any $ S' \to S$ where $S'$ is affine and smooth. The splittings and projectors constructed are compatible with the base change. 

(ii) If $Z$ is smooth, then for any $S$-special correspondence $\theta$ one has $Sp (\theta )=\id_Z$.
\medskip

5.1.2.   By  3.1.2, the above results  translate into statements about $\Bbb A^1$-homotopy invariant presheaves with transfers. For example, we have the next proposition which implies the lemma in 4.6.4: 

\proclaim{\quad Proposition}  Let $F$ be an $\Bbb A^1$-homotopy invariant presheaf with transfers. Let $\eta$ be  the generic point of a smooth variety and $X\subset \Bbb A^1_\eta $ an open subscheme. Then   $F(X)\iso F_{\text{Zar}}(X)\iso F_{\text{Nis}}(X )$ and $H^{>0}_{\text{Zar}}(X, F_{\text{Zar} })=H^{>0}_{\text{Nis}}(X, F_{\text{Nis} })=0$. 
\endproclaim

{\it Proof.}  For every open $U\subset  X$ the datum $(X/S,U)$ where $S=\eta$,   is a Voevodsky triple. If $Z:=X\smallsetminus U= \{ z_i \}$, then, by (i) of the proposition in 5.1.1, the sequence $\Z^\Delta [U/X]\to \Z^\Delta [U]\to \Z^\Delta [X]$ in $\CD_{p\CM}^{\text{eff\, tri}}$ is split exact. Thus, in particular, the map $F(X)\to F(U)$ is injective.

Let $\pi : X' \to X$ be any \' etale map such that $\pi^{-1}(Z)\iso Z$. Then $(X'/S,U')$, $U':=\pi^{-1}(U)$, is a Voevodsky triple, and
for every $a\in X (k)\smallsetminus Z$ the correspondence $\theta :=\pi - n\cdot a \in\CC or (X',X)$, where $a :X' \to X$ is a constant map with value $a$, $n$ the degree of $Z$, is special. Notice that $\theta^{rel}: \Z^\Delta [U'/X']\to \Z [U/X]$ equals $\pi$. Since $\theta^{rel}$ is a homotopy equivalence by (the proof of) (ii) in the proposition in 5.1.1,  one has $F(U)/F(X)\iso F(U')/F(X')$. Passing to the inductive limit by all such $X'$, resp.~by those $X'$ that are open subsets of $X$, we see that $F(U)/F(X)\iso H^1_Z (X_{\text{Zar}}, F_{\text{Zar}})\iso H^1_Z (X_{\text{Nis}},F_{\text{Nis}})$. This implies the assertion. 
\hfill$\square$

\medskip

{\bf 5.2. Nice pairs.} We use the notation from 3.4.1. A Voevodsky triple $(X/S,U)$ is said to be {\it nice} if $Z:=X\smallsetminus U$ is \' etale over $S$.
A {\it Voevodsky pair} is a pair $(X,U)$ (see 3.4.1) which can be included in a Voevodsky triple $(X/S,U)$. A {\it nice Voevodsky pair}, or simply {\it nice pair}, is one that can be included in a nice $(X/S,U)$. We denote by $\CC or (\CP a^{nice})$ the full subcategory of $\CC or (\CP a)$ formed by nice pairs.
Notice that for any affine smooth $Z$ the pair $Z\times \Bbb T$ is nice.

\proclaim{\quad Proposition} (i) For every Voevodsky pair $(X,U)$ the morphism $\Z^\Delta [U]\to\Z^\Delta [X]$ is a split surjection  in $\CD_{p\CM}^{\text{eff\, tri}}$  with kernel $\Z^\Delta [U/X]\to \Z^\Delta [U]$.

(ii) Every nice pair $(X,U)$ admits a link correspondence (see 4.6.1).

(iii) For every nice pair $(X,U)$ there is a canonical isomorphism $$\Z^\Delta [Z](1)[1] \iso \Z^\Delta [U/X] \tag 5.2.1$$  in $\CD_{p\CM}^{\text{eff\, tri}}$ which is the identity map for $(X,U)=Z\times\Bbb T$. For every correspondence $\theta : (X',U')\to (X,U)$ between nice pairs  it identifies $ \theta^{rel} : \Z^\Delta [U'/X'] \to \Z^\Delta [U/X]$ with $Sp (\theta )(1)$.
\endproclaim

{\it Proof.} (i) is the same as (i) of the proposition from 5.1.1. 

(ii) Since for a nice Voevodsky triple $(X/S,U)$ the triple $(Z\times \Bbb A^1 /S, /Z\times\Bbb G_m )$ is  a (nice) Voevodsky triple as well, there is 
an $S$-special $\theta \in \CC or (Z\times\Bbb T, (X,U))$ (see the proof in 5.1.1). This is a link correspondence by Remark (ii) in 5.1.1.

 (iii) It suffices to provide for each nice $(X,U)$ {\it some} isomorphism
 $\alpha_{(X,U)} :\Z^\Delta [Z](1)[1] $ $\iso \Z^\Delta [U/X]$  in $\CD_{p\CM}^{\text{eff\, tri}}$ which is identity map for $(X,U)=Z\times\Bbb T$ and such that
 for every correspondence $\theta : (X',U')\to (X,U)$ between nice pairs the map $ \alpha_{(X,U)}^{-1}\theta^{rel}\alpha_{(X',U')} : \Z^\Delta [Z'](1) \to \Z^\Delta [Z](1)$ equals $Sp (\theta )(1)$.
 Indeed,  $\alpha_{(X,U)}$ coincides then with 
 $\theta^{rel}: \Z^\Delta [Z](1)[1]\to \Z^\Delta [U/X]$ for every link correspondence $\theta$, hence is canonical. 

For a given nice $(X,U)$, choose a nice Voevodsky triple $(X/S,U)$ and an $S$-special $\rho_{(X,U)} \in \CC or (Z\times\Bbb T, (X,U))$. Set 
 $\alpha_{(X,U)}:=\rho_{(X,U)}^{rel}$. By (the proof of) (ii) in 4.5.1, this is an isomorphism  in $\CD_{p\CM}^{\text{eff\, tri}}$ whose inverse is the map $\mu_{(X,U)}^{rel}$ for any
 $S$-special $\mu_{(X,U)}\in \CC or ((X,U),Z\times \Bbb T)$. Therefore $ \alpha_{(X,U)}^{-1}\theta^{rel}\alpha_{(X',U')} = (\mu_{(X,U)} \theta \rho_{(X',U')})^{rel}$. By Remark (ii) in 5.1, one has $Sp ( \mu_{(X,U)} \theta \rho_{(X',U')})=Sp (\theta )$. We are done by  3.4.2 applied to $\mu_{(X,U)} \theta \rho_{(X',U')}\in \CC or (Z'\times\Bbb T ,Z\times\Bbb T )$.  \hfill$\square$  

\medskip

{\it Example.} Let $(X,U)$ and $(X',U')$ be nice pairs and $\pi: X'\to X$ be an \' etale map such that $\pi^{-1}(Z)=Z'$. Consider $\pi$  as a correspondence $(X',U')\to (X,U)$; then $Sp (\pi )$ equals $\pi|_{Z'} : Z'\to Z$. By the proposition,  for every link correspondences $\theta$, $\theta'$ for $(X,U)$, $(X',U')$ the diagram $$\spreadmatrixlines{2\jot} \matrix \Z^\Delta [Z'](1)[1] 
 & \buildrel{\theta^{\prime rel}}\over\lra & \Z^\Delta [U'/X' ] \\ \pi   |_{Z'} (1) \downarrow&&\pi\downarrow \\  \Z^\Delta [Z](1)[1]  & \buildrel{\theta^{rel}}\over\lra & \Z^\Delta [U/X] \endmatrix \tag 5.2.2$$  is homotopically commutative.

\medskip

{\bf 5.3.} One can use 5.1.1, 5.2 to treat every situation Zariski locally:

\proclaim{ \quad Proposition}   Let $X$ be a 
smooth quasi-projective variety, $Z\subset X$ a divisor,  $Q\subset X$ a finite set of  smooth points.  If $k$ is infinite, then $Q$ admits a Zariski neighborhood $V\subset X$ such that  $(V,V\smallsetminus Z)$ is a Voevodsky pair. 
  If $Z$ is  smooth at $Q\cap Z$, then  one can find $V$ such that $(V,V\smallsetminus Z)$ is nice.
\endproclaim

{\it Proof.} We can assume that $X$ is affine and equidimensional, and (replacing each point in $Q$ by a smooth closed point from its closure) that points of $Q$ are closed. 

Consider a closed embedding  $X\hra K$ where $K$ is a vector space, and compose it with an embedding $K\hra L$, where $L$ is another vector space such that every polynomial of degree $\le 2$ on $K$ is the restriction to $K$ of some
 linear function on $L$. 

Let $L\subset \bar{L}:= \Bbb P (L\times  \Bbb A^1 )$ be the open embedding, so $\bar{L}\smallsetminus L =\Bbb P (L)$. Let $\bar{X}$ be the closure of $X$ in $\bar{L} $; same for $\bar{Z}$, etc.

Consider a  linear projection $\pi : L\twoheadrightarrow D$, $\dim D =\dim X-1$. Let $C$ be the closure of $X\subset \bar{L}\times D$, so we have a proper projection $ C \to D$; set $X^\infty := C\smallsetminus X$.

Suppose  $\pi$ is  generic, i.e., it lies in a sufficiently small open subset in the space of all linear projections (since $k$ is infinite, we can  choose such a $\pi$ to be defined over $k$).
Then $\bar{Z}\smallsetminus Z$ does not intersect the closure of a fiber of $\pi$, hence the map $Z\to D$ is finite. Similarly, $\bar{X}\smallsetminus X$ intersects it by finitely many points, hence the projection $X^\infty \to D$ is finite.  
By a Bertini theorem of \cite{SGA4} XI  2.1,\footnote{It is here that we need to replace $K$ by $L$.} the fibers of $\pi$ passing through points of $Q$ intersect $X$  transversally, and if $Z$ is a divisor smooth at $Q\cap Z$, then these fibers intersect $Z$ transversally at smooth points of $Z$. 

Set $\tilde{Q}:= Q\cup (Z\cap \pi^{-1}\pi (Q))\subset X$; this is a finite set of  points.
Consider the line bundle $\CO (\Delta )$ on $X\mathop\times\limits_D  X$. Choose an  affine open $W\subset X\mathop\times\limits_D  X$ which contains $ \tilde{Q}\mathop\times\limits_D \tilde{Q} $  where $\CO (\Delta )$ is trivial. Let $U\subset X$ be an affine open which contains $\tilde{Q}$ and such that  $U\mathop\times\limits_S \tilde{Q} \subset W$.

Let $S\subset D$ be an affine Zariski neighborhood  of $\pi (Q)$; write $X_S$, etc., for the pull-back of our schemes to $S$. If $S$ is sufficiently small, then  $U^\infty_S := C_S \smallsetminus U_S$ is finite over $S$, $U_S^\infty \cap Z_S =\emptyset$ and  
$U_S^\infty \cup Z_S$ admits an affine neighborhood in $C_S$; if $Z$ is a  smooth at $Q\cap Z$, then $Z_S$ is \'etale over $S$. Our $V/S$ is such an $U_S$.
  \hfill$\square$

\medskip

{\bf 5.4. Proof of the proposition in 4.6.1.} By Remark (b) in 4.6.1, we can assume that $k$ is infinite.

(i) Follows from the proposition in 5.3 and (ii) of the proposition in 5.2.

(ii)  
By  5.3, one can find an open $V\subset U$ which contains the set  $\theta (P)\cup \theta' (P )$, hence $P$, such that the pair $(V,V\smallsetminus Z)$ is nice. Set $W:=\{ w\in (V\cap Z)\times \Bbb A^1 : \theta (w)\cup\theta' (w)\subset V\}$; this is an open neighborhood of $P$ in $ (V\cap Z)\times \Bbb A^1$. For an open affine neighborhood $Z'  $ of $P$ in $W\cap Z$ set $V':= W\cap (Z'\times \Bbb A^1 )$. If $Z'$ is sufficiently small, then $(V'/Z', V'\smallsetminus Z' )$ is a nice Voevodsky triple, so  $(V',V'\smallsetminus Z')$ is a nice pair. Then $\theta|_{V'},\theta'|_{V'} \in \CC or ((V',V'\smallsetminus Z'), (V,V\smallsetminus Z))$. Since $Sp(\theta|_{V'})=Sp(\theta'|_{V'})$, the morphisms $\theta^{rel} |_{V'}, \theta^{\prime rel} |_{V'} : \Z^\Delta [Z'](1)[1] \to \Z^\Delta [(V\smallsetminus Z)/V]$ are homotopic by (iii) of the proposition in 5.2. So their compositions with $\Z^\Delta [(V\smallsetminus Z)/V]\to\Z^\Delta [ (U\smallsetminus Z)/U]$ are homotopic, q.e.d.  \hfill$\square$

\medskip

{\bf 5.5. } The proof of the theorem in 4.6.3 uses another Bertini-type statement:  

  Let  $L$ be a finite-dimensional vector space. Suppose we have a datum $\{ (Y^m , S_m )\}$, 
 $m=1,\ldots , \dim L $,   where $Y^m \subset L$ is a
 closed reduced subscheme purely of codimension $m$, and $S_m \subset Y^m$ is a finite subset of closed smooth points of $Y^m$ which meets every irreducible component of $Y^m$.

A flag $\Phi =( L_1   \subset L_2 \subset \ldots  )$  of linear subspaces of $L$, $\dim L_m = m$,  is said to be {\it cotransversal} to  $\{ (Y^m , S_m )\}$ if for every $m =1,\ldots ,\dim L-1$ the map $Y^m \to L/L_m$ is finite and is  \'etale at $S_m$, and over a Zariski neighborhood of the image of $S_{m+1}$ in $L/L_m$ the projection $Y^{m+1} \to L/L_m$ is an embedding.

\proclaim{\quad Proposition} 
If $k$ is infinite, then one can find $\Phi$  cotransversal to $\{ (Y^m , S_m )\}$.
\endproclaim

{\it Proof.} We can assume that $\dim L >1$.  Let us construct a non-empty Zariski open $U\subset \Bbb P (L)$ such that for every line $\ell $ in $ U$  all the projections $Y^m \to L/\ell$ are finite, for every $m\ge 2$ the morphism $Y^m \to L/\ell$  is an embedding over a Zariski neighborhood of the image of $S_m$, and the projection $Y^1 \to L/\ell$ is \'etale at $S_1$.

Consider the open embedding $L \subset \Bbb P (L \times \Bbb A^1  )$, so $\Bbb P (L)= \Bbb P (L \times\Bbb A^1 )\smallsetminus L$. Let $\bar{Y}^m $ be the closure of $Y^m$, $Y_\infty^m := \bar{Y}^m \smallsetminus Y^m \subset \Bbb P (L)$.  Let $T^m \subset L$ be the union of tangent spaces to $Y^m$ at points of $S_m$ (considered as affine planes in $L$); set $T^m_\infty := \bar{T}^m \smallsetminus T^m  \subset \Bbb P (L)$.\footnote{I.e., $T^m_\infty =\cup_{s\in S_m}  \Bbb P (T_s )$ where $T_s$ is the tangent plane to $Y^m$ at $s\in S_m$.} 
Let $Z^m \subset \Bbb P (L) $ be the closed subvariety of those points $x$ that for some $s\in S_m$  the line  passing through $x$ and  $s$ intersects $\bar{Y}^m $ {\it not} only at $s$. Since $\dim Z^m \le \dim Y^m$ and $\dim Y^m_\infty \cup T^m_\infty <\dim Y^m$,  the complement  to 
$(\mathop\cup\limits_{m\ge 2} Z^m ) \cup Y^{1}_\infty \cup T^{1}_\infty  $ in $\Bbb P (L)$ is  non-empty; this is our $U$. 

We construct  $\Phi =(L_1 \subset L_2 \subset \cdots )$ by induction by $\dim L$. Take for $L_1$ any element of $U(k)$ (which is non-empty since $k$ is infinite). Let $ (L_2 /L_1 \subset L_3 /L_1 \subset \ldots )$ be any flag in $L/L_1$   cotransversal to the datum of the images of $\{ (Y^m ,S_m )\}$, $m\ge 2$, in $L/L_1$; one can find it by the induction assumption.
 \hfill$\square$

\medskip

{\bf 5.6. Proof of the theorem in 4.6.3.}  As in loc.~cit.,  $k$ is perfect and $F$ is an $\Bbb A^1$-homotopy invariant presheaf with transfers. By an argument similar to that in Remark (b) in 4.6.1, we can assume that $k$ is infinite.

5.6.1.  Let us check that the Cousin differential  commutes with the pull-back  map $f^*$ for a smooth $f: Y\to X$. The assertion is local, so we can assume that $f$ is the composition of an \' etale map and a direct product projection; then it suffices to treat  either of these types of maps separately. Replacing $X$ by the smooth part of the normalization of its reduced  irreducible subscheme of codimension $n$, and $F$ by $F_{-n}$, we see that it suffices to consider $d: \CC ous (F)^0 \to \CC ous (F)^1$. We need to check that $f$ commutes with the link correspondences at the generic point of a divisor on $X$ and that of its preimage on $Y$. The \' etale case follows from (5.2.2). For a direct product projection $Y=P\times X \to X$ we  choose a link correspondence on $Y$ as the direct product of that on $X$ and $\id_P$; then the assertion is evident.

5.6.2. Let us check that for a {\it proper} $g: Z\to X$ the Cousin differential commutes with the trace map $g_*$.
Looking at  the components $d_{xy}$ (and replacing $F_{-n}$ by $F$), we see that the assertion amounts to the next lemma:

\proclaim{\quad Lemma} (i) Let $T$, $T'$ be semi-local schemes as in 4.6.2, and $ T'\to T$ be a finite  surjective morphism. Then $tr_{s_{T'}/s_{T}}\Res^{T'} = \Res^{T} tr_{\eta_{T'}/\eta_{T}} : F(\eta_{T'} )\to F_{-1}(s_{T})$.

(ii) Let $\eta$ be the generic point of a smooth variety, $C/\eta$ a proper irreducible curve, and  $Z\subset C$ a reduced subscheme finite over $\eta$ such that $U:= C\smallsetminus Z$ is smooth. Then the composition $F(U)\buildrel{Res_Z}\over\longrightarrow F_{-1}(Z) \buildrel{tr_{Z/\eta}}\over\longrightarrow F_{-1}(\eta )$ vanishes.
\endproclaim

{\it Proof.} (i) Passing to the normalizations (as in 4.6.2), we can assume that  $T$ and $T'$ are regular. 
One can find a finite morphism $\pi : V' \to V$ of smooth affine varieties and a smooth divisor $Z\subset V$ with $Z' := \pi^{-1}(Z)_{red}$  smooth, such that $T'/T$ is the localization of $V'/V$ at the generic point(s) $\eta_Z$ of $Z$. Thus $\eta_Z = s_T$, $\eta_{Z'}=s_{T'}$.

Set $\pi_s :=\pi|_{Z'}: Z' \to Z$.  The transposed correspondence $\pi^\vee \in \CC or (V,V')$ can be seen as a morphism of pairs $(V,V\smallsetminus Z) \to (V',V'\smallsetminus Z')$; then $Sp(\pi^\vee )=\pi^\vee_s$.

Take any element of $ F(\eta_{T'} )$. Shrinking $V$, we can assume that it comes from some $\phi \in F(V'\smallsetminus Z' )$. One can find an open $U' \subset V'$  that contains $\eta_{Z'}$ such that the pair $(U',U'\smallsetminus Z')$ is nice, and an open $U\subset V$ with $\eta_Z \subset U$ and $\pi^{-1}(U)\subset U'$ such that the pair  $(U, U\smallsetminus Z)$ is nice (see 5.3). Notice that $\pi^\vee|_U \in \CC or ((U, U\smallsetminus Z),(U',U'\smallsetminus Z'))$. 

Let $\theta$, $\theta'$ be link correspondences for $(U,U\smallsetminus Z)$ and $(U',U'\smallsetminus Z')$ (see 5.2(ii)). 
Set $\alpha := \pi^\vee \theta , \beta :=\theta' (\pi_s^\vee \times\id_{\Bbb T})\in \CC or (Z\times\Bbb T, (U',U'\smallsetminus Z'))$. Then $Sp(\alpha )=\pi_s^\vee =Sp (\beta )$, hence, by 5.2(iii), $\alpha^{rel}=\beta^{rel}$.  Since
$ \Res^{T} tr_{\eta_{T'}/\eta_{T}}(\phi )= \alpha^{rel*} (\phi )$,  
$tr_{s_{T'}/s_{T}}\Res^{T'} (\phi )= \beta^{rel*} (\phi )$, we are done.

 (ii) Choose a finite map $C\to \Bbb P^1_\eta$ which maps $Z$ to $\infty$. Applying (i) to the localizations of $C$ and $\Bbb P^1_\eta$ at $Z$ and $\infty$ respectively, we reduce our statement to the case $C=\Bbb P^1_\eta$, $Z=\infty$. Here it is evident since $F(U)=F(\eta)$.   \hfill$\square$
 
 \medskip

5.6.3.  Let us check that   the square of the Cousin differential vanishes.  The assertion is local, so we consider an affine  variety $X$ and $\phi \in \CC ous (F)^n (X)$; we want to check that $d_X^2 (\phi )=0$. Let $Y\subset X$ be the support of $\phi$ and $\xi$ be the union of the generic points of the support of $d^2 ( \phi)$.

Let $X\hra T$ be any closed embedding  in a vector space.  By 5.5, one can find a linear projection $T \twoheadrightarrow L$,  such that  $Y\to L$ is finite and surjective and $\xi \to L$ is an embedding. 
By (i) in the lemma in 5.6.2 applied to $Y/L$ and $F_{-n}$, we see that $d^2_X (\phi )\in F_{-n-2}(\xi )$ equals $d^2_{L} tr_{Y/L} (\phi )$. So, replacing $F_{-n}$ by $F$, $X$ by $L$, and $ tr_{Y/L} (\phi )$ by $\phi$, we are reduced to the case $\phi \in \CC ous (F)^0 (L)=F(\eta_L )$.

Let $Y^1 \subset L$ be a reduced hypersurface such that $\phi$ comes from some $\phi\,\tilde{} \in F (L\smallsetminus Y^1 )$, and $Y^2 \subset Y^1$ be a codimension 1 reduced subvariety such that $Y^1 \smallsetminus Y^2$ is smooth and $d(\phi ) \in F_{-1}(\eta_{Y^1} )$ comes from some $d(\phi )\,\tilde{}\in F_{-1}(Y^1 \smallsetminus Y^2 )$. By 5.5, we can find a decomposition $L=  \Bbb A^1 \times D$ such that $Y^1$ is finite over $D$ and $\eta_{Y^2} \to D$ is an embedding. By (i) of the lemma in 5.6.2 applied to $Y^1 /D$ and $d (\phi )\in F_{-1}(\eta_{Y^1})$, one has $d^2_L (\phi )= d_{D}tr_{Y^1 /D}d_L (\phi )\in F_{-2}(\eta_{Y^2})$.

Consider the open embedding $L= \Bbb A^1 \times D \subset \bar{L}:= \Bbb P^1 \times D$. Then $d_{\bar{L}} (\phi )= d_L (\phi )+ \Res_\infty (\phi )$ where $\Res_\infty$ is the residue at $L^\infty :=  \bar{L}\smallsetminus L= D$. By (ii) of the lemma in 5.6.2 applied to $\bar{L}/D$, one has $
tr_{Y^1 /D}d_L (\phi )= -\Res_\infty (\phi )\in F_{-1}(\eta_{D})$. Now 
$\Res_\infty (\phi )$ is the restriction to the generic point of $D=L^\infty$ of  $\Res_D (\phi\,\tilde{}\, )\in (F_{-1})_{\text{Zar}}(D)$ (see the remark in 4.6.2). Therefore
$d_{D} \Res_\infty (\phi )=0$,  q.e.d.

\medskip

5.6.4.  Let us check that  for a smooth semi-local  $\hat{X}$ one has $H^{>0} \CC ous (F)(\hat{X})$ $=0$. Our $\hat{X}$ is localization of some smooth affine variety $X$ at finitely many points, so $\hat{X}$ is the intersection of semi-local schemes $X_Q$,  where $Q\subset X$ are some finite collections of {\it closed} points. Thus $\CC ous (F)(\hat{X}) =\limright \CC ous (F)(X_Q )$, and we can assume that $\hat{X}=X_Q$. 

Let $\phi \in \CC ous (F)^m (X_Q )$, $m>0$, be any cycle. Since $\CC ous (F)(X_Q )$ is the inductive limit of the Cousin complexes of the Zariski neighborhoods of $Q$, we can assume (shrinking $X$ if necessary) that $\phi$ is a cycle in $\CC ous (F)(X)$. 
Let $Y\subset X$ be the support of $\phi$. We look for an open neighborhood $U$ of $Q$ such that $\phi |_U$ is exact.

Let $X\hra T$ be a closed embedding into a vector space. By 5.5, one can find a linear projection $T\twoheadrightarrow L$ such that the projection  $\pi : X \to L$ is finite, surjective, and \'etale at $Q$, and $\eta_Y \to L$ is an embedding. Consider the cycle $tr_{X/L} (\phi )\in \CC ous (F)^m (L)$. Suppose we know that it is exact, i.e., $tr_{X/L} (\phi )=d_L (\psi )$. One can find an open neighborhood $U\subset
X$ of $Q$ such that $U/L$ is \'etale and $\pi^{-1}\pi (\eta_Y )\cap U=\eta_Y$. Set $\psi_U := (\pi |_U )^* (\psi )\in \CC ous (F)^{m-1}(U)$. Then $d_U (\psi_U )= \phi |_U$; we are done.

It remains to show that $H^{>0} \CC ous (F)(L)=0$. We do this by induction by $\dim L$. Write $L=\Bbb A^1 \times D$ and consider the projection $p : L\twoheadrightarrow  D$. We know that 
$H^{>0} \CC ous (F)(D)=0$ by the induction assumption, so it remains to check that $p^* : \CC ous (F)(D)\to\CC ous (F)(L)$ is a quasi-isomorphism. Filtering the Cousin complexes as in 4.6.4, we see that $p^*$ is, in fact, a filtered quasi-isomorphism (since for any $\xi\in D$ the 2-term complex $\CC ous (F)(\Bbb A^1_\xi )$ is a resolution of $ F_{-m}(\xi )=F_{-m}(\Bbb A^1_\xi )$, say, by  5.1.2), q.e.d.

\medskip

5.6.5. It is clear that the image of the map $F(X)\to \CC ous (F)^0 (X)$ is killed by the Cousin differential.
Let us show that for $\hat{X}$ as in 5.6.4, $F(\hat{X})\iso H^0 \CC ous (F)(\hat{X} )$.

Our map is injective: this follows from the proposition in 5.3 and (i) of the proposition in 5.1.1. 
 It remains to show that any $\phi_\eta \in H^0 \CC ous (F)(\hat{X} )$ comes from $F(\hat{X})$. We can assume that $\hat{X}$ is the localization of a smooth affine $X$ at a finite subset $Q\subset X$, and $\phi_\eta \in H^0 \CC ous (F)(X)$. There is a divisor $Z\subset X$ such that $\phi_\eta$ comes from some  $\phi \in F(X\smallsetminus Z)$. Replace $X$ by $V$ from the proposition in 5.3. Choose a correspondence $s$ from (i) of the proposition in 5.1.1; its action $s^*$ is left inverse to the restriction morphism $F(V)\to F(V\smallsetminus Z )$. So we have $s^* (\phi )\in F(V)$. Let us check that $\phi_\eta $ equals $ s^* (\phi )_\eta$ (the restriction of $s^* (\phi )$ to $\eta =\eta_X$); this will prove our assertion.
 
Replacing $\phi$ by  $ \phi - s^* (\phi )|_{V\smallsetminus Z}$, we can assume that $s^* (\phi )=0$;  we want to check  that $\phi_\eta =0$.  Let $\kappa$ be the generic point subscheme of the base  of the projection $V\to S$. Pulling back   $(V/S, V\smallsetminus Z)$ to $\kappa$, we get a $\kappa$-curve $V_{\kappa}$, a divisor $Z_{\kappa} \subset V_{\kappa}$,   $\phi_{\kappa}\in F(V_{\kappa} \smallsetminus Z_{\kappa})$, and, by Remark (i)  in 5.1.1,
 $s_{\kappa} \in \CC or (V_{\kappa},  V_{\kappa} \smallsetminus Z_{\kappa})$ such that $s^*_{\kappa}(\phi_{\kappa})=0$. By (i) of the proposition in 5.1.1 and 
5.2, we have a  decomposition $F(V_{\kappa} \smallsetminus Z_{\kappa})\iso F(V_{\kappa})\oplus F_{-1}(Z_{\kappa})$ where the first projection is $s_{\kappa}^*$ and the second one is $\Res_{Z_{\kappa}}$ (see 4.6.2).  Since $\phi$ is killed by the Cousin differential, $\Res_{Z_{\kappa}} (\phi_{\kappa} )=0$. Thus $\phi_{\kappa}=0$, so $\phi_\eta =0$, q.e.d.  \hfill$\square$

 \bigskip
 
\centerline{\bf \S 6. Playing with motives: first results}

\medskip
 
 We assume that the base field $k$ is perfect.  This section corresponds roughly to the material of lectures 3, 5.9, 12.21--12.28, 13.22--13.27, 16.13--16.25, 20 in
  \cite{MVW},  3.5,  4.1, 4.3 of \cite{Vo2}, \cite{Deg1}, and \S 6 of \cite{Bo1}.
 
\medskip

{\bf 6.1. The DG category of motives.} Theorem 3.3 remains valid for motives: 

\proclaim{\quad Proposition} The Tate motive is homotopically quasi-inversible: on $\CD_\CM^{\text{eff}}{}^{\text{tri}}$ (or $\mathop{\vtop{\ialign{#\crcr
  \hfil\rm $\CD^{\text{eff}}_{\CM}$\hfil\crcr
  \noalign{\nointerlineskip}\rightarrowfill\crcr
  \noalign{\nointerlineskip}\crcr}}}{\!}^{\text{tri}}$) the Tate twist endofunctor  is  fully faithful.
\endproclaim

{\it Proof.} For $F\in D\CP\CS h_{\tr}$ consider a morphism  $\mu_F : F\to (\CC^\CM (F(1)))_{-1}[-1]$ that comes from the canonical morphism $F(1)\to \CC^\CM (F(1))$ by adjunction. As in the remark in 3.3,  it suffices to check that $\mu_F$ is a quasi-isomorphism for $F\in I_{\tr}^{\Delta \bot}\cap I_{\tr}^{\text{Nis}\bot}$. 
By (4.4.3), this  follows from 3.3 and the next lemma:

 \proclaim{\quad Lemma} For $G\in I^{\Delta\bot}_{\tr}$ the evident arrow $ \CC^{\text{Nis}}(G_{-1})\to  (\CC^{\text{Nis}} (G))_{-1}$ is a quasi-isomorphism. 
 \endproclaim
 
 {\it Proof of Lemma.} We compute the Nisnevich cohomology using the Cousin resolution (see 4.6.3).
 The pull-back for the projection $\Bbb G_m \times X \to X$ yields an embedding
 $\CC ous (G)(X)\hra \CC ous (G)(\Bbb G_m \times X )$.  Consider
 the natural morphism  $\CC ous ( G_{-1} )(X)\to \CC ous (G)( \Bbb G_m \times X)/\CC ous (G)(X )$. Both  complexes are filtered:  the first one by the stupid filtration, the second one by the codimension of the image in $X$ of the support.  Our arrow is a filtered quasi-isomorphism 
 by the lemma in 4.6.4 and the theorem in 4.6.3 (cf.~4.6.4), hence  is a quasi-isomorphism, q.e.d.
\hfill$\square$

\proclaim{\quad Corollary} For  $a,b \ge 0$ there is a canonical homotopy equivalence in $D\CP\CS h_{\tr}$ $$\CC^\CM  (R_{\tr} (a) )\iso (\CC^\CM (R_{\tr} (a+b)))_{-b}[-b]. \tag 6.1.1$$ If $c>a$, then $(\CC^\CM (R_{\tr} (a)))_{-c}=0$. Same is true for $\CC^\CM$ replaced by $\CC^\Delta$.  \hfill$\square$
\endproclaim

  {\it Definition.} The  DG category $\CD_\CM$ of motives is the localization of $\CD_\CM^{\text{eff}}$ by the Tate motive. Explicitly, any motive is represented as $M(a)$, $M\in\CD_\CM^{\text{eff}}$, $a\in\Bbb Z$, and $\Hom_{\CD_\CM }(M(a), N(b)):= \limright \Hom_{\CD_\CM^{\text{eff}}} (M (a+n ), N(b+n ))$; here the inductive limit is taken for $n\to +\infty$ and  $a+n ,b+n $ are assumed to be $\ge 0$. 
 
 The  homotopy category $\CD_\CM^{\text{tri}}$ is naturally a tensor category: the tensor product is  $M(a)\otimes N(b):= ( M\otimes N )(a+b)$, and the commutativity constraint is defined using (ii) of the proposition in 3.2.
 
 {\it Remarks.} (i) A positive answer to the question in 3.2 would permit us to lift the above tensor structure to a homotopy tensor structure on the DG category $\CD_\CM$.
 
(ii) The above construction make sense for premotives as well.

\medskip

{\bf 6.2. The absolute motivic cohomology.} 6.2.1. For $X\in\CS m$ and $n\ge 0$ set\footnote{See 2.3 for the notation.} $$R\Gamma_\CM (X, R (n)):= \Hom_{\CD_\CM^{\text{eff}}} (M(X),R_\CM (n))=\CC^\CM (R_{\tr} (n))(X)  \tag 6.2.1$$ and  $H^{m,n}_\CM (X, R ):= H^m R\Gamma_\CM (X, R (n))$.
So  $R\Gamma_\CM (X, R (* ))= R\Gamma (X_{\text{Zar}} ,R^\Delta_{\tr} (* ))= R\Gamma (X_{\text{Nis}} ,R^\Delta_{\tr} (* ))$ by 4.4, 4.5. For a morphism of smooth varieties $Y\to X$ we have $R\Gamma_\CM (Y/X, R (n)):= \CC one (  \CC^\CM (R_{\tr} (n)) (X)\to  \CC^\CM (R_{\tr} (n)) (Y))$, etc.

{\it Example.} One has $\CC^\CM (R_{\tr} (0))=R$, so $R\Gamma_\CM (X,R (0))= R^{\pi_0 (X)}$.  According to 3.2(i), 4.4, one has  $\CC^\CM (R_{\tr} (1))=\CO^{\times }_{\text{Nis}}\mathop\otimes\limits^L R[-1]$, so $$R\Gamma_\CM (X,R (1))= R\Gamma (X_{\text{Zar}}, \CO^\times )\mathop\otimes\limits^L R[-1]= R\Gamma (X_{\text{Nis}}, \CO^\times )\mathop\otimes\limits^L R[-1]. \tag 6.2.2$$

As in 3.2, $R^\Delta_{\tr} (* ):= \mathop\oplus\limits_{n\ge 0} R^\Delta_{\tr}  (n)$  is a commutative unital algebra in the tensor category $ 
\mathop{\vtop{\ialign{#\crcr 
  \hfil\rm $\CD^{\text{eff\, tri}}_{p\CM}$\hfil\crcr
  \noalign{\nointerlineskip}\rightarrowfill\crcr
  \noalign{\nointerlineskip}\crcr}}} $. Therefore   $\CC^\CM (R_{\tr} (* ))$ is a commutative unital algebra in $D\CP\CS h_{\tr}$ (see Remark (ii) in 2.3), hence  in $D\CP\CS h$ (by
 the remark in 2.2). Thus $R\Gamma_\CM (X, R (* )) = \oplus R\Gamma_\CM (X,R (n))$ is a  commutative unital  $R$-algebra (equipped with an extra grading by $n$) in the tensor category $D(R)$.

{\it Remark.} The positive answer to Question in 3.2 would imply that $R\Gamma_\CM (X, R (* ))$ can be naturally realized as a homotopy object of the category of $E_\infty$-algebras.

6.2.2. The  motivic cohomology product can  also be seen in the next two ways:

(a)  The motive $M(X)$ is  naturally a cocommutative counital coalgebra in $\CD^{\text{eff\, tri}}_\CM$:  the coproduct $\delta :M(X)\to M(X)\otimes M(X)$ comes from the diagonal $X\to X\times X$, the counit is the augmentation morphism. Now the product of $f_i \in  H^{m_i , n_i }_\CM (X, R )=\Hom_{\CD_\CM^{\text{tri}}} (M(X),R_\CM (n_i )[m_i ])$, $i=1,2$, is  $ f_1 \cup f_2 = \pi (f_1 \otimes f_2 )\delta$ where $\pi$ is the product on the algebra $R_\CM (* ):= M(R_{\tr}(*))$ in $\CD^{\text{eff\, tri}}_\CM$. 

We see that if $P\in   \CD^{\text{eff\, tri}}_\CM $ is an $M(X)$-comodule, then $P(*)[\cdot ]$ is naturally an $H^{\cdot ,*}_\CM (X,R)$-module:  $\alpha \in H^{m,n}(X,R)$ acts as the composition $P\to M(X)\otimes P \buildrel{\alpha\otimes\id_P}\over\lra P(n)[m]$. This is the  {\it $\cap$-product} action. 

(b)   For any $N\in \CD_\CM^{\text{tri}}$ the tensor product $M(X)\otimes N$ is naturally a counital $M(X)$-comodule. For $M(X)$-comodules $P$, $Q$ in   $ \CD_\CM^{\text{tri}} $ let $\Hom^{M(X)}(P,Q)\subset \Hom (P,Q)$ be the subgroup of the $M(X)$-comodule morphisms. If $Q=M(X)\otimes N$, then the counit  $M(X)\to R$ defines a projection $\Hom (P,M(X)\otimes N)\to \Hom (P,N)$ which yields an identification $\Hom^{M(X)} (P,M(X)\otimes N)\iso \Hom (P,N)$.
 Similarly,  for $R_\CM (* )$-modules $K$, $L$ we have the subspace $\Hom_{R_\CM (* )} (K,L)\subset \Hom (K,L)$ of $R_\CM (* )$-module morphisms; if $K=T(* ):= R_\CM (* )\otimes T$ is an induced $R_\CM (* )$-module, then the unit $R_\CM \to R_\CM  (* )$  defines a projection 
 $\Hom (T (* ) ,L)\to \Hom (T,L)$ which yields an identification
  $\Hom_{R_\CM (* )} (T (* ) ,L)\iso \Hom (T,L)$.  Consider now the  subalgebra  $\End^{\cdot M(X)}_{R_\CM (* )} (M(X)(*))\subset \End^\cdot (M(X)(* ))$ of  $M(X)$-comodule and $R_\CM (* )$-module endomorphisms. So the $\cap$-product action is  an isomorphism of graded algebras $$
 H^{m,n}_\CM (X, R)\iso\Hom^{ M(X)}_{R_\CM  (* )}(M(X)(*),M(X) (*+n )[m]). \tag 6.2.3$$ 

For any map $f: Y\to X$ we have the algebra morphism $ f^* : H^{\cdot ,*}_\CM (X, R)\to H^{\cdot, *}_\CM (Y, R)$, and
 the morphism $M(Y)(* )\to M(X)(* )$ is a morphism of $H^{\cdot,*}_\CM (X,R )$-modules. Since $M(Y)\to M(X)$ is a morphism of coalgebras, $M(Y)$ and 
 $M(Y/X):= \CC one (M(Y)\to M( X))[-1]$ are $M(X)$-comodules,  hence $M(Y)(*)[\cdot ]$,  $M(Y/X)(* )[\cdot ]$
are $H^{\cdot,*}_\CM (X,R )$-modules. As in the previous paragraph, one has a natural isomorphism of  $H^{\cdot,*}_\CM (X,R )$-modules $$ H^{m,n}_\CM (Y/X ,R )\iso \Hom^{M(X)}_{R_\CM (* )} ( M(Y/X)(* ), M(X)(*+n )[m]). \tag 6.2.4$$ 

{\it Remarks.} (i) For $P$, $N$ as above the identification $\Hom (P,N)\iso \Hom^{M(X)} (P, $ $M(X)\otimes N)$ is $\varphi\mapsto (\id_{M(X)}\otimes\varphi )\delta_P$ where $\delta_P$ is the $M(X)$-coaction on $P$. If $P=M(X)$, $N=M(Y)$, then for morphisms given by correspondences it can be described as a map $ \CC or (X,Y)\hra \CC or (X, X\times Y)$ which assigns to a cycle $c$ on $X\times Y$ its image by $\Delta_X \times \id_Y : X\times Y  \hra X \times (X\times Y)$.

(ii) Replacing  motives by premotives, we get  premotivic cohomology complexes  which map naturally to the motivic ones. The above discussion remains valid in this context.

(iii)  If $a>0$, then $\Hom (M(X),R (- a ))=0$ (see the corollary in 6.1).

\medskip

{\bf 6.3. The Gysin  equivalence.} We refer to \cite{Deg1} for a detailed exposition.  

6.3.1. Let  $i: Z\hra X$ be a closed codimension $n$ embedding   of smooth varieties. Set $M(X, X\smallsetminus Z):= M((X\smallsetminus Z)/X) [1]= \CC one (M(X\smallsetminus Z)\to M(X))\in \CD_\CM^{\text{eff}}$.

\proclaim{\quad Proposition} There is a canonical homotopy equivalence $$i^{Gys} : M(X, X\smallsetminus Z) \iso M(Z)(n)[2n] . \tag 6.3.1$$ 
\endproclaim

{\it Remarks.} (i)  Were the question in 4.6.5 to have a positive answer, this would be immediate.

(ii) By 4.4, 4.5, the proposition means that for any $\Bbb A^1$-homotopy invariant complex  $F$ of presheaves with transfers there is a natural quasi-isomorphism $$ i^{Gys *}   \colon R\Gamma (Z, F_{-n})[-n] \iso R\Gamma_Z (X,F).
 \tag 6.3.2$$ Here the cohomology is that of associated  Zariski or  Nisnevich sheaves. Notice that if $F$ is a single presheaf, then such a quasi-isomorphism is provided by the morphism $i_*$ of the Cousin resolutions (see 4.6.3). One checks easily that the two quasi-isomorphisms coincide (see Exercise (i) below).

{\it Proof.} Let $N$ be the normal bundle to $Z$, and  $Z\hra N$ be the zero section. Now (6.3.1) is the composition of  two canonical homotopy equivalences to be defined, respectively,   in 6.3.2 and 6.3.3 below:  $$M(X, X\smallsetminus Z)\iso M(N, N \smallsetminus Z )\iso M( Z)(n)[2n].  \tag 6.3.3$$

6.3.2.  For any scheme $S$ we write $X_S $ for $X\times S$ considered as an $S$-scheme; as always,  $\Bbb G_m := \Bbb A^1 \smallsetminus \{0\}\subset \Bbb A^1$.
Let $Q=Q_t$ be the deformation to normal cone for $Z\subset X$ (see e.g.~\cite{F} ch.~5). So $Q$ is a smooth $\Bbb A^1$-scheme together with  identifications $Q|_{\Bbb G_m} \iso X_{\Bbb G_m}$, $Q_0 \iso N$, and  a closed embedding $Z_{\Bbb A^1} \hra Q$ 
which are compatible in the evident manner. Thus $Q_1 =X$, so the embeddings $Q_1 \hra Q  \hookleftarrow Q_0$ yield morphisms of motives $$M(X, X\smallsetminus Z)\to M(Q, Q\smallsetminus Z_{\Bbb A^1} ) \leftarrow M(N, N \smallsetminus Z  ).  \tag 6.3.4$$

 We will see in a moment that these  are homotopy equivalences.  The first morphism in (6.3.3) is defined then as their composition.

Let us show that the  arrows in (6.3.4)  yield quasi-isomorphisms between the $\Hom$ complexes with values in  any motive. This means that for any $\Bbb A^1$-homotopy invariant complex  $F$ of presheaves with transfers 
the pull-back maps $$R\Gamma_Z (X,F)\leftarrow R\Gamma_{Z_{\Bbb A^1}}(Q,F)\to R\Gamma_Z (N ,F) \tag 6.3.5$$ are quasi-isomorphisms; here  the cohomology are taken with respect to either the Zariski or the Nisnevich topology (see 4.4, 4.5). The topologies have finite cohomological dimension, so we can assume that $F$  is a single sheaf.   

For any closed embedding $i : T\hra Y$ of smooth varieties of codimension $n$ one has a canonical identification $i_* : R\Gamma (Y , F_{-n}) \iso R\Gamma_T (Y, F)[n]$, see the theorem in 4.6.3 (we represent both complexes  by means of the Cousin resolutions). The next general lemma (applied to $Y= Q$, $T= Z_{\Bbb A^1}$, and $D= Q_0 ,Q_1$) together with 4.5 (applied to $F_{-n}$ and $X=Z$) imply that the arrows in (6.3.5) are, indeed, quasi-isomorphisms:

\proclaim{\quad Lemma} Let $i : T\hra Y$ be as above, $r: D\hra Y$ be a smooth divisor transversal to $T$; denote by $i' $, $r'$  the embeddings  $ T\cap D  \hra D$, $T\cap D \hra T$. Then the diagram  
$$\spreadmatrixlines{2\jot}
\matrix
R\Gamma_T (Y, F)[n]&
\buildrel{r^*}\over\lra & R\Gamma_{T \cap D} (D, F)[n]
 \\
i_* \uparrow&&i'_* \uparrow \\
R\Gamma (T, F_{-n})  &
    \buildrel{r{\prime *}}\over\lra & R\Gamma (T\cap D, F_{-n}), 
\endmatrix
\tag 6.3.6$$
where the horizontal  arrows are the pull-back maps, commutes.
\endproclaim

{\it Proof of Lemma.} Since $Ri^! F_Y [n]= H^n R i^! F_Y$, $R i^{\prime !} F_D [n]=H^n Ri^{\prime !}F_D$, (6.3.6) comes by applying the functor $R\Gamma (T,\cdot )$ to the diagram 
$$\spreadmatrixlines{2\jot}
\matrix
i^* H^n Ri^! F_Y&\buildrel{r^*}\over\lra & i^* r_* H^n Ri^{\prime !} F_D
\\ i_* \uparrow && i'_* \uparrow \\
F_{-n\, T} & \buildrel{r^{\prime *}}\over\lra &  r'_* F_{-n T\cap D}
\endmatrix \tag 6.3.7$$
of {\it sheaves} on $T$, and we will check the commutativity of (6.3.7). Since $i^{\prime}_*$ is an isomorphism and, by (ii) of the theorem in 4.6.3, every local section of $F_{-n T\cap D}$ is uniquely determined by its value at the generic point $\eta_{T\cap D}$, it suffices to check the assertion on a neighborhood of $\eta_{T\cap D}$.

Fit $T$ into a flag of smooth subvarieties $T=T^n \buildrel{i_n}\over\hra \ldots  \buildrel{i_2}\over\hra T^{1} \buildrel{i_1}\over\hra T^0 =Y$ such that $\codim\,  T^a = a$ and $T^i$'s are transversal to $D$. Then $i_* : F_{-n T}\to H^n Ri^! F_{Y}$ can be written as the composition of the !-pull-back to $T$ of the  arrows $i_{a*}  : F_{-a T^a} \iso Ri_a^! F_{1-a T^{a-1}}[1]$ shifted by $n-a$. Same holds if we replace $i$ by $i'$. Thus our statement for $i$ follows if we prove it for each $i_a$ (and $F_{a}$). We are reduced to the situation when $T$ is a divisor, i.e., $n=1$.

Shrinking $X$, we can assume that the pair $(X,X\smallsetminus T)$ is nice; let $V$ be a neighborhood of $\eta_{T\cap D}$ in $D$ such that the pair $(V,V\smallsetminus T)$ is also nice. If $\theta$, $\theta'$ are link correspondences for these pairs, then $Sp( \theta r') =ir' =ri'=Sp ( r\theta')$, hence, by 5.2.1(iii), $\theta^{rel} r'=r\theta^{\prime rel}$. Since $\theta^{rel}$, $\theta^{\prime rel}$ are maps inverse to the vertical arrows in (6.3.7), we are done. \hfill$\square$

6.3.3. Let us construct  the second arrow in (6.3.3). We can forget about $X$, so $N$ is any vector bundle of rank $n$ on a smooth $Z$. Let $Z\hra N$ be its zero section, $[Z]$ be the set of connected components of $Z$. First, we define  a natural quasi-isomorphism
 $$R^{[Z]} \iso  R\Gamma_{\CM }(N, N \smallsetminus Z, R (n))[2n]  \tag 6.3.8$$
compatible with the base change.  It suffices to define (6.3.8) locally with respect to $Z$ (in a way compatible with the base change). So we can assume that $N$ is a trivial vector bundle. Choose a trivialization $\Bbb A^n_Z \iso N$. Consider a covering $\{ U_i \}$ of $\Bbb A^n \smallsetminus \{ 0\}$,  
$U_i := \Bbb A^{i-1}\times\Bbb G_m \times \Bbb A^{n-i} $, and its  \v Cech complex $V_\cdot$. We get identifications $M(\Bbb A^n ,\Bbb A^n \smallsetminus \{ 0\})\buildrel{\sim}\over\leftarrow M(\Bbb A^n , V_\cdot )= M( \Bbb A^1 ,\Bbb A^1 \smallsetminus \{ 0\})^{\otimes n}= (R_\CM (1)[2])^{\otimes n}= R_\CM (n)[2n]$. So the trivialization of $N$ yields an identification $M(N,N\smallsetminus Z)\iso M(Z)\otimes M(\Bbb A^n ,\Bbb A^n \smallsetminus \{ 0\})\iso M(Z)(n)[2n]$. We define (6.3.8) as the composition $
R^{[Z]} \iso \Hom (M(Z),R_\CM )\iso\Hom (M(Z)(n)[2n], R_\CM (n)[2n])\iso \Hom (M(N,N\smallsetminus Z),
R_\CM (n)[2n])=R\Gamma_{\CM }(N, N \smallsetminus Z, R (n))[2n]$.
 Its independence of the trivialization follows since $GL(n)$ is connected. 

The generator $1^{[Z]} \in R^{[Z]}$ can be seen, by (6.3.8) and  (6.2.4), as a  morphism $M(N, N\smallsetminus Z )\iso M((N\smallsetminus Z )/Z )[1]\to M(Z)(n)[2n]$ of $M(Z)$-comodules. This is the second arrow in (6.3.3). A moment ago we have checked  that 
it is a homotopy equivalence when $N$ is trivial; the general case reduces to this if we
replace $Z$ by a hypercovering $U_\cdot$ such that $N$ is trivial on $U_0$. \hfill$\square$

\proclaim{\quad Corollary} Let $F$ be an $\Bbb A^1$-homotopy invariant complex of presheaves with transfers and  $X$ be a smooth variety. Then for the coniveau filtration\footnote{I.e., the filtration by the codimension of support.}  on $R\Gamma (X, F)$ one has a canonical quasi-isomorphism $\gr^n R\Gamma (X, F)\iso \mathop\oplus\limits_{x\in X^{(n)}} F_{-n}(\eta_x )[-n]$.  \hfill$\square$
\endproclaim

6.3.4. {\it Exercises.} (i) If the above $F$ is a single sheaf, then the filtered complex  $R\Gamma (X, F)$ identifies canonically with  $\CC ous (F)(X)$ equipped with the stupid filtration, and isomorphism (6.3.2)
coincides with the identification $i_*$ from the theorem in 4.6.3 (we compute both complexes by the Cousin resolution).

(ii) Let $i' : Z' \hra X' $ be another embedding as in 6.3.1,  and $f : X \to X'$ be a morphism transversal to $Z'$ such that $ Z= f^{-1}(Z')$. Then the next diagram commutes: 
$$\spreadmatrixlines{2\jot}
\matrix
M(X,X\smallsetminus Z )&
\buildrel{f}\over\lra & M(X',X'\smallsetminus Z')
 \\
i^{Gys} \downarrow&&i^{\prime Gys} \downarrow \\
M(Z)(n)[2n]  &
  \buildrel{f|_{Z'}}\over\lra & M(Z')(n)[2n].
\endmatrix
\tag 6.3.9$$

E.g., looking at $X' = X\times X$, $Z' = X\times Z$,  $f=\Delta$, we see that $i^{Gys}$ is a morphism of $M(X)$-comodules, hence it commutes with the $\cap$-product action of $H^{\cdot ,*}_\CM (X,R)$.

(iii)   Gysin's identification (6.3.2) for $F= R^\Delta_{\tr} (*)$ is an isomorphism  $$R\Gamma_\CM (Z, R(*-n))[-2n] \iso R\Gamma_{\CM\, Z}(X, R(*)) \tag 6.3.10$$ 
of $R\Gamma_\CM (X, R(*))$-modules.  Let $$\Res_Z : R\Gamma_\CM (X\smallsetminus Z, R(*))\to R\Gamma_\CM (Z, R(*-n))[-2n+1] \tag 6.3.11$$ be the composition of the inverse to (6.3.10) and the boundary map for $R\Gamma_\CM$; this is  a morphism of $R\Gamma_\CM (X, R(*))$-modules as well. 

(iv) Suppose that $Z$ is a divisor in $X$, so $n=1$, and $R=\Bbb Z$. Then (6.2.2) identifies the residue map $\Res_Z : H^{1,1}_\CM (X\smallsetminus Z, \Bbb Z (1))\to H^{0,0}_\CM (Z, \Bbb  Z)$ with the map $div_Z : \CO^\times (X\smallsetminus Z )\to \Bbb Z^{\pi_0 (Z)}$. 

{\bf 6.4. The class of a cycle.} Let $X$ be a smooth variety. Let $Z\subset X$ be a closed subvariety of codimension $\ge n$; denote by $[Z]$ the set of $n$-dimensional irreducible components of $Z$.  

\proclaim{\quad Proposition}  (i)  There is a canonical quasi-isomorphism $$cl : R^{[Z]}\iso  R\Gamma_\CM (X,X\smallsetminus Z, R (n) )[2n] .\tag 6.4.1$$

(ii) One has $H^{>2n,n}_\CM (X,R (n) )=0$. There is a canonical isomorphism $$cl : CH^n (X)\otimes R \iso H^{2n,n}_\CM (X,R). \tag 6.4.2$$
\endproclaim

{\it Proof.} Consider the coniveau filtration on $T:=R\Gamma_\CM (X, R (n) ) $. The corollaries in 6.3 and 6.1 imply that $\gr^{>n} T=0$ and $\gr^n T [2n] = \CZ^n (X)\otimes R$ (the group of codimension $n$ cycles on $X$).  
Since $\gr^{<n} T$ is the same for $X$ and $X\smallsetminus Z$, we get (i).  Since $H^{>m}\CC^\Delta (R (m)) =0$, one has $H^{\ge 2n}\gr^{<n} T=  H^{ 2n-1}\gr^{<n-1} T=0$. Therefore $H^{2n,n}_\CM (X,R )$ equals the cokernel of the
 differential $H^{2n-1} \gr^{n-1}T \to H^{2n}\gr^n T =\CZ^n (X)\otimes R$.
 
 Now $H^{2n-1} \gr^{n-1}T = \mathop\oplus\limits_{x\in X^{(n-1)}} \CO^\times (\eta_x )\otimes R$ by (6.2.2).  By Exercise (iv) in 6.3.4, the above differential  sends a rational function $\phi \in \CO^\times (\eta_x )$ to its divisor $div (\phi )\in \CZ^n (X)$. Hence the cokernel equals $CH^n (X)\otimes R$, q.e.d.
\hfill$\square$

\medskip

{\it Exercise.} Suppose that $ Z$ is smooth irreducible  of codimension $n$. Then 
$cl (Z)\in H^{2n,n}_\CM (X, X\smallsetminus Z,R)$ equals the composition $M(X,X\smallsetminus Z) \iso M(Z)(n)[2n] \to R (n)[2n],$ where the first arrow is the Gysin map, the second one is the augmentation  for $M(Z)$. Thus the composition $M(X,X\smallsetminus Z) \iso M(Z)(n)[2n] \to M(X) (n)[2n]$ is the morphism of $M(X)$-comodules that corresponds to $cl(Z)$ (see 6.2.2).
The $\cap$-product action $ M(X)\to M(X)(n)[2n]$ of $cl (Z)$ (see (6.2.3)) equals the composition $M(X)\to M(X, X\smallsetminus Z) \iso M(Z)(n)[2n] \to M(X)(n)[2n],$  the middle arrow is the Gysin map, the other two are  evident morphisms.

{\it Remark.} According to \cite{MVW} 19.1, for every $m,n$ there is a natural isomorphism $H^{m,n}_\CM (X,\Bbb Z )\iso CH^n (X,2n-m)$ where the r.h.s.~is  Bloch's higher Chow group.

\medskip

{\bf 6.5. The motive of a projective bundle.} For a line bundle $\CL$ on $X$ denote by $c_1^\CM (\CL )$  its class in $Pic (X)\otimes R\cong H^{2,1}_\CM (X,R )$, see (6.2.2).

{\it Exercise.} If $\lambda$ is a non-zero rational section of $\CL$, then $c_1^\CM (\CL )= cl (div (\lambda ))$.

Let now $\CE$ be any vector bundle of rank $n+1$ on $X$, and $p: \Bbb P (\CE )\to X$ be the corresponding projective bundle. Set $c:= c_1^\CM (\CO (1)) \in H^{2,1}_\CM (\Bbb P (\CE ),R )$. As in (6.2.3), let us  consider $c^a \in H^{2a,a}_\CM (\Bbb P (\CE ),R )$ as morphisms  $M (\Bbb P (\CE ))\to M (\Bbb P (\CE ))(a)[2a]$. Composing with $p$, we get the morphisms $p_* c^a : M (\Bbb P (\CE ))\to M (X)(a)[2a]$.

\proclaim{\quad Proposition} One has a natural homotopy equivalence $$(p_* c^a ): M (\Bbb P (\CE ))\iso \mathop\oplus\limits_{0\le a \le n} M (X)(a)[2a]. \tag 6.5.1$$\endproclaim

{\it Proof.} Choose a covering $\{ U_\alpha \}$ of $X$ such that $\CE $ is trivial on $U_\alpha $'s. Let $U_\cdot$ be the corresponding \v Cech hypercovering, so $M(U_\cdot )\iso M(X)$, $M(\Bbb P (\CE )_{U_\cdot})\iso M(\Bbb P (\CE ))$. Consider the corresponding filtrations on $M(X)$, $M(\Bbb P (\CE ))$. The operators $c^i$  and $p_*$ act naturally on the filtered objects. To show that (6.5.1) is a homotopy equivalence, it suffices to check that such is the associated graded map, which is $(p_* c^a ): M(\Bbb P (\CE )_{U_m}) \to \oplus M(U_m )(a)[2a]$. Thus we can assume that $\CE$ is a trivial vector bundle, so $\Bbb P (\CE )= \Bbb P^{n}_X$.

Consider the embeddings $ \Bbb P^{n-1}_X\hra \Bbb P^n_X$ and  $ \Bbb A_X^n := \Bbb P^n_X \smallsetminus \Bbb P^{n-1}_X \buildrel{j}\over\hra \Bbb P^n_X$. Since $c= cl (\Bbb P^{n-1}_X )$,  the exercise in 6.5 implies that $p_* c^{a}$ for $a\ge 1$ equals the composition $M(\Bbb P^n_X )\buildrel\alpha\over\to M (\Bbb P^n_X , \Bbb A^n_X )\iso M (\Bbb P_X^{n-1})(1)[2] \buildrel{c^{a-1}}\over\lra M (\Bbb P_X^{n-1})(a)[2a] \to M(X)(a)[2a]$ where the right arrow comes from the projection $\Bbb P^{n-1}_X \to X$.  Since the composition  $M(\Bbb A^n_X )\buildrel{j_*}\over \to M(\Bbb P^n_X )\buildrel{p_*}\over\to M(X)$ is a homotopy equivalence, we see that the morphism $M(\Bbb P^n_X )\buildrel{(\alpha ,p_* )}\over\lra M(\Bbb P^n_X , \Bbb A^n )\oplus M(X)$ is also a homotopy equivalence. Together with the previous assertion, this proves (6.5.1) by induction by $n$.  \hfill$\square$

\medskip

We see that there is a canonical quasi-isomorphism $$ \mathop\oplus\limits_{0\le a\le n} R\Gamma_\CM (X, R (m-a))[-2a] \iso R\Gamma_\CM (\Bbb P (\CE ), R (m)) \tag 6.5.2$$
which sends $(\phi_a )$ to $ \sum \phi_a c^a $. So one can define the Chern classes $c_i^\CM (\CE )\in H^{2i,i}_\CM (X,R )$ by the usual Grothendieck procedure.

\medskip

{\bf 6.6. The motive of a blow-up.} Let $i: Z\hra X$ be a closed embedding of smooth varieties of codimension $m$ and  $f: \tilde X\to X$ be the  blow-up. Set $\tilde{Z}:= f^{-1}(Z)$, $U: = \tilde X \smallsetminus \tilde Z= X\smallsetminus Z$; then $\tilde{Z} = \Bbb P (N)$ where $N$  is the normal bundle to $Z$ in $X$.   Let   $\alpha$  be the composition  
$ M(\tilde{ X}) \to M(\tilde{X},U)\iso M(  \Bbb P (N))(1)[2]\iso \mathop\oplus\limits_{0 < a \le m} M(Z)(a)[2a] \to   \mathop\oplus\limits_{0 < a < m} M(Z)(a)[2a], $
 where the second morphism is the Gysin equivalence (6.3.1), the third one  is (6.5.1), the last arrow is the projection.
\proclaim{\quad Proposition}  The map $(f,\alpha )$ is a quasi-isomorphism:
 $$ M(\tilde {X}) \iso M(X)\oplus \mathop\oplus\limits_{0 < a < m} M(Z)(a)[2a]. \tag 6.6.1$$
\endproclaim
{\it Proof.}  $(f,\alpha )$ is compatible with the evident morphisms from $M(U)$ to the l.h.s.~and r.h.s.~of (6.6.1). So it suffices to show that the corresponding morphism of cones $M(\tilde {X},U) \iso M(X,U)\oplus \mathop\oplus\limits_{0 < a < m} M(Z)(a)[2a]$ is a quasi-isomorphism.

 By the construction of $\alpha$, the assertion follows if we check that the composition $\xi_1$ of  $
 M(\tilde{X},U)\iso M(  \Bbb P (N))(1)[2]\iso \mathop\oplus\limits_{0 < a \le m} M(Z)(a)[2a] \to    M(Z)(m)[2m] $ equals the composition  $\xi_2$ of $M(\tilde{X},U)\buildrel{f}\over\to M(X,U) \buildrel{i^{Gys}}\over\lra M(Z)(m)[2m]$. 
 
Looking at the first step of the construction of the Gysin equivalence (the first isomorphism in (6.3.3)) described in 6.3.2, we see that it is compatible with the above constructions. So we can assume that $X$ is a vector bundle $N$ over $Z$, $i$ is its zero section, $\tilde{X}$ is the space $L$ of the line bundle $\CO (-1)$ over $\Bbb P (N)$, $f$ is the corresponding standard map.
 One has a commutative diagram 
  $$\spreadmatrixlines{1\jot}
\matrix
\tilde{X}=L&
\buildrel{f}\over\lra &
 X=N \\  
\uparrow\downarrow&&\uparrow\downarrow \\
\Bbb P (N) &
    \buildrel{p}\over\lra & Z.
\endmatrix
\tag 6.6.2$$ 

We are playing with $Z$-schemes, so our motives are $M(Z)$-comodules (see 6.2.2). By the exercise in 6.4 and 6.5, $\xi_1$ is the morphism of $M(Z)$-comodules that corresponds to $cl (\tilde{Z})\cup c_1 (\CO (1))^{\cup n-1} \in H^{2n,n}_\CM (\tilde{X},U;R)$, and $\xi_2$ to $f^* cl (Z)$. It remains to check that these classes are equal. As in 6.3.3, one has $H^{2n,n}_\CM (\tilde{X},U;R)=R^{[Z]}$, so we can assume that $Z=\Spec\, k$. The rest is left to the reader.  \hfill$\square$

{\it Remark.} The  proposition implies that  after the motivic localization  the complex
$$    R_{\tr}[ \tilde {Z}] \to R_{\tr}[\tilde{ X}]\oplus R_{\tr}[Z]\to R_{\tr}[X] \tag 6.6.3$$
becomes quasi-isomorphic to $0$. We will see in 6.9.3 that,  if  resolution of singularities is available,  this statement is valid  without smoothness assumption on $X$ and $Z$.

\medskip

{\bf 6.7. Poincar\'e duality.} 6.7.1. We live in a tensor category with unit object $\b 1$. As in \cite{D2},  a {\it duality datum} consists of a pair of objects $V$, $V^*$ and morphisms $\epsilon : V^* \otimes V \to \b 1$, $\delta : \b 1 \to V\otimes V^*$ such that both compositions $V\buildrel{\delta \otimes \id_V}\over\lra V\otimes V^* \otimes V \buildrel{\id_V \otimes\epsilon}\over\lra V$ and $V^* \buildrel{\id_{V^*}\otimes \delta }\over\lra V^* \otimes V \otimes V^* \buildrel{\epsilon\otimes \id_{V^*}}\over\lra V^*$ are the identity morphisms. For given $V$ the triple $(V^*,\epsilon ,\delta )$ is uniquely defined (if it exists), and called the {\it dual} object  to $V$.  Since $\delta$ is determined uniquely by $\epsilon$ (and vise versa),  it suffices to specify either of them.

Now our tensor category is $ \CD_\CM^{\text{tri}}$. For a smooth variety $X$ of dimension $n$ let 
$$  \epsilon'_X   :\,  M(X)\otimes M(X)\to R_\CM (n)[2n]  \tag 6.7.1$$  
 be the composition $M(X)\otimes M(X)\iso M(X\times X)\buildrel{cl(\Delta )}\over\lra R_\CM (n)[2n]$ where $cl (\Delta )\in H^{2n,n}_\CM (X\times X,R)$ is the class of the diagonal cycle $\Delta : X\hra X\times X$. Equivalently, this is the composition
  $M(X\times X)\to M(X\times X ,   (X\times X)\smallsetminus \Delta (X)  )\iso M(X)(n)[2n]\to R_\CM (n)[2n]$; here $\iso$ is the Gysin map, the last arrow is the augmentation.  Set $\epsilon_X := \epsilon'_X(-n)[-2n]  :\,  M(X)(-n)[-2n]\otimes M(X)\to R_\CM .$

\proclaim{\quad Proposition} If $X$ is proper, then  $(M(X)(-n)[-2n],\epsilon_X )$  is the dual to $M(X)$.
\endproclaim

{\it Proof.} We can assume that $X$ is irreducible and, by an argument similar to one in Remark (b) in 4.6.1,  that $k$ is infinite. 

Let $\delta_X$ be the composition 
 $R_\CM (n)[2n]  \buildrel{\kappa}\over\to M(\Bbb P^n ) \buildrel{\nu}\over\to M(X) \buildrel{\Delta}\over\to M(X\times X)$,  where $\kappa$ is a component of (6.5.1),  and $\nu \in \CC or (  \Bbb P^n , X)$ is any correspondence whose graph is an irreducible cycle $\Gamma_\nu \subset  X\times \Bbb P^n  $ such that at the generic point of $\Gamma_\nu$  the projection $\Gamma_\nu \to X$ is an isomorphism and the one $\Gamma_\nu \to \Bbb P^n$ is \'etale  ($\nu$ exists by  Chow's lemma \cite{EGA II} 5.6.1, and 5.5). Set  $\delta_X :=\delta'_X (-n)[-2n]:\,  R_\CM \to M(X)\otimes M(X)(-n)[-2n]$.

Let us check that $(\epsilon_X ,\delta_X )$ is a duality datum. Our $\epsilon'_X$ and $\delta'_X$ are symmetric, so we need to show that the composition $$M(X)(n)[2n] \buildrel{\beta }\over\longrightarrow M(X\times X\times X ) \buildrel{\alpha }\over\lra M(X)(n)[2n] , \tag 6.7.2 $$ where $\alpha := \id_{M(X)}\otimes \epsilon'_X$ and
$\beta := \delta'_X \otimes \id_{M(X)}$, is the identity morphism. 

(a) Let $p_a : X\times X \to X$, $a=1,2$, and $pr_i ,pr_{ij}: X\times X \times X \to X, X\times X$, $i,j=1,2,3$ be the projections. Then $\alpha$ equals the $\cap$-product action of $ pr_{23}^* (cl (\Delta ))$ followed by $pr_1$. The pull-back of $ pr_{23}^* (cl (\Delta ))$ by $\Delta\times \id_X : X\times X\to X\times X\times X$ equals $cl (\Delta )$, so 
 the composition $\gamma$ of $$M(X\times X)\,\, \buildrel{\Delta\times \id_X}\over\lra \,\, M(X\times X\times X)\buildrel{\alpha }\over\lra M(X)(n)[2n]$$ equals $\tilde{cl}(\Delta )$:= the $\cap$-product action of $ cl (\Delta )$, followed by $p_1$.

Notice that $\tilde{cl}(\Delta )$
equals the composition $$M(X\times X) \to M(X\times X , (X\times X)\smallsetminus \Delta (X)) \iso M(X)(n)[2n] \buildrel\Delta\over\lra M(X\times X)(n)[2n]$$  (see Exercise in 6.4). The transposition $\sigma$ of factors of $X\times X$ acts on our objects; since it acts trivially on the term $M(X)(n)[2n]$, we see that $\tilde{cl} (\Delta )$ is invariant with respect to the left and right compositions with the transposition. Therefore $\gamma =  p_1 \sigma \tilde{cl}(\Delta )= p_2 \tilde{cl}(\Delta )$.

Now  $p_2 $ defines on  $ M(X\times X)$ an $ M(X)$-comodule structure (see 6.2.2). Our $\gamma$
  is a  morphism of $M(X)$-comodules which corresponds to $cl (\Delta )$ via the identification $\Hom^{M(X)} (M(X\times X), M(X)(n)[2n])\iso  \Hom (M(X\times X), R_\CM (n)[2n])$ $= H^{2n,n}_\CM (X\times X ,R )$ (see 6.2.2(b); here $\iso$ is the composition with the augmentation).

(b) Set $\nu_X := \nu \otimes \id_{M(X)}: M(\Bbb P^n \times X)\to M(X\times X)$; this is a morphism of $M(X)$-comodules. By (a),  the composition $ M(\Bbb P^n \times X)\buildrel{\nu_X}\over\to  M(X\times X) \buildrel{\gamma}\over\to M(X)(n)[2n]$ is a  morphism of $M(X)$-comodules that  corresponds to $\nu_X^* (cl (\Delta ))$  via the identification $\Hom^{M(X)} ( M(\Bbb P^n \times X), M(X)(n)[2n])\iso \Hom ( M(\Bbb P^n \times X), R_\CM (n)[2n]) = H^{2n,n}_\CM ( \Bbb P^n \times X, R )$. By (6.3.9), the conditions on $\nu$ imply that $\nu_X^* (cl (\Delta ))= cl (\Gamma_\nu ).$ 

(c) Set $j_X := j\otimes \id_{M(X)} : M(X)(n)[2n] \to M(\Bbb P^n \times X)$; this  is a morphism of $M(X)$-comodules. One has $\alpha\beta = \gamma \nu_X j_X :  M(X)(n)[2n]\to M(X)(n)[2n]$. By (b), this is
the $M(X)$-comodule endomorphism that corresponds to the class $j_X^* (cl (\Gamma_\nu )) \in \Hom (M(X)(n)[2n] ,R_\CM (n)[2n])= H^0_\CM (X,R (0))=R$. 

Consider the  K\"unneth decomposition for the Chow groups $H^{2n,n}_\CM (\Bbb P^n \times X ,R )= \sum H^{2(n-a), n-a}_\CM (\Bbb P^n , R )\otimes H^{2a, a}_\CM (X,R )= \sum H^{2a,a}_\CM (X,R )$ (see (6.5.2)). It shows that for a cycle $Z$ the number $j_X^* (cl (Z))$ is simply the degree of $Z$ over the generic point of $X$. Thus $j_X^* (cl (\Gamma_\nu )) =1$, so $\alpha\beta$ is the identity morphism, q.e.d. \hfill$\square$

6.7.2. If a tensor category is Karoubian, then its objects that admit a dual form a Karoubian subcategory.
So for a Karoubian tensor triangulated category its objects that
 admit a dual form a thick subcategory. If resolution of singularities is available, then $ \CD_\CM^{\text{tri}}$ is the idempotent completion of the triangulated subcategory strongly generated by
motives of proper smooth varieties, hence 

\proclaim{ \quad  Corollary} If $k$ admits resolution of singularities, then $ \CD_\CM^{\text{tri}}$ is a rigid tensor category, i.e., each of its objects admits a dual.  \hfill$\square$    \endproclaim

See 6.9.4 for an explicit construction of the dual to $M(X)$ for non-compact $X$.

6.7.3. If $P,Q$ are objects of a tensor category such that $P$ admits dual $P^*$, then $\Hom (Q,P)=\Hom (Q\otimes P^* ,\b 1 )$. This is  formula (6.7.3) below; the rest of the corollary is an exercise for the reader:

 \proclaim{\quad  Corollary} For  $X$, $Y$  smooth, $X$  proper, one has a natural identification $$\Hom (M(Y), M(X))\iso R\Gamma_\CM (X\times Y, R(\dim X ))[2\dim X].  \tag 6.7.3$$  Thus  $H^0 \Hom (M(Y), M(X))= CH^{\dim X} (X\times Y)\otimes R$ and $H^{>0} \Hom (M(Y), M(X))$ $=0$. 
For   $c\in \CC or (X,Y)$ the class in $CH^{\dim X} (X\times Y)$ of  $M(Y)\buildrel{c}\over\to M(X)$ equals the class of $c$ as of a cycle on $X\times Y$. The composition of morphisms of motives of proper varieties equals the composition of Chow correspondences. \hfill$\square$
\endproclaim

6.7.4. Consider an additive $R$-category of smooth projective varieties and  Chow correspondences: $\Hom (Y,X)= CH^{\dim X}(X\times Y )\otimes R$. Its idempotent completion is the category $\CC \CH^{\text{eff}}_\CM $ of {\it effective Chow motives}. Let $ \CC \CH_\CM^{\text{eff\,pretr}}$ be the DG tensor category of finite complexes in $\CC \CH^{\text{eff}}_\CM$, and $\CC \CH_\CM^{\text{eff\, tri}}$ be its homotopy category.
The next result is due to M.~Bondarko \cite{Bo1} \S 6:

\proclaim{\quad Proposition} For $k$ that admits resolution of singularities, there is a natural homotopy tensor  DG functor $$\epsilon : \CD_\CM^{\text{eff}}\to \CC \CH_\CM^{\text{eff\, pretr}}.\tag 6.7.4$$ The  functor $\epsilon^{\text{tri}}:   \CD_\CM^{\text{eff\, tri}}\to \CC \CH_\CM^{\text{eff\, tri}} $ is conservative and yields an isomorphism of the $K_0$-groups.
\endproclaim

{\it Proof.} Let $\CB$ be the full DG subcategory in $\CD_\CM^{\text{eff}}$ formed by those motives that are isomorphic in $\CD_\CM^{\text{eff\, tri}}$   to direct summands of motives of smooth proper varieties. By  
6.7.3, $B:= \Ho\, \CB $ equals $   \CC \CH^{\text{eff}}_\CM$, and for every $M,N\in \CB$ and $n>0$ one has $H^{>0}\Hom (M,N)=0$. By (iii) of the proposition in 1.5.6, the triangulated subcategory of  $\CD_\CM^{\text{eff\, tri}}$ strongly
generated by $B$ is Karoubian. Therefore, 
by the assumption on $k$, it equals   $\CD_\CM^{\text{eff\, tri}}$. Now we can apply 1.5.6 to $\CA = 
 \CD_\CM^{\text{eff}}$ and $B$. The compatibility of $\epsilon$ with the tensor products is evident.
 \hfill$\square$

\medskip

 {\bf 6.8. The Steinberg relation.} Below $R=\Bbb Z$. Let $X$ be a smooth variety, $f,g \in\CO^\times (X)= H^{1,1}_\CM (X,\Bbb Z )$ be invertible functions; we get $f\cup g \in H^{2,2}_\CM (X, \Bbb Z )$.

 \proclaim{\quad Proposition}  If $f+g =1$, then $ f\cup g  =0$.
  \endproclaim
 
{\it Proof.} An observation: Suppose  $X$ is the complement to a smooth divisor $D$ in a smooth variety $\bar{X}$. We have the map $\Res_D : H^{\cdot ,\cdot }_\CM (X,\Bbb Z)\to R\Gamma_\CM^{\cdot -1 ,\cdot -1}(D,\Bbb Z )$ (see (6.3.11)) which is a morphism of $H^{\cdot ,\cdot }_\CM (\bar{X},\Bbb Z)$-modules. By Exercise (iv) in 6.3.4, for $g\in H^{1,1}_\CM (X,\Bbb Z)=\CO^\times (X)$ one has $\Res_D (g)= div_D (g)$. Therefore if $f\in \CO^\times (\bar{X} )$, then $\Res_D (f\cup g)\in 
H^{1,1}_\CM (D,\Bbb Z )=\CO^\times (D)$ equals
$(f|_D )^{div_D (g)}$ (see Exercise below for a general statement). In particular, if $f$ and $g$ are such that on each irreducible component of $D$ either $f$ or $g$ is regular and equals  1, then $\Res_D (f\cup g)$ vanishes.

To  prove the proposition, 
 it suffices to consider the universal situation of $X= \Bbb A^1 \smallsetminus D$, where $D=\{ 0,1\}$, $f=t$, $g=1-t$. Set $\alpha := f\cup g$. By the  observation, $\Res_D (\alpha )=0$, so 
the  exact sequence $0 \to H^{2,2}_\CM (\Bbb A^1 ,\Bbb Z)\to H^{2,2}_\CM (X, \Bbb Z)\buildrel{\Res_D}\over\lra k^\times \times k^\times \to 0$  implies that $\alpha \in  H^{2,2}_\CM (\Bbb A^1 ,\Bbb Z) =H^{2,2}_\CM (k,\Bbb Z )$. Since $\cup$ is skew-commutative, the symmetry $t\mapsto 1-t$ of $X$ sends $\alpha$ to $-\alpha$. Since symmetries of $X$ act on the subgroup $H^{2,2}_\CM (k,\Bbb Z )$ of $H^{2,2}_\CM (X,\Bbb Z )$ as identity, one has $2\alpha =0$. Thus $t^2 \cup (1-t)=0=t^2 \cup (1+t )$, so the element $t^2 \cup (1-t^2 )= t^2 \cup (1-t) + t^2 \cup (1+t) $ of $ H_\CM^{2,2}(\Bbb A^1 \smallsetminus \{ 0,\pm 1\} ,\Bbb Z )$ vanishes. But this
 is the pull-back of $\alpha$ by the map $t\mapsto t^2$. The latter map acts as identity on the subgroup $H^{2,2}_\CM (k,\Bbb Z )$, so $\alpha =0$, q.e.d. \hfill$\square$

{\it Exercise.} In the situation of the above observation, for any $f,g\in\CO^\times (X)$ the function $\Res_D (f\cup g)\in\CO^\times (X)$ coincides with the tame symbol $\{ f,g\}_D$.

Recall that if $X$ is a local scheme, then the {\it Milnor ring} $K^{M}_\cdot (X)$ is a graded associative ring generated by $K^{M}_1 (X):= \CO^\times (X)$ modulo the Steinberg relation $\{ f,g \} =0$ if $f, g\in \CO^\times (X)$, $f+g =1$; here $\{ \, \}$ denotes the product in $K^{M}_\cdot (X)$. By the proposition, the identification $\CO^\times (X)\iso H_\CM^{1,1}(X,\Bbb Z )$ extends to a morphism of graded algebras $K^{M}_a (X)\to H^{a,a}_\CM (X,\Bbb Z )$. According to \cite{Ke} (and \cite{NS}, lecture 5 of \cite{MVW}, where the case of a field is treated), this is an isomorphism.

\medskip

{\bf 6.9. The cdh localization.} 
Below ``scheme" means ``separated $k$-scheme of finite type".

6.9.1. The {\it cdh topology} is  the weakest Grothendieck topology on the category of schemes whose covers  include every proper or \'etale map which  admits a constructible section (see 4.2.2). Thus every cdh cover admits a constructible section. The cdh topology 
 is stronger than the  Nisnevich one. 

A proper  map $p: Y\to X$ is called an {\it abstract blow-up} if there is an open dense 
 $U\subset X$ such that  $p^{-1}(U)_{\text{red}}\to U_{\text{red}}$ is an isomorphism.

{\it Example.} Let $ Y\to X$ be an abstract blow-up and $Z\hra X$ be a closed subscheme such that $U:=X\smallsetminus Z$ satisfies the above condition. Then $Y\sqcup Z \to X$ is a cdh cover. Such cdh covers are called {\it special}.

The next lemma\footnote{This is a special case of \cite{RG} Th.~5.2.2, which says that for any  $ T/ X$ one can find an abstract blow-up $Y/ X$ such that the strict transform $\bar{T}_{\eta_X} $ is flat over $Y$.} an essential tool:
\proclaim{\quad Lemma} Let $T/X$ be an $X$-scheme proper and generically finite  over $X$. Then there exists an abstract blow-up $Y/ X$ such that the strict transform 
$\bar{ T}_{\eta_X} \subset T\times _X Y$ (which is the closure in  $  T\times _X Y$ of the generic fibers
$T_{\eta_X}$) is finite over $Y$. 
\endproclaim

{\it Proof of Lemma.} Since the projection to $X$ of the disjoint union of normalizations of its reduced irreducible components is an abstract blow-up, we can assume that $X$ is integral and normal.  By Chow's lemma, we can assume that $T/X$ is projective.
Let $U\subset X$ be an open dense subset such that $ T_U / U$ is finite flat of some degree $d$.   We get a section $ U\to \Sym^d (T_U /U)\subset \Sym^d (T/X)$ (see   (a) of the proof in 2.1.3).  Our $Y$ is the  closure of its image; the promised property holds by (b) in loc.~cit.  \hfill$\square$

Below ``proper cdh cover" means a proper map which is a cdh cover (i.e., admits a constructible section).

\proclaim{\quad Proposition} (i) Every proper cdh cover $  T/X$ can be refined to   a composition of special cdh covers $F/X$. Thus
the cdh topology is generated by the Nisnevich covers and special cdh covers. 

(ii) Every cdh cover $V/X$ admits a refinement $W/X$ which can be factored as $W\to Q\to X$ where $W/Q$ is a Nisnevich cover, $Q/X$ is a proper cdh cover.
\endproclaim

{\it Proof.} (i) We use induction by $\dim X$.  Let $Z\subset X$ be a nowhere dense closed subscheme such that 
 $T/X$ admits a section $s$ over $(X\smallsetminus Z)_{\text{red}}$. Let $Y$ be the closure of $s((X\smallsetminus Z)_{\text{red}})$ in $T$. Then $Y\sqcup Z /X$ is a special cdh cover. By the induction assumption, we find a composition of special cdh covers $F_Z /Z$ which refines the cover
 $T_Z := T\times_X Z $ of $ Z$. Then  $F:= Y\sqcup F_Z \to X$ is a composition of special cdh covers which factors through $Y\sqcup T_Z /X$, hence through $T/X$, q.e.d.  
 
 (ii)  It suffices to consider the case when $V/X$ is a composition $V\to U\to X$ where $U/X$ is a Nisnevich cover, $V/U$ is a proper cdh cover. By (i) we can assume that $V / U$ is an isomorphism at the generic points of $U$; by Chow's lemma, we can assume that $V/U$ is projective.
 
 It suffices to find an abstract blow-up $Y/X$ such that the cover  $V\times_X Y /Y$   admits a Nisnevich refinement. Indeed, let $Z\subset X$ be a nowhere dense closed subscheme such that 
 $Y\sqcup Z /X$ is a special cdh cover. The cdh cover $V|_Z$ of $Z$ admits a refinement as in (ii) by induction by $\dim X$, so  $V\times_X (Y\sqcup Z )/Y\sqcup Z$ admits  a refinement $W/Y\sqcup Z$ as in (ii). Then $W/X$ is the promised refinement of $V/X$.
  
We will find an abstract blow-up $Y/X$ such that each connected component of $Y$ is normal and maps onto an irreducible component of $X$,  and  the strict transform of $V $ in $V\times _X Y$  is finite  over $U\times_X Y$. Then   $V\times _X Y /  U\times _X Y$ admits a (unique) section,\footnote{For $U\times _X Y$ is normal  and $V\times _X Y \to U\times _X Y$ is an isomorphism over its generic points.  } so $V\times_X Y/Y$  admits a Nisnevich refinement  (namely, $U\times _X Y /Y$) as needed.
 
  To construct  $Y/X$, we compactify $V \to U \to X$, i.e., 
find  projective maps $ \bar{V} \to  \bar{ U} \to X$ with $U$  an open dense subset of $\bar{U}$,  $\bar{V} \times_{\bar{ U}} U =V$,  and then apply the lemma to $\bar{V}/X$. \hfill$\square$
 
\medskip

6.9.2. Most  applications of the cdh topology require the resolution of singularities, so from now on we assume that char $ k=0$. First implications:

(i) {\it  For every separated scheme $X$ of finite type there exists a proper cdh cover $ \tilde{X}/X$  with smooth $\tilde{X}$.}
Indeed, there is a special cdh cover $Y\sqcup Z \to X$ as in Example in 6.9.1 with smooth $Y$, and, by induction by $\dim X$, one can find a proper cdh cover $\tilde{Z}\to Z$ with $\tilde{Z}$ smooth; set $\tilde{X}:= Y\sqcup\tilde{Z}$.

 (ii) {\it Every abstract blow-up $ Y\to X$ of a  smooth scheme $X$ is a cdh cover.\footnote{The claim  might be false if $X$ is singular: Let $Z\subset \Bbb P^n$ be any  closed subvariety with $Z(k)=\emptyset$; let $X\subset \Bbb A^{n+1}$ be the cone over $Z$, $Y$ be the blow-up of $X$ at $0$. Then $Y/X$ is {\it not} a cdh cover.  }} Indeed,  there exists a sequence of blow-ups $X_k\to \cdots 
\to X_0 =X$ with smooth centers such that the composition factors through $p$, and  all $X_i \to X_{i-1}$ are  cdh covers for obvious reasons.  

(iii) {\it Each proper cdh cover of a  smooth  $X$  admits a refinement which can be factored as  $ X_k\to  \cdots\to X_0=X$, where  $X_i / X_{i-1}$ are  blow-ups  with smooth centers.} This follows from (ii), since $Y/X$ admits an abstract blow-up refinement.

\medskip

Consider the cdh topology on $\CS m$. As in 1.11, we have the subcategories $I^{\text{cdh}}, I^{\text{cdh}\bot} \subset D\CP \CS h$ of, respectively, cdh locally acyclic and cdh local complexes. We have the following standard objects of $I^{\text{cdh}}$:

(a) Nisnevich Mayer-Vietoris complexes, see 4.2.1.

(b) Complexes $R [\tilde Z] \to R[\tilde{X}]\oplus R[Z]\to R[X]$ where   $\tilde{X}/X$ is the blow-up with a smooth center $Z\subset X$ and  $\tilde{Z}\subset\tilde{X}$ is the preimage of $Z$; the differentials are the difference and the sum of the evident maps.
 
\proclaim{\quad Proposition} $I^{\text{cdh}}$ is  generated by all complexes of types (a), (b).
\endproclaim

{\it Proof.} As in the proof of the proposition in 4.2.1, we may apply 1.4.2(i), with $S$ the set of standard objects in (a) and (b). Thus it suffices to show that if 
 $F\in I^{\text{cdh}}$ is such that  the  total complexes  $F(X)\to F(\tilde{X})\oplus F(Z)\to F(\tilde{Z})$ for  all data of (b) and similar complexes for all data of (a) are acyclic, then $F=0$. 
 
 We prove that $F(X)$ is acyclic by induction by $\dim X$. Assuming that $F(Y)$ is acyclic for every $Y$ of dimension $< m$, let us check that for every $X$ of dimension $m$ and $h\in H^a F(X)$ there is an open $V\subset X$ with $X\smallsetminus V$ of codimension $\ge n$ such that $h|_V \in H^a F(V)$ vanishes. Since  for $n=m+1$ this means that $F(X)$ is acyclic, the claim yields the proposition.

For any $\tilde{X}/X$ as in (b) one has $\dim Z,$ $ \dim\tilde{Z}<m$, hence $F(X)\iso F(\tilde{X})$. Thus, by  implication (iii) above,  every proper cdh cover of $X$ admits a refinement $T/X$, $\dim T=m$, such that $F(U)\iso F(T_U )$ for every open $U\subset X$.

Since $F\in I^{\text{cdh}}$, there is a cdh cover $W/X$ such that $h|_W \in H^a F (W)$ vanishes. By (ii) of the proposition in 6.9.1, we can, refining $W$ if needed,  factor it as $W\to T\buildrel{\pi}\over\to X$, where $T/X$ is as above, $W/T$ is a Nisnevich cover. By the argument from the proof of the proposition in 4.2.1 (with $X$ from loc.~cit.~equal to  $T$), there is an open $V_T \subset T$ with $P:=T\smallsetminus V_T$ of codimension $\ge n$ such that $h|_{V_T}=0$. The promised $V$ is $X\smallsetminus \pi (P)$ (here $h_V =0$ since $F(V)\iso F(T_V )$ and $T_V \subset V_T$). \hfill$\square$

\medskip

6.9.3.  Let us show that 
replacing the Nisnevich topology by  the cdh one  does  not change the category of motives. 
\proclaim{\quad Lemma}  The cdh coverings define a Grothendieck topology on $R_{\tr}[\CS m]$.
\endproclaim

{\it Proof.}  As in the proposition in 4.3, the lemma follows from the next assertion (and implication (i) in 6.9.2):  {\it Let $X$, $Y$ be schemes, $\pi : U\to X$ be a cdh cover; suppose $Y$ is normal. Then for any
 $  \gamma \in \CC or(Y,X)$  there exist a cdh cover $\pi' : V \to Y$ with normal $V$ and $\tilde{ \gamma} \in \CC or(V,U)$ such that $\pi \tilde{\gamma}= \gamma \pi'$. }

Since  $U/X$ can be refined to a cdh cover which is a composition of Nisnevich covers and proper cdh covers, it is enough to check the claim when $\pi$ is in either of the above classes of covers.
For the Nisnevich covers, see the lemma in 4.3.

 Suppose   $\pi$ is proper. By the argument from the proof of the lemma in 4.3, we can  assume that $\gamma$ is given by an irreducible $\Gamma\subset X\times Y$.  Let $\eta_\Gamma$ be its generic point. By the  cdh property,  $\eta_\Gamma \hra X\times Y$ lifts to  $ \eta_\Gamma \hra U\times Y$; let $\Gamma' \subset U\times Y$ be the closure of its image. Then $ \Gamma'$ is  proper and generically finite  over $Y$. By  the lemma in 6.9.1 and implication (i) in 6.9.2, there is an abstract blow-up $\pi' : V\to Y$ with smooth $V$, such that the strict transform $\tilde{\Gamma}$ of $\Gamma'$ is finite over $V$. It yields  $\tilde{ \gamma} \in\CC or(V,U)$  such that $\pi \tilde{\gamma}= \gamma \pi'$, q.e.d.   \hfill$\square$

 Following 4.4 and 2.3, we define $\CI_{\tr}^{\text{cdh}} \subset \CP_{\tr}$ as the homotopy idempotent completion of the full pretriangulated subcategory strongly generated by the Nisnevich Mayer-Vietoris complexes and the complexes  $R_{\tr} [\tilde Z] \to R_{\tr}[\tilde{X}]\oplus R_{\tr}[Z]\to R_{\tr}[X]$ for data of type (b) in  6.9.2.  The categories $\CI^{\Delta\text{cdh}}_{\tr}$, $I^{\text{cdh}}_{\tr}$, etc., are defined then as in 4.4 replacing index ``Nis" by ``cdh".

\proclaim{\quad Proposition} (i) One has $\CI^{\Delta\text{cdh}}_{\tr}=\CI^{\Delta\text{Nis}}_{\tr}$, so $\CD_\CM^{\text{eff}}\iso \CP_{\tr}/^\kappa\, \CI^{\Delta\text{cdh}}_{\tr}$, etc.

(ii) The subcategories  $  I^\Delta_{\tr}, I^{\text{cdh}}_{\tr}\subset D\CP\CS h_{\tr}$ are compatible (see 1.3). 

(iii) For any $\Bbb A^1$-homotopy invariant complex $F$ of presheaves with transfers one has $H^\cdot_{\text{Nis}}(X,F)\iso H^\cdot_{\text{cdh}}(X,F)$.
\endproclaim

{\it Proof.} (i) By 6.6 and the proposition in 6.9.2,  the generators of $\CI^{\text{cdh}}_{\tr}$ lie in $\CI^{\Delta\text{Nis}}_{\tr}$, and we are done. 

(iii) One has $\CC^{\text{Nis}} (F)\in 
I^{\Delta\text{Nis}\bot}_{\tr}= I^{\Delta\text{cdh}\bot}_{\tr}\subset  I^{\text{cdh}\bot}_{\tr}$ (here ``$\in$"  is by  4.4(i), ``=" by (i) above) and $\CC one (F\to \CC^{\text{Nis}}(F))\in   I^{\text{Nis}}_{\tr} \subset I^{\text{cdh}}_{\tr}$, so  $\CC^{\text{Nis}} (F)=\CC^{\text{cdh}}(F)$.

(ii)  By (iii) and  the theorem in 4.5, for $F$ as in (iii) $H^\cdot_{\text{cdh}}(X,F)$ is $\Bbb A^1$-homotopy invariant. Now a cdh variant of the argument in 4.5(a) does the job. \hfill$\square$

\proclaim{\quad Corollary}Let $p: Y\to X$ be an abstract blow-up, $Z\subset X$ be a closed subscheme such that  $p^{-1}(X-Z)_{\text{red}}\iso (X-Z)_{\text{red}}$. 
Then the motivic localization of the complex
$   R_{tr}[ p^{-1} (Z)] \to R_{tr}[Y]\oplus R_{tr}[Z]\to R_{tr}[X]$
is quasi-isomorphic to $0$.  \hfill$\square$
\endproclaim

In turn,  the corollary implies that  $M(X)$ is a geometric motive (see 2.3) for any (not necessary smooth) scheme $X$.  

6.9.4. For a $k$-scheme $X$  we denote by $R^c_{tr}[X]$ the presheaf with transfers which assigns to a smooth connected scheme $Y$  the free 
$R$-module generated by
irreducible cycles $\Gamma\subset X\times Y$  such that the projection $p_1: \Gamma \to Y$ is quasi-finite and dominant over  $Y$.
Notice that $X\mapsto    R^c_{tr}[X]$ is covariant functorial with respect to proper morphisms  and contravariant functorial with respect to  flat quasi-finite morphisms.
Let $M^c(X)$ denote the motivic localization of $ R^c_{tr}[X]$.
\proclaim{\quad Lemma} For any  open  $U\subset X$ the motivic localization of the complex 
$$    R^c_{tr}[X-U] \to   R^c_{tr}[X] \to   R^c_{tr}[U] \tag 6.9.1$$
is quasi-isomorphic to $0$.
\endproclaim 

{\it Proof.} By the proposition in 6.9.3, it suffices to show that the cdh localization of (6.9.1) is a short exact sequence of cdh sheaves. We only need to check surjectivity of the right arrow. One can assume that $X$ is proper. Take any $\gamma\in   R^c_{tr}[U](Y)$. We look for  a cdh cover $p: V\to Y$ with  smooth $V$ such that $p^*(\gamma)\in   R^c_{tr}[U](V)$ 
comes from some $\overline  {p^*(\gamma)}\in   R^c_{tr}[X](V)$. One can assume that $\gamma$ is given by
 an irreducible cycle $\Gamma\subset U\times Y$.  Let $\bar{\Gamma}$ be its closure in 
$X\times Y$. Then $\bar{\Gamma}/ Y$ is proper and generically finite. By  Lemma in 6.9.1 we find an abstract  blow-up $p:V\to Y$ with  smooth $V$ such that
the  strict transform $\bar{\Gamma}_V$  of $\bar{\Gamma}$ is finite over $V$. Our
$\overline  {p^*(\gamma)}$ is  $\bar{\Gamma}_V$.  \hfill$\square$

For any smooth $X$   of dimension $n$ we define a canonical morphism in $\CD_\CM^{\text {eff\, tri}}$
$$\epsilon^{\prime c}_X : M^c (X)\otimes M(X)\to R_\CM (n)[2n] \tag 6.9.2
$$  as follows. 
Choose any smooth compactification 
$ X \subset\bar{X}$; set $D= \bar{X }\smallsetminus X$. By the lemma, $M^c (X)\iso M(\bar{X},D):= C^{{\CM}}((R_{tr}[\bar{ X}]/R_{tr}[D])$. Notice that $D\times X$ is disjoint from $\Delta (X)$ of $X$
 in $\bar{X}\times X$. Now (6.9.2) is the composition $M(\bar{X},D)\otimes M(X)= M(\bar{X} \times X ,D\times X)\to M(\bar{X}\times X, \bar{X}\times X \smallsetminus \Delta (X))$ 
$\buildrel{cl (\Delta )}\over\lra  R_\CM (n)[2n]$; here $cl (\Delta )\in H^{2n,n}_\CM (\bar{X}\times X, \bar{X}\times X \smallsetminus \Delta (X),R)$ is the class of the diagonal cycle, see (6.4.1). It
does not depend on the choice of compactification. 
 
\proclaim{\quad Proposition} The pairing $\epsilon_X^c =\epsilon^{\prime c}_X(-n)[-2n] : M^c (X)(-n)[-2n]\otimes M(X)\to R_\CM $ identifies $M^c (X)(-n)[-2n]$ with the dual to $M(X)$.
\endproclaim 
{\it Proof.} We want to check that $\epsilon^c$ is non-degenerate, i.e., that   the morphism $M^c \to M(X)^* (n)[2n]$ in $\CD_\CM^{\text{tri}}$ it defines (see 6.7.2)  is an isomorphism. Let us compute $\epsilon^c_X$ using 
 a smooth  compactification $X\subset   \bar{ X}$ by a divisor $D$ with normal crossings such that all the irreducible components $D_i$, $i=1,\ldots ,k$, are smooth. 

Let $F$ be a complex of presheaves with transfers with $F^{>0}=0$, $F^0 =R_{\tr}[\bar{X}]$, $F^{-a}= 
\mathop\oplus\limits_{i_1 <\ldots <i_a} R_{\tr}[D_{i_1 \ldots i_a }]$ for $a>0$; the differential is the evident one. Here $D_{i_1 \ldots i_a } := D_{i_1}\cap \ldots\cap D_{i_a}$.   The  morphism $F\to R_{\tr}[\bar{X}]/R_{\tr}[D]$ is a resolution of $R_{\tr}[\bar{X}]/R_{\tr}[D]$.

Let $G$ be a complex of presheaves with transfers with $G^{<0}=0$, $G^0 = R_{\tr}[\bar{X}]$, $G^a := \mathop\oplus\limits_{i_1 <\ldots <i_a} R_{\tr}[\bar{X}]/R_{\tr} [\bar{X}\smallsetminus D_{i_1 \ldots i_a} ]$ for $a>0$; the differential is the evident one. The  morphism $R_{\tr}[X]\to G$ is  a resolution of  $R_{\tr}[X]$ in the category of Nisnevich sheaves.

The evident pairing 
$R_{\tr}[\bar{X}]^{\otimes 2}\to     R_{\tr}[\bar{X}\times\bar{X}]/ R_{\tr}[(\bar{X}\times\bar{X})\smallsetminus \Delta (\bar{X})]$ 
extends naturally to  $\tilde{\epsilon}_{FG} : F\otimes G \to R_{\tr}[\bar{X}\times\bar{X}]/ R_{\tr}[(\bar{X}\times\bar{X})\smallsetminus \Delta (\bar{X})]$ whose other non-zero components are   $\tilde{\epsilon}_{i_1 \ldots i_a }: R_{\tr}[D_{i_1 \ldots i_a} ]\otimes R_{\tr} [\bar{X}]/R_{\tr}[\bar{X}\smallsetminus D_{i_1 \ldots i_a} ]\to   R_{\tr}[\bar{X}\times\bar{X}]/ R_{\tr}[(\bar{X}\times\bar{X})\smallsetminus \Delta (\bar{X})]  $. 
It also yields a morphism  $\tilde{\epsilon}: (R_{\tr}[\bar{X}]/R_{\tr}[D])\otimes R_{\tr}[X]\to   R_{\tr}[\bar{X}\times\bar{X}]/ R_{\tr}[(\bar{X}\times\bar{X})\smallsetminus \Delta (\bar{X})]  $, so $\tilde{\epsilon}_{FG}$ is a pairing of the resolutions that lifts $\tilde{\epsilon}$.

Let us pass to $\CD_\CM^{\text{eff\,tri}}$. Composing $\tilde{\epsilon}_{FG}$, $\tilde{\epsilon}$ with $cl(\Delta )$, we get  $R_\CM (n)[2n]$-valued pairings  $\epsilon_{FG}$, $\epsilon$ identified by the resolution quasi-isomorphisms. One has $\epsilon =\epsilon^{\prime c}_X$, so it remains to show that $\epsilon_{FG}$ is non-degenerate. It is compatible with the stupid filtrations on $F$, $G$, so it suffices to check that each $\gr^a \epsilon_{FG}$ is non-degenerate. By 6.7.1, this is true for $a=0$, since $\gr^0 \epsilon_{FG}$ equals $\epsilon'_{\bar{X}}$ from (6.7.1). For $a>0$, $\gr^a \epsilon_{FG}$ is the direct sum of pairings $\epsilon_{i_1 \ldots i_a}=cl (\Delta ) \tilde{\epsilon}_{i_1 \ldots i_a} : M(D_{i_1 \ldots i_a} )\otimes M(\bar{X}, \bar{X}\smallsetminus D_{i_1 \ldots i_a} )\to R_\CM (n)[2n]$. The non-degeneracy comes again from 6.7.1, since
the Gysin isomorphism $M(\bar{X}, \bar{X}\smallsetminus D_{i_1 \ldots i_a} )\iso  M ( D_{i_1 \ldots i_a} )(a)[2a]$ identifies $\epsilon_{i_1 \ldots i_a}$ with $\epsilon'_{D_{i_1 \ldots i_a} }$ from (6.7.1), as follows from the next exercise:

{\it Exercise.} Let $Z\hra Y\hra X$ be closed embeddings of smooth varieties of codimensions  $n$ and $m$. Then the compositions $M (X,X\smallsetminus Y) \iso M(Y)(n)[2n]\to M(Y,Y\smallsetminus Z)(n)[2n]\iso M(Z)(m+n)[2(m+n)]$ and $M (X,X\smallsetminus Y)\to M(X,X\smallsetminus Z)\iso M(Z)(m+n)[2(m+n)]$, where $\iso$ are Gysin maps (6.3.1),  are homotopic. \hfill$\square$

 \bigskip

{\bf Notation.} $\CA^{\text{op}}$ 1.1; $\CA^\kappa$ 1.1, 1.5.1; $\CA^{\text{perf}}$ 1.4.2, 1.7; 
$\CA^{\text{pretr}}$ 1.5.4; $\CA^{\text{tri}}$ 1.5.4; $ \LA$ 1.6.1; $\CA$-mod 1.1; $\CA$-dgm 1.6.3; $\CA /\CI$, $\CA /^\kappa \, \CI$ 1.4.2, 1.8; $\CA r$ 2.4; $cl$ 6.4;  $C\CA$  1.5.1; $C^b \CA$ 1.5.4;  $\CC^\CT$ 1.11;   $\CC^\CM  $ 2.3; $\CC^\Delta$ 3.1; $\CC^{\text{Nis}}$ 4.2.2; $\CC^{\text{Zar}}$ 4.1; $\CC\CH_\CM$ 6.7.3;  $\CC or$, $\CC or (\CS m)$ 2.1.2; $\CC or^\Delta (X,Y)$, $\CC or^\Delta (\CS m)$ 3.1.2; $\CC or (\CP a)$, $\CC or (\CP a^{\text{sm}})$, $\CC or ((X',U'), (X,U))$ 3.4.1; $\CC or (\CP a^{\text{nice}})$ 5.2; $\CC ous (F)$ 4.6.3; $deg $ 2.1.2;     $\CD_\CM$, $\CD_\CM^{\text{tri}}$, $\CD_{p\CM}$ 6.1; $\CD_\CM^{\text{eff}}$ 2.3; $\CD_{p\CM}^{\text{eff}}$ 3.1.2;
 $F_{-n}$ 2.2; $F_{\CT }$ 1.11; $F_{\text{Zar}}$ 4.1.1; $F_{\text{Nis}}$ 4.2.1;  $F^\cdot_{\text{Nis}}$  4.2.2; $H^{m,n}_\CM (X,R  )$ 6.2.1; $hocolim$ 1.4.1; $\CH om$ 1.9; $\Ho \,\CA$ 1.5.1;  $I^\bot$ 1.1; $i^{Gys}$ 6.3.1; $I^\CT$ 1.11; $I^{\text{cdh}} $ 6.9.2;  $I^{\text{cdh}}_{\tr} $ 6.9.3;
 $\CI^\Delta_{\tr}$ 2.3; $\CI^{\text{Nis}}_{\tr}$ 4.4; $I_{\tr}^{\text{Nis}}$ 4.3; $\CI^{\text{Zar}}_{\tr}$ 2.3; $K^M_\cdot$ 6.8; $M(X)$ 2.3; $M^c (X)$ 6.9.4; $o$ 2.1.2;  $\CP$, $\CP_{\tr}$ 2.1.4; $\CP\CS h$ 1.10, 2.1.1; $\CP\CS h_{\tr}$ 2.1.2; $\CS h_\CT$ 1.11,  $R$, $R$-mod 1.1; $R$-dgm 1.6.3; $R_\CM $ 2.3; $R [T]$, $R [U_\cdot ]$, $R [X/Y]$ 1.10; $R (n)$ 2.2; $R[G]$-perm 2.4; $\Z $,  $\Z [X]$ 2.1.2; $\Z (n)$ 2.2;  $\Z^\Delta  [X]$, $\Z^\Delta (n)$ 3.1.2; $R^c_{\tr}[X]$ 6.9.4; Res 4.6.2; $R\Gamma_\CM (X,R (n))$ 6.2; $\CS h^\CT$ 1.11; $\CS h^{\text{Nis}}$ 4.2.1;  $\CS h^{\text{Nis}}_{\tr}$ 4.3; $\CS m $ 2.1; $\CS m_0$ 2.4; $Sp $ 3.4.1; $\Sym^\cdot X$ 2.1.3; $\Bbb T$ 3.4.2;  $X^{(n)}$ 4.6.3; $X^{\text{Nis}}_\cdot$ 4.2.2; $(X,U)$ 3.4.1, $(X/S,U)$ 5.1.1 $\Pi^\varphi$ 6.1.2; $\pi_0$ 2.4.

  \end